\pdfoutput=1

\documentclass[preprint,3p,times]{elsarticle}

\usepackage{lineno,hyperref}
\usepackage{bm}
\usepackage{amsmath}
\usepackage{color}
\usepackage{algorithm}
\usepackage{algorithmic}
\usepackage{multirow}
\modulolinenumbers[5]

\journal{Journal of Computational Physics}

\begin{document}

\begin{frontmatter}

\title{A high-order semi-explicit discontinuous {G}alerkin solver for {3D} incompressible flow with application to {DNS} and {LES} of turbulent channel flow}

\author{Benjamin Krank}
\ead{krank@lnm.mw.tum.de}
\author{Niklas Fehn}
\ead{fehn@lnm.mw.tum.de}
\author{Wolfgang A. Wall}
\ead{wall@lnm.mw.tum.de}
\author{Martin Kronbichler\corref{correspondingauthor1}}
\cortext[correspondingauthor1]{Corresponding author at: Institute for Computational Mechanics, Technische Universit\"at M\"unchen, Boltzmannstr. 15, 85748 Garching, Germany. Tel.: +49 89 28915300; fax: +49 89 28915301}
\ead{kronbichler@lnm.mw.tum.de}
\address{Institute for Computational Mechanics, Technische Universit\"at M\"unchen,\\ Boltzmannstr. 15, 85748 Garching, Germany}

\begin{abstract}
We present an efficient discontinuous Galerkin scheme for simulation of the incompressible Navier--Stokes equations including laminar and turbulent flow. We consider a semi-explicit high-order velocity-correction method for time integration as well as nodal equal-order discretizations for velocity and pressure. The non-linear convective term is treated explicitly while a linear system is solved for the pressure Poisson equation and the viscous term. The key feature of our solver is a consistent penalty term reducing the local divergence error in order to overcome recently reported instabilities in spatially under-resolved high-Reynolds-number flows as well as small time steps. This penalty method is similar to the grad-div stabilization widely used in continuous finite elements. We further review and compare our method to several other techniques recently proposed in literature to stabilize the method for such flow configurations. The solver is specifically designed for large-scale computations through matrix-free linear solvers including efficient preconditioning strategies and tensor-product elements, which have allowed us to scale this code up to 34.4 billion degrees of freedom and 147,456 CPU cores. We validate our code and demonstrate optimal convergence rates with laminar flows present in a vortex problem and flow past a cylinder and show applicability of our solver to direct numerical simulation as well as implicit large-eddy simulation of turbulent channel flow at $Re_{\tau}=180$ as well as $590$. 
\end{abstract}


\end{frontmatter}

\section{Introduction}
The discontinuous Galerkin (DG) method has attained increasing popularity for simulation of the compressible Navier--Stokes equations due to a series of highly desired properties, which are stability in the convection-dominated regime, high-order capability using unstructured meshes, geometrical flexibility on curved boundaries as well as efficiency on massively parallel high-performance computers. This unique combination makes DG a very attractive approach for many high-Reynolds number applications, e.g., direct numerical simulation (DNS)~\cite{Hindenlang12,Bassi16}, large-eddy simulation (LES)~\cite{Beck14,Wiart15,Bassi16} as well as RANS and URANS~\cite{Landmann08,Wang15,Bassi14} of turbulent compressible flows. Applications range from internal turbomachinery flows~\cite{Wiart15b}, computation of high-lift configurations of an entire aircraft~\cite{Hartmann16} to environmental flows \cite{Giraldo02,Giraldo08}.

Extension of DG to incompressible flows has however been limited. Applications governed by the incompressible Navier--Stokes equations are frequently computed using standard compressible codes at small Mach numbers to avoid compressibility effects, see, e.g., \cite{Bassi16,Wiart14,Collis02}, coming along with significant time step restrictions \cite{Wiart14}, or artificial compressibility methods, see, e.g.,~\cite{Marek15} (LES),~\cite{Noventa16} (URANS) and~\cite{Crivellini13} (RANS). Fully incompressible numerical schemes in the context of LES have so far only been employed, to the authors' best knowledge, by~\cite{Ferrer12b} in the context of a 2D DG solver coupled with a spectral vanishing viscosity approach in the third space dimension and by~\cite{Tavelli16} using a space-time method.

An efficient time integration scheme is the key to large-scale simulations of incompressible turbulent flows. Coupled solvers applied within DG for example in \cite{Cockburn05,Cockburn09,Schoetzau02,Klein15,Rhebergen13} require the solution of a saddle point problem and include non-linear iterations within each time step; they have so far only been applied to small-scale academic examples. On the contrary, temporal splitting schemes allow for problem-tailored solution procedures regarding the respective terms contained in the Navier--Stokes equations, which renders them much more efficient in many applications. There are four main branches of splitting methods, namely pressure-correction, velocity-correction, algebraic-splitting and consistent-splitting schemes, see~\cite{Guermond06} for an overview. With respect to DG, a pressure-correction method with discontinuous velocity and continuous pressure has for example been proposed in~\cite{Botti11}. Pressure-correction schemes are however limited to second order accuracy in time~\cite{Guermond06}. Algebraic splitting schemes have for example been proposed in~\cite{Shahbazi07} for DG and in \cite{Lehrenfeld16} for hybrid DG discretizations.

In this paper, we consider the high-order velocity-correction (also termed velocity-projection or dual splitting) scheme by Karniadakis et al.~\cite{Karniadakis91}, which has been applied within DG in a series of recent publications \cite{Hesthaven07,Ferrer12,Ferrer12b,Steinmoeller13,Ferrer14} and within the related spectral multidomain penalty method (SMPM) in \cite{EscobarVargas14,Joshi16}. The scheme is inf-sup stable for equal-order interpolations of velocity and pressure~\cite{Karniadakis13} and splits each time step into three substeps: The non-linearity present in the convective term is first handled explicitly, a Poisson problem is subsequently solved for the pressure which is used to make the velocity divergence-free and the viscous term is taken into account in the third step.

A downside of many splitting schemes are limitations coming along with the splitting approach~\cite{Guermond06}. A primary concern of this article is to provide remedies to two such limitations observed with the present scheme when employed in conjunction with the DG method as reported in a series of recent papers~\cite{Ferrer12,Steinmoeller13,Ferrer14,Joshi16}. These issues are not related to the well-understood aliasing errors induced by under-integration of non-linear terms \cite{Hesthaven07} or convection-dominated flow regimes. The two sources of instabilities are:
\begin{itemize}
\item Ferrer and Willden~\cite{Ferrer12} and Ferrer et al.~\cite{Ferrer14} discuss instabilities encountered for small time steps with this scheme, both for continuous and discontinuous Galerkin discretizations. We show in this paper that these instabilities arise due to spurious divergence errors as a consequence of the finite spatial resolution. Several remedies to this issue are reviewed, among them a consistent div-div penalty term within the local projection reducing the point-wise divergence error. This term may be seen as a much simpler variant of the post-processing proposed in \cite{Steinmoeller13} and is similar to the grad-div term frequently used in continuous Galerkin~\cite{Olshanskii09}.
\item Violation of the mass balance through velocity discontinuities across element boundaries triggers another instability recently described by Joshi et al.~\cite{Joshi16}. This instability becomes especially relevant in spatially under-resolved simulations such as large-eddy simulation of turbulent flow. We review and benchmark two remedies to this problem, which are a supplementary jump-penalty term included in the projection according to~\cite{Joshi16} on the one hand or partial integration of the right-hand side of the Poisson equation on the other hand.
\end{itemize}

Based on these enhancements, we develop an efficient high-order DG solver applicable to laminar and turbulent incompressible flow including implicit large-eddy simulation. We choose the local Lax--Friedrichs numerical flux for discretization of the convective term and the interior penalty method \cite{Arnold82} for the Poisson problem as well as the viscous term. Local conservativity of the overall method is attained by using the divergence form of the incompressible Navier--Stokes equations \cite{Shahbazi07}. The method further employs nodal Lagrangian shape functions. In this case, nodal polynomials can be more efficient than modal ones for high polynomial orders as those used, e.g., in \cite{Ferrer12,Steinmoeller13,Ferrer14}. On the one hand, face integrals involving shape values of nodal polynomials with support on the element boundary are cheaper to evaluate. On the other hand, we observed fewer iterations of the multigrid solver/smoother combination with nodal polynomials. The method is implemented based on matrix-free operator evaluation that relies on sum factorization \cite{Kronbichler12,Kormann16} in a solver called INDEXA (INcompressible Discontinuous Galerkin towards the EXA scale) that targets modern and emerging high performance computing systems while using generic finite element programming tools.

The remainder of this article is organized as follows. We start by introducing the incompressible Navier--Stokes equations as well as the time discretization scheme in Section~\ref{sec:ins}. Thereafter, the spatial discretization is presented in Section~\ref{sec:spatial}. Section~\ref{sec:varrhs} is devoted to the instabilities in the small-time-step as well as spatially under-resolved limit and presents four variants that are benchmarked and compared providing numerical evidence for our hypothesis. Section~\ref{sec:impl} presents our implementation in the framework of tensor-product elements and matrix-free linear solvers. We verify our solver with two laminar examples in Section~\ref{sec:ex_ver} where we also compare the variants numerically. In Section~\ref{sec:dles}, the solver is applied to compute direct numerical simulation and implicit large-eddy simulation of turbulent channel flow at two Reynolds numbers. Conclusions close the article in Section~\ref{sec:conclusion}.

\section{Incompressible Navier--Stokes equations}
\label{sec:ins}
Our solver builds upon a discretization of the incompressible Navier--Stokes equations. The temporal splitting scheme outlined in Subsection~\ref{sec:timesplitting} is the starting point for discussion of the spatial discretization in Sections~\ref{sec:spatial} and~\ref{sec:varrhs}.
\subsection{Problem statement}
The incompressible Navier--Stokes equations are given in conservative form as
\begin{equation}
\frac{\partial \bm{u}}{\partial t} + \nabla \cdot \big(\bm{\mathcal{F}}^c(\bm{u}) + p \bm{I} - \bm{\mathcal{F}}^{\nu}(\bm{u})\big) = \bm{f} \text{ \hspace{0.5cm} in } \Omega \times [0,\mathcal{T}]
\label{eq:mom}
\end{equation}
with the incompressibility constraint
\begin{equation}
\nabla \cdot \bm{u} = 0  \text{ \hspace{0.5cm} in } \Omega \times [0,\mathcal{T}],
\label{eq:conti}
\end{equation}
where $\bm{u}=(u_1,u_2,u_3)^T$ is the velocity, $p$ the kinematic pressure, $\mathcal{T}$ the simulation time, $\Omega$ the domain size and  $\bm{f}=(f_1,f_2,f_3)^T$ the right-hand-side vector. The convective flux is defined as $\bm{\mathcal{F}}^c(\bm{u})=\bm{u}\otimes\bm{u}$ and the viscous flux as $\bm{\mathcal{F}}^{\nu}(\bm{u})=2 \nu \bm{\epsilon} (\bm{u})$ with the symmetric rate-of-deformation tensor $\bm{\epsilon}(\bm{u})=1/2(\nabla \bm{u}+(\nabla \bm{u})^T)$ as well as the kinematic viscosity $\nu$. At $t=0$ a divergence-free velocity field is imposed as initial condition with
\begin{equation}
\bm{u}(t=0)=\bm{u}_0   \text{ \hspace{0.5cm} in } \Omega.
\end{equation}
Boundary conditions on the Dirichlet and Neumann boundaries $\partial \Omega^D$ and $\partial \Omega^N$, with $\partial \Omega^D \cup \partial \Omega^N = \partial \Omega$ and $\partial \Omega^D \cap \partial \Omega^N = \emptyset$, are defined as
\begin{equation}
\bm{u}=\bm{g}_{\bm{u}} \text{ \hspace{0.5cm} on } \partial \Omega^D
\end{equation}
and
\begin{equation}
(-p \bm{I}+\bm{\mathcal{F}}^{\nu}(\bm{u})) \cdot \bm{n}=\bm{h} \text{ \hspace{0.5cm} on } \partial \Omega^N,
\end{equation}
where the outward unit normal vector with respect to $\partial \Omega$ is denoted by $\bm{n}$.

\subsection{Temporal velocity-correction scheme}
\label{sec:timesplitting}
For integration in time, we employ a semi-explicit multi-time-step scheme proposed by Karniadakis et al. \cite{Karniadakis91}. Herein, the transient term is discretized by a  backward-differencing formula (BDF), the non-linear convective term is treated with an extrapolation scheme (EX) and the pressure, viscous and body force terms are handled implicitly. The time-discretized momentum equation~\eqref{eq:mom} becomes
\begin{equation}
\frac{\gamma_0 \bm{u}^{n+1}-\sum_{i=0}^{J-1}\big(\alpha_i \bm{u}^{n-i}\big)}{\Delta t} + \sum_{i=0}^{J-1}\beta_i \nabla \cdot \bm{\mathcal{F}}^c(\bm{u}^{n-i})+\nabla p^{n+1}-\nabla \cdot \bm{\mathcal{F}}^{\nu}(\bm{u}^{n+1})=\bm{f}^{n+1},
\label{eq:semidis}
\end{equation}
with the solution $\bm{u}^{n+1}$ and $p^{n+1}$ at time level $t^{n+1}=(n+1)\Delta t$ with $n$ indicating the time step and $\Delta t$ the increment in time. The time integrator constants~$\gamma_0, \alpha_i$ and~$\beta_i$ of the BDF and EX schemes are given in~\cite{Karniadakis91} of which we consider the temporal orders of accuracy~$J=\{1,2,3\}$. As the scheme is not self-starting for $J=\{2,3\}$, the first time steps are performed either by successively increasing the BDF order or by interpolation of the solution to the discrete time instants $t^{n-1},...,t^{n-J+1}$ if an analytical solution is available.

Due to reasons of efficiency, we use the stiffly stable velocity-correction scheme (dual-splitting scheme) by Karniadakis et al. \cite{Karniadakis91}. In the framework of discontinuous Galerkin methods, similar approaches have been investigated by Hesthaven and Warburton \cite{Hesthaven07}, Ferrer and Willden \cite{Ferrer11}, Steinmoeller et al. \cite{Steinmoeller13} as well as Ferrer et al. \cite{Ferrer14}. It splits Equation~\eqref{eq:semidis} into three substeps: (i.) the non-linear convective term is advanced in time explicitly, (ii.) the pressure is computed by solving a pressure Poisson equation and the result is used to project the velocity onto a solenoidal space and (iii.) the viscous term is handled implicitly. The substeps are:

\renewcommand{\theenumi}{\roman{enumi}}
\begin{enumerate}
\item {\it Explicit convective step}

In the first substep, the non-linear convective term is handled efficiently by explicit time advancement
\begin{equation}
\frac{\gamma_0 \hat{\bm{u}}-\sum_{i=0}^{J-1}\big(\alpha_i \bm{u}^{n-i}\big)}{\Delta t} = - \sum_{i=0}^{J-1}\beta_i \nabla \cdot \bm{\mathcal{F}}^c(\bm{u}^{n-i})+\bm{f}^{n+1}
\label{eq:convstep}
\end{equation}
yielding the first intermediate velocity $\hat{\bm{u}}$.

\item {\it Pressure Poisson equation and projection}

The pressure step consists of solving a Poisson equation for the pressure at time $t^{n+1}$ given as
\begin{equation}
-\nabla^2 p^{n+1} = -\frac{\gamma_0}{\Delta t} \nabla \cdot \hat{\bm{u}}.
\label{eq:poisson}
\end{equation}
Consistent boundary conditions for this problem on $\partial \Omega^D$ may be derived according to \cite{Orszag86,Karniadakis91} by multiplication of the momentum equation~\eqref{eq:mom} with the normal vector. The resulting transient term is treated using the Dirichlet values given, likewise the body force, while the convective and viscous contributions are handled explicitly to avoid dependency on the velocity solution, yielding
\begin{equation}
\nabla p^{n+1} \cdot \bm{n}= -\bigg(\frac{\partial \bm{g}_{\bm{u}}(t^{n+1})}{\partial t} +\sum_{i=0}^{J_p-1}{\beta_i\big(\nabla \cdot \bm{\mathcal{F}}^c(\bm{u}_h^{n-i}) + \nu\nabla \times (\nabla \times \bm{u}^{n-i})\big)-\bm{f}^{n+1}}\bigg)\cdot \bm{n} \text{ \hspace{0.5cm} on } \partial \Omega^D .
\label{eq:bc_d_pres}
\end{equation}
Herein,~$J_p$ denotes the extrapolation order of the convective and viscous terms in the Neumann pressure boundary condition and we restrict ourselves to the case~$J_p=J$ in order to obtain optimal temporal convergence rates both in velocity and pressure. On the contrary, the mixed-order case $J_p=J-1$ discussed for example in~\cite{Guermond06} yields sub-optimal convergence rates for the pressure~\cite{Guermond06}. Note that only the solenoidal part in form of the rotational formulation of the viscous term is accounted for, which has been reported to be essential for reducing boundary divergence errors as well as high-order temporal accuracy of the overall methodology \cite{Orszag86,leriche2000high}. A suitable boundary condition on $\partial \Omega^N$ for Neumann outflow is $p^{n+1} = g_p(t^{n+1})$, prescribing the desired pressure value $g_p(t^{n+1})$ directly.

Utilizing the new pressure field $p^{n+1}$, the first intermediate velocity $\hat{\bm{u}}$ is projected onto the space of divergence-free vectors by
\begin{equation}
\hat{\hat{\bm{u}}}=\hat{\bm{u}}-  \frac{\Delta t}{ \gamma_0} \nabla p^{n+1},
\end{equation}
resulting in the second intermediate velocity $\hat{\hat{\bm{u}}}$.
\item {\it Implicit viscous step}

The final solution $\bm{u}^{n+1}$ at time $t^{n+1}$ is computed implicitly due to stability considerations by a Helmholtz-like equation reading
\begin{equation}
\frac{\gamma_0}{\Delta t}(\bm{u}^{n+1}-\hat{\hat{\bm{u}}}) = \nabla \cdot \bm{\mathcal{F}}^{\nu}(\bm{u}^{n+1}).
\label{eq:visc}
\end{equation}
\end{enumerate}

The system is closed by specifying boundary conditions for the velocity according to
\begin{equation}
\begin{array}{ll}
\bm{u}^{n+1}=\bm{g}_{\bm{u}}(t^{n+1})&\text{ \hspace{0.5cm} on } \partial \Omega^D \text{ and}\\
\bm{\mathcal{F}}^{\nu} (\bm{u}^{n+1}) \cdot \bm{n} = \bm{h}(t^{n+1}) + g_p(t^{n+1})\bm{n} &\text{ \hspace{0.5cm} on } \partial \Omega^N.
\end{array}
\label{eq:bc_vel}
\end{equation}

\section{Spatial discretization}
We commence the discussion on the spatial discretization in the first Subsection~\ref{sec:prel} by introducing the notation used. Subsequently, the variational formulation is presented in Subsection~\ref{sec:gal}. Several alternative variants of the resulting weak form are discussed and compared with regard to stability in marginally resolved simulations as well as for small time steps in the subsequent Section~\ref{sec:varrhs}.
\label{sec:spatial}
\subsection{Preliminaries}
\label{sec:prel}
In this work we consider a tessellation of the $d$-dimensional domain $\Omega \subset \mathbb{R}^d$ into $N_e$ non-overlapping hexahedral finite elements $\Omega_h = \bigcup_{e=1}^{N_e} \Omega_e$. The subscript $(\cdot)_h$  used here indicates identification of the respective variable with a characteristic element length $h$. The exterior boundaries of $\Omega_h$ are denoted by $\partial \Omega_h$. They are partitioned into a Dirichlet and Neumann boundary $\partial \Omega_h=\partial \Omega_h^D \cup \partial \Omega_h^N$ with $\partial \Omega_h^D \cap \partial \Omega_h^N = \emptyset$. Interior boundaries $\partial \Omega_e^- \cap \partial\Omega_e^+$ between two adjacent elements $\Omega_e^-$ and $\Omega_e^+$ are named $\partial \Omega_e^{\Gamma}$. The unit normal vectors of such interior boundaries are oriented outwards of the respective element yielding $\bm{n}_{\Gamma}^-=-\bm{n}_{\Gamma}^+$, i.e., $\bm{n}_{\Gamma}^-$ is oriented from $\Omega_e^-$ to $\Omega_e^+$ as well as outwards with respect to $\Omega_e^-$ on exterior boundaries accordingly. In the element-wise notation of the weak form presented in the following subsections, we refer to the current element by the superscript $(\cdot)^-$ and to the neighboring element by $(\cdot)^+$, i.e., $\bm{n}_{\Gamma}=\bm{n}_{\Gamma}^-$.

The discontinuity of the primary variables across element interfaces may be expressed in terms of jump operators $[\cdot]$ and $\llbracket \cdot \rrbracket$ defined as $[\phi]= \phi^- - \phi^+$ and $\llbracket \phi \rrbracket = \phi^- \otimes \bm{n}_{\Gamma}^- + \phi^+ \otimes \bm{n}_{\Gamma}^+$, respectively, where the latter is given for the multiplication operator $\otimes$ applicable to scalars, vectors as well as tensors and which increases the tensor rank by one. Similarly, an averaging operator is defined as $\{\{\phi\}\}=1/2(\phi^-+\phi^+)$. We also use extensions of these definitions to the boundaries, which are $\{\{\phi\}\}^{ND}=\{\{\phi\}\}$ and $[\phi]^{ND}=[\phi]$ on $\partial \Omega_e^{\Gamma}$ and $\{\{\phi\}\}^{ND}=\phi^-$ as well as $[\phi]^{ND}=0$ on $\partial \Omega_e^N\cup \partial \Omega_e^D$.

Further notation used in the weak formulations below includes $L^2$-inner products abbreviated as usual with $(a,b)_{\Omega_e}=\int_{\Omega_e}ab d\Omega$ for scalars, $(\bm{a},\bm{b})_{\Omega_e}=\int_{\Omega_e}\bm{a}\cdot \bm{b} d\Omega$ for vectors and $(\bm{a},\bm{b})_{\Omega_e}=\int_{\Omega_e}\bm{a}: \bm{b} d\Omega$ for tensors of rank~2 (the double dot product is evaluated as $\bm{a}: \bm{b} = a_{ij}b_{ij}$ in index notation). Boundary integrals are defined accordingly.

Approximate solutions are to be found within spaces of the form
\begin{equation*}
\mathcal{V}_k^p = \{p_h \in L^2 : p_h |_{\Omega_e} \in P_k(\Omega_e), \forall e \in \Omega_h \}
\end{equation*}
for the pressure and the equivalent vector-valued version $\mathcal{V}_k^{\bm{u}}=(\mathcal{V}_k^p)^d$ for the velocity. Herein, $P_k(\Omega_e)$ denotes the space of polynomials of tensor degree up to $k$ where the polynomial order is restricted by $k\ge1$ while the high-order pressure boundary conditions may first be represented sufficiently well with $k\ge2$. The polynomials are given by a tensor product of one-dimensional Lagrange polynomials based on Legendre--Gauss--Lobatto nodes for good conditioning at arbitrary polynomial degrees \cite{Karniadakis13}. Further details on aspects related to efficient implementation of these elements are discussed in Section~\ref{sec:impl}.

\subsection{Variational formulation}
\label{sec:gal}
The variational formulation is derived for each substep of the time integration scheme by multiplying Equations~\eqref{eq:convstep} to \eqref{eq:bc_vel} with appropriate weighting functions $\bm{v}_h\in \mathcal{V}_k^{\bm{u}}$ and $q_h\in \mathcal{V}_k^{p}$, respectively, and integrating over one element volume $\Omega_e^-$. If possible without ambiguity, the superscript $(\cdot)^-$ is dropped in the following for simplicity.
\begin{enumerate}

\item {\it Explicit convective step}

The flux formulation of the convective step is derived by integration by parts and subsequent application of the divergence theorem, yielding
\begin{equation}
\bigg(\bm{v}_h,\frac{\gamma_0 \hat{\bm{u}}_h-\sum_{i=0}^{J-1}\big(\alpha_i \bm{u}_h^{n-i}\big)}{\Delta t}\bigg)_{\Omega_e} =
\sum_{i=0}^{J-1}\beta_i \bigg(\big(\nabla\bm{v}_h, \bm{\mathcal{F}}^c(\bm{u}_h^{n-i})\big)_{\Omega_e}-\big(\bm{v}_h,\bm{\mathcal{F}}^{c*}(\bm{u}_h^{n-i})\bm{n}_{\Gamma}\big)_{\partial \Omega_e}\bigg)+\big(\bm{v}_h,\bm{f}_h^{n+1}\big)_{\Omega_e}
\label{eq:convstepspat}
\end{equation}
where $\hat{\bm{u}}_h \in \mathcal{V}_k^{\bm{u}}$. The local Lax--Friedrichs numerical flux is applied as it provides a stable formulation of the convective term. We have
\begin{equation}
    \bm{\mathcal{F}}^{c*}(\bm{u}_h^{n-i})=\left\{
                \begin{array}{ll}
                  \{\{\bm{\mathcal{F}}^{c}(\bm{u}_h^{n-i})\}\}+\Lambda/2\llbracket \bm{u}_h^{n-i}\rrbracket &\text{ \hspace{0.5cm}  on }\partial \Omega_e^{\Gamma},\\
                  \bm{\mathcal{F}}^{c}(\bm{u}_h^{n-i})&\text{ \hspace{0.5cm}  on } \partial \Omega_e^{N} \text{ and}\\
                  1/2\big(\bm{\mathcal{F}}^{c}(2\bm{g}_{\bm{u}}(t^{n-i})-\bm{u}_h^{n-i})+\bm{\mathcal{F}}^{c}(\bm{u}_h^{n-i})\big)+\Lambda\big(\bm{u}_h^{n-i}-\bm{g}_{\bm{u}}(t^{n-i})\big)\otimes \bm{n}_{\Gamma} &\text{ \hspace{0.5cm}  on } \partial \Omega_e^{D}
                \end{array}
              \right.
\end{equation}

using $\Lambda = \max(\lambda^-,\lambda^+)$ according to, e.g., \cite{Shahbazi07} and the maximum eigenvalue of the respective flux Jacobian,
\begin{equation}
\begin{array}{ll}
\lambda^-=\max_j\bigg|\lambda_j\bigg(\frac{\partial \bm{\mathcal{F}}(\bm{u})\cdot \bm{n}_{\Gamma}}{\partial \bm{u}}\big|_{\bm{u}_h^{-,n-i}}\bigg)\bigg| = 2|\bm{u}_h^{-,n-i} \cdot \bm{n}_{\Gamma}| \text{ and}\\
\lambda^+=\max_j\bigg|\lambda_j\bigg(\frac{\partial \bm{\mathcal{F}}(\bm{u})\cdot \bm{n}_{\Gamma}}{\partial \bm{u}}\big|_{\bm{u}_h^{+,n-i}}\bigg)\bigg| = 2|\bm{u}_h^{+,n-i} \cdot \bm{n}_{\Gamma}|.
\end{array}
\end{equation}
The latter is defined as $\lambda^+=2|(2\bm{g}_{\bm{u}}(t^{n-i})-\bm{u}_h^{-,n-i}) \cdot \bm{n}_{\Gamma}|$ on $\partial \Omega_e^{D}$. Note that we evaluate $\Lambda$ on each quadrature point, which differs from several other studies, where it is suggested to evaluate $\Lambda$ based on mean values of $\bm{u}^+$ and $\bm{u}^-$ across the entire respective element (see, e.g., \cite{Klein13,Shahbazi07}).

The explicit treatment of the convective step restricts the time step size according to the Courant-Friedrichs-Lewy (CFL) condition to
\begin{equation}
\frac{\mathrm{CFL}}{k^2} = \frac{U \Delta t}{h}
\label{eq:cfl}
\end{equation}
with $\mathrm{CFL}=\mathcal{O}(1)$, a characteristic element length $h$ and velocity $U$.

\item {\it Pressure Poisson equation and projection}

For discretization of the pressure term, the symmetric interior penalty method by Arnold \cite{Arnold82} is considered. The weak form of the right-hand side of Equation~\ref{eq:bc_d_pres}, denoted by $a(q_h,\hat{\bm{u}}_h)$ in the following, is one of the primary objects of study in this article and will be discussed in detail in Section~\ref{sec:varrhs}. In the simplest variant, we have

\begin{equation}
(\nabla q_h,\nabla p_h^{n+1})_{\Omega_e} - (\nabla q_h,(p_h^{n+1}-\mathcal{P}^*)\bm{n}_{\Gamma} )_{\partial \Omega_e} - (q_h,\Pi^{*}\cdot \bm{n}_{\Gamma})_{\partial \Omega_e}
=\underbrace{-(q_h,\frac{\gamma_0}{\Delta t} \nabla \cdot \hat{\bm{u}}_h)_{\Omega_e}}_{=:a(q_h,\hat{\bm{u}}_h)},
\label{eq:poissonspat}
\end{equation}
where the numerical flux $\mathcal{P}^*$ is given as
\begin{equation}
    \mathcal{P}^*=\left\{
                \begin{array}{ll}
                  \{\{ p_h^{n+1}\}\} &\text{ \hspace{0.5cm}  on }\partial \Omega_e^{\Gamma},\\
                  g_p(t^{n+1}) &\text{ \hspace{0.5cm}  on } \partial \Omega_e^{N} \text{ and}\\
                  p_h^{n+1} &\text{ \hspace{0.5cm}  on } \partial \Omega_e^{D}.
                \end{array}
              \right.
\end{equation}
The interior penalty flux $\Pi^{*}$ includes a stabilization term according to
\begin{equation}
    \Pi^*=\left\{
                \begin{array}{ll}
                  \{\{\nabla p_h^{n+1}\}\}-\tau_{IP}\llbracket p_h^{n+1} \rrbracket  &\text{ \hspace{0.5cm}  on }\partial \Omega_e^{\Gamma},\\
                  \nabla p_h^{n+1}-2\tau_{IP} (p_h^{n+1}-g_p(t^{n+1})) \otimes \bm{n}_{\Gamma} &\text{ \hspace{0.5cm}  on } \partial \Omega_e^{N}\text{ and}\\
                  -\bigg(\frac{\partial \bm{g}_{\bm{u}}(t^{n+1})}{\partial t} +\sum_{i=0}^{J-1}{\beta_i\big(\nabla \cdot \bm{\mathcal{F}}^c(\bm{u}_h^{n-i})+\nu\nabla \times \bm{\omega}_h^{n-i}\big)-\bm{f}^{n+1}}\bigg) &\text{ \hspace{0.5cm}  on } \partial \Omega_e^{D},
                \end{array}
              \right.
              \label{eq:ipflux}
\end{equation}
with the interior penalty parameter for hexahedra \cite{Hillewaert13}
\begin{equation}
\tau_{IP,e}=(k+1)^2\frac{A(\partial \Omega_e^{\Gamma})/2+A(\partial \Omega_e^{N} \cup \partial \Omega_e^{D})}{V(\Omega_e)}
\label{eq:taue}
\end{equation}
including surface area $A$ and element volume $V$. As usual, we apply the maximum penalty parameter across element boundaries
\begin{equation}
\tau_{IP} = \left\{
\begin{array}{ll}
\max(\tau_{IP,e}^-,\tau_{IP,e}^+) &\text{ \hspace{0.5cm}  on } \partial \Omega_e^{\Gamma} \text{ and}\\
\tau_{IP,e}^- &\text{ \hspace{0.5cm}  on } \partial \Omega_e^{N} \cup \partial \Omega_e^{D}.
\end{array} \right.
\label{eq:tau}
\end{equation}
The vorticity $\bm{\omega}_h \in \mathcal{V}_k^{\bm{u}}$ present in the boundary condition~\eqref{eq:ipflux} is evaluated employing a local $L^2$-projection by
\begin{equation}
(\bm{v}_h,\bm{\omega}_h)_{\Omega_e}=(\bm{v}_h,\nabla \times \bm{u}_h)_{\Omega_e}
\end{equation}
in order to avoid the necessity of computing the second derivatives directly, evaluating $\nabla \times\nabla \times \bm{u}_h = \nabla \times \bm{\omega}_h$.

Using the newly computed pressure $p_h^{n+1}$, the velocity is projected onto the solenoidal space by an element-wise operation according to
\begin{equation}
(\bm{v}_h,\hat{\hat{\bm{u}}}_h)_{\Omega_e} = (\bm{v}_h,\hat{\bm{u}}_h)_{\Omega_e} \underbrace{-(\bm{v}_h,\frac{\Delta t}{\gamma_0} \nabla p_h^{n+1})_{\Omega_e}}_{=:b(\bm{v}_h,p_h^{n+1})},
\label{eq:projectionspat}
\end{equation}
with $\hat{\hat{\bm{u}}}_h \in \mathcal{V}_k^{\bm{u}}$.
This equation is studied in detail in this work and we are going to modify the right-hand-side term $b(\bm{v}_h,p_h^{n+1})$ as well as add supplementary terms to ensure stability in the small-time-step limit and under-resolved case in Section~\ref{sec:varrhs}.

\item {\it Implicit viscous step}

The weak form of the viscous step~\eqref{eq:visc} reads
\begin{equation}
\begin{split}
(\bm{v}_h,\frac{\gamma_0}{\Delta t}\bm{u}_h^{n+1})_{\Omega_e} + (\bm{\epsilon}(\bm{v}_h),\bm{\mathcal{F}}^{\nu}(\bm{u}_h^{n+1}))_{\Omega_e} - s (\bm{\mathcal{F}}^{\nu}(\bm{v}_h),(\bm{u}_h^{n+1}-\bm{\mathcal{U}}^*)\otimes\bm{n}_{\Gamma})_{\partial \Omega_e}
 - (\bm{v}_h,\bm{\mathcal{F}}^{\nu*}(\bm{u}_h^{n+1})\cdot \bm{n}_{\Gamma})_{\partial \Omega_e} = 
 (\bm{v}_h,\frac{\gamma_0}{\Delta t}\hat{\hat{\bm{u}}}_h)_{\Omega_e}.
 \end{split}
\end{equation}
The parameter $s$ is chosen either as $1$ or $-1$ corresponding to symmetric or non-symmetric interior penalty methods (see e.g. \cite{Riviere08} for an overview). In the remainder of this article, we will solely consider the symmetric version $s=1$ since it qualifies for efficient solution procedures of the linear system via a conjugate gradient solver and enables optimal convergence rates of orders $k+1$ in the $L^2$ norm~\cite{Hartmann07}. The non-symmetric variant may nevertheless be useful since it provides a stable method with relaxed requirements on the penalty parameter~\eqref{eq:tau}. The numerical flux function~$\bm{\mathcal{U}}^*$ is defined as
\begin{equation}
    \bm{\mathcal{U}}^*=\left\{
                \begin{array}{ll}
                  \{\{ \bm{u}_h^{n+1} \}\} &\text{ \hspace{0.5cm}  on }\partial \Omega_e^{\Gamma},\\
                  \bm{g}_{\bm{u}}(t^{n+1}) &\text{ \hspace{0.5cm}  on } \partial \Omega_e^{D}\text{ and}\\
            \bm{u}_h^{n+1} &\text{ \hspace{0.5cm}  on } \partial \Omega_e^{N}.      
                \end{array}
              \right.
\end{equation}
The numerical flux $\bm{\mathcal{F}}^{\nu*}(\bm{u}_h^{n+1})$ includes a penalty term as already used for the pressure Poisson equation~\eqref{eq:ipflux}
\begin{equation}
    \bm{\mathcal{F}}^{\nu*}(\bm{u}_h^{n+1})=\left\{
                \begin{array}{ll}
                  \{\{\bm{\mathcal{F}}^{\nu}(\bm{u}_h^{n+1})\}\} -\tau_{IP}\nu\llbracket \bm{u}_h^{n+1}\rrbracket &\text{ \hspace{0.5cm}  on }\partial \Omega_e^{\Gamma} ,\\
                  \bm{\mathcal{F}}^{\nu}(\bm{u}_h^{n+1})-2\tau_{IP} \nu(\bm{u}_h^{n+1}-\bm{g}_{\bm{u}}(t^{n+1})) \otimes \bm{n}_{\Gamma} &\text{ \hspace{0.5cm}  on } \partial \Omega_e^{D}\text{ and}\\
                   (\bm{h}(t^{n+1})+g_p(t^{n+1})\bm{n}_{\Gamma})\otimes \bm{n}_{\Gamma} &\text{ \hspace{0.5cm}  on } \partial \Omega_e^{N},
                \end{array}
              \right.
\end{equation}
with $\tau_{IP}$ as defined in~\eqref{eq:taue} as well as~\eqref{eq:tau}.
\end{enumerate}

\section{Four variants of the projection step and impact on stability}
\label{sec:varrhs}

We pay special attention to the small-time-step as well as the spatially under-resolved limit and associated instabilities of the ``standard'' version of the present scheme as discussed in a series of recent papers~\cite{Ferrer11,Ferrer14,Steinmoeller13,Joshi16}. Our investigations have confirmed that these aspects are of high relevance regarding an accurate and robust numerical method, especially considering under-resolved turbulent flows at high Reynolds number. In the first Subsection~\ref{sec:instabilities} we give numerical evidence that the instabilities may be traced back to two particular sources within the right-hand side of the Poisson equation~\eqref{eq:poisson}. Remedies presented in literature as well as several extensions are reviewed and compared in Subsection~\ref{sec:var}. The two most promising stabilization techniques are selected in Subsection~\ref{sec:varconclusion} which will be validated and compared thoroughly in Sections~\ref{sec:ex_ver} and \ref{sec:dles} for laminar and turbulent flow, respectively. The matrix formulation is outlined in Subsection~\ref{sec:mat}.

The preliminary numerical investigations regarding small time steps shown in this section are performed using the laminar vortex problem according to \cite{Hesthaven07} and described in detail in Section~\ref{sec:vo} with a domain size of $ [-0.5,0.5] \times [-0.5,0.5]$ and discretizations of $N_e=4^2$ elements as well as the polynomial degrees $k=\{2,3,4\}$ and a kinematic viscosity of $\nu=0.025$. The cases are labeled accordingly specifying the number of elements $N4^2$, the polynomial order $k\{2,3,4\}$ and the variant under investigation, i.e. \emph{VHW} for the ``standard'' variant of Hesthaven and Warburton~\cite{Hesthaven07} given in Section~\ref{sec:spatial} or \emph{V1} to \emph{V4} discussed in the following. We use a first-order time integration scheme with $J=J_p=1$ (BDF1 and EX1) for these developments to show the influence of the temporal discretization error, which would be negligible for $J=\{2,3\}$, while results presented in the subsequent sections will mainly employ the third-order accurate scheme.

The second numerical example in this section investigating marginal spatial resolution is turbulent channel flow at a friction Reynolds number of $Re_{\tau}=180$ employing a discretization of $8^3$ elements of degree $k=3$ and a CFL number of the order of unity for the BDF3 scheme. Further details on the configuration are given in Section~\ref{sec:iles}.

\subsection{Sources of instabilities}
\label{sec:instabilities}
We commence the discussion by introducing the two modes of instability and give an overview of remedies considered in this work.
\subsubsection{The small-time-step limit}
\label{sec:instabilities_dt}
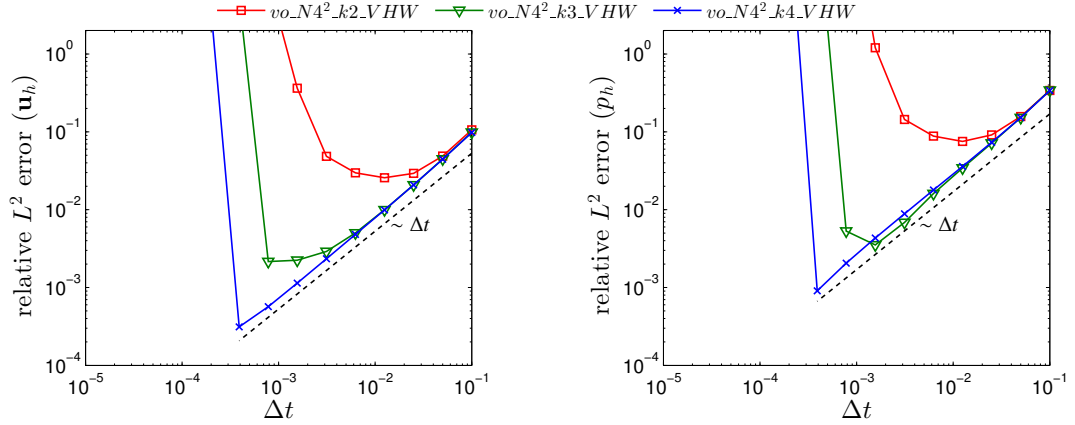
\begin{figure}[htb]
\centering
\includegraphics[trim= 0mm 0mm 0mm 0mm,clip,scale=0.6]{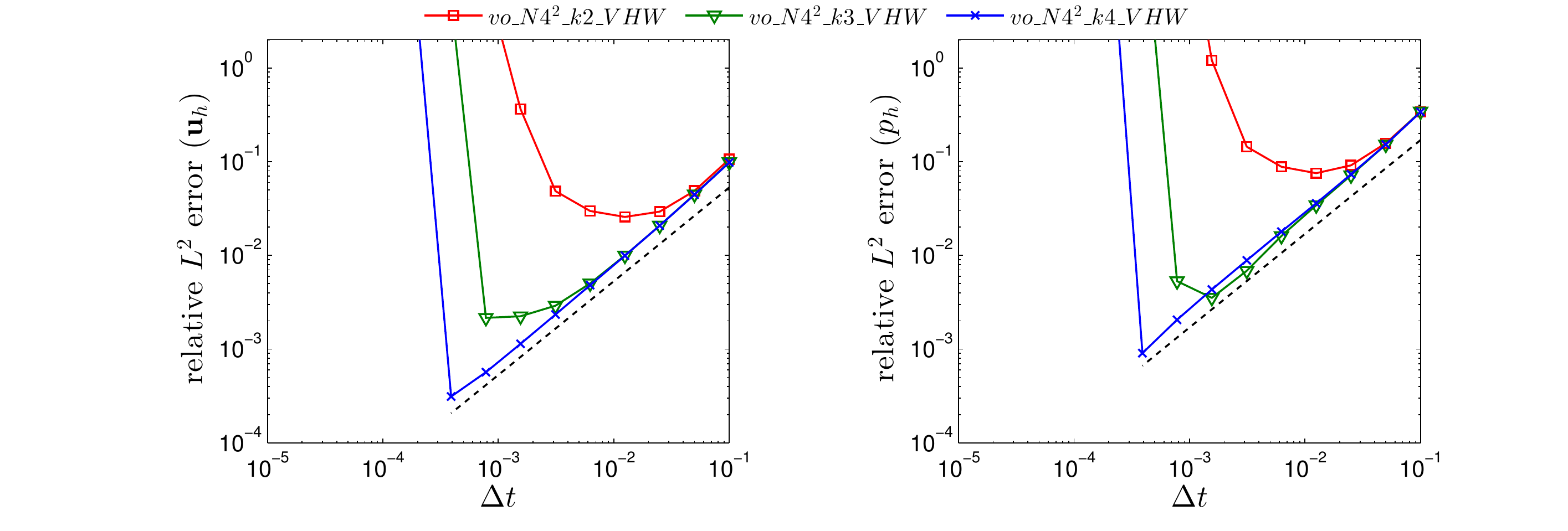}
\begin{picture}(100,0)
\put(15,90){\footnotesize $\sim \Delta t$}
\put(213,90){\footnotesize $\sim \Delta t$}
\end{picture}
\caption{Instability behavior for small time steps using the ``standard'' variant by Hesthaven and Warburton~\cite{Hesthaven07} (\emph{VHW}) for three spatial discretizations.}
\label{fig:stab_hw}
\end{figure}

The occurrence of instabilities for small time step sizes has first been described by Ferrer and Willden~\cite{Ferrer11} and has been investigated further by Ferrer et al.~\cite{Ferrer14}. The source of these instabilities may be identified by rewriting the strong form of the pressure Poisson equation~\eqref{eq:poisson} by inserting the first intermediate velocity $\hat{\bm{u}}$ of the convective step~\eqref{eq:convstep}:
\begin{equation}
-\nabla^2 p^{n+1} = -\sum_{i=0}^{J-1}\left(\frac{\alpha_i}{\Delta t} \nabla \cdot \bm{u}^{n-i}\right) -\nabla \cdot \left( - \sum_{i=0}^{J-1}\beta_i \nabla \cdot \bm{\mathcal{F}}^c(\bm{u}^{n-i})+\bm{f}^{n+1} \right).
\label{eq:PPE_instabilities}
\end{equation}
The first term on the right-hand side includes a scaling of the divergence of $\bm{u}^{n-i}$ with $1/\Delta t$, resulting in an amplification of spurious divergence errors introduced by the discretization for small time steps, which may result in an inaccurate scheme and eventually render the method unstable for very small time steps if no additional measures are taken. In order to estimate the velocity divergence errors~$\nabla\cdot\bm{u}^{n-i}$ a Helmholtz equation for the velocity divergence is considered, derived by taking the divergence of Equation~\eqref{eq:visc}
\begin{equation}
\frac{\gamma_0}{\Delta t}\nabla\cdot\bm{u}^{n+1} - 2\nu \nabla^2 (\nabla\cdot\bm{u}^{n+1})=\frac{\gamma_0}{\Delta t}\nabla\cdot\hat{\hat{\bm{u}}} .
\end{equation}
For a discrete-in-time but continuous-in-space formulation we have~$\nabla\cdot\hat{\hat{\bm{u}}}=0$ and as shown by boundary divergence analysis in~\cite{Karniadakis91,Karniadakis13} the boundary divergence errors are of order~$\left(\nabla\cdot\bm{u}^{n+1}\right)_{\partial \Omega^D}=\mathcal{O}(\Delta t^{J_p})$. For a discrete-in-space formulation, however, the intermediate velocity~$\hat{\hat{\bm{u}}}_h$ obtained in the projection step is not exactly divergence-free~\cite{Steinmoeller13}. Hence, the velocity~$\bm{u}^{n+1}$ does not only contain divergence errors related to the temporal splitting method, but also spurious divergence errors that originate from the fact that the spatial resolution is finite. Since the latter type of divergence errors does not tend to zero as~$\Delta t \rightarrow 0$, the first term on the right-hand side of Equation~\eqref{eq:PPE_instabilities} grows as~$1/\Delta t$ for small time steps which may result in instabilities. This behavior is observed in Figure~\ref{fig:stab_hw} where the laminar vortex problem has been computed using the ``standard'' method as presented in Section~\ref{sec:gal}. 

We consider the following remedies to this problem, which will be discussed in more detail in the subsequent Subsection~\ref{sec:var}:
\begin{itemize}
\item  \emph{V1}: Ferrer et al.~\cite{Ferrer14} propose to increase the penalty parameter $\tau_{IP}$ of the discrete Laplace operator of the pressure Poisson equation to circumvent the instabilities. In the same publication, it is stated that this type of instability is related to the inf-sup condition, which we cannot confirm as mixed-order elements of degrees $k$ and $k-1$ for velocity and pressure, respectively, also lead to instabilities in the limit of small time steps.
\item  \emph{V2}: The problematic term is dropped using the condition $\nabla \cdot \bm{u}^{n-i}=0$, which corresponds to another type of time integration scheme, namely the high-order non-splitting method (HONS) of Leriche and Lambrosse~\cite{leriche2000high} as well as Leriche et al.~\cite{leriche2006numerical}.
\item \emph{V3} and \emph{V4}: The divergence error is controlled by an additional and consistent div-div penalty term inspired by works of Steinmoeller et al.~\cite{Steinmoeller13} as well as Joshi et al.~\cite{Joshi16}. The term is similar to the popular grad-div term in continuous Galerkin and enables a stable numerical method.
\end{itemize}

\subsubsection{Conservation of mass in the under-resolved limit}
\label{sec:instabilities_mass}
\begin{figure}[htb]
\centering
\includegraphics[trim= 0mm 0mm 0mm 0mm,clip,scale=0.6]{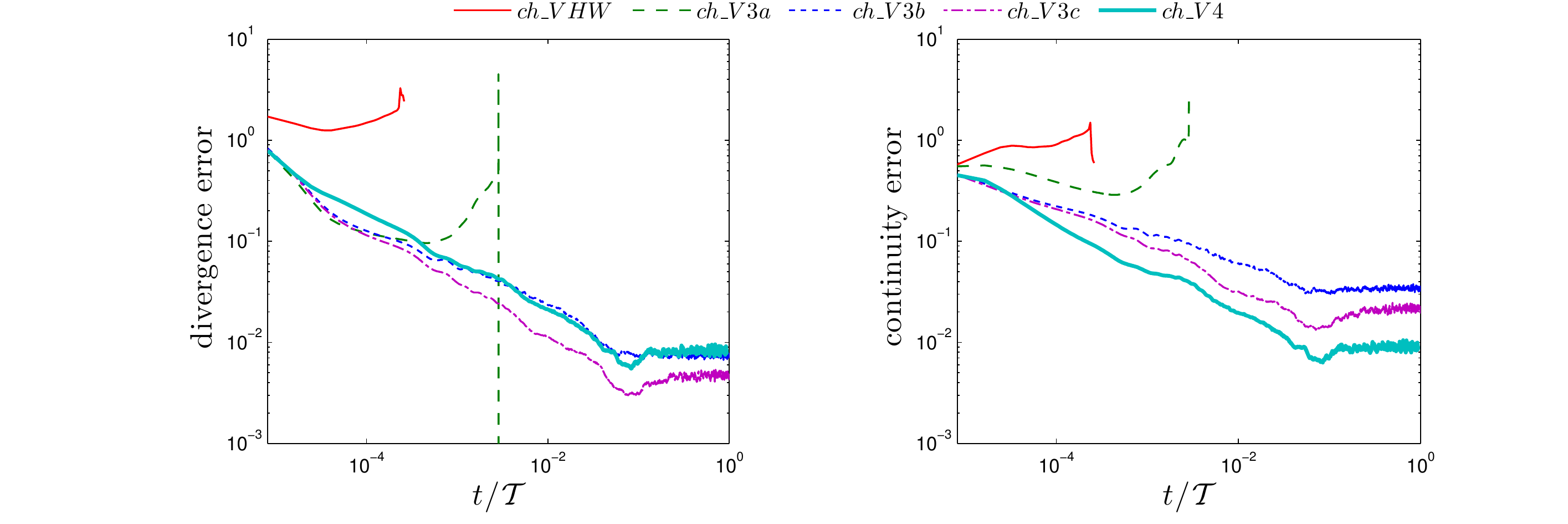}
\caption{Stability investigation for the spatially under-resolved limit for variants \emph{VHW}, \emph{V3a} to \emph{V3c} and \emph{V4} (from left to right) of turbulent channel flow at $Re_{\tau}=180$ and a spatial discretization of $8^3$ $k=3$ elements. The divergence as well as continuity errors are defined in Equations~\eqref{eq:masserror1} and \eqref{eq:masserror2}, respectively.}
\label{fig:stab_mass}
\end{figure}

The second aspect of the present discontinuous Galerkin method requiring special attention is mass conservation in the under-resolved limit. The error stemming from the continuity equation in the discontinuous context may for example be described in an element-wise sense as
\begin{equation}
e_e^{continuity}=\int_{\Omega_e}|\nabla \cdot \hat{\hat{\bm{u}}}_h| d\Omega_e + \int_{\partial \Omega_e}\frac{1}{2}|[\hat{\hat{\bm{u}}}_h] \cdot \bm{n}_{\Gamma}| d\Gamma
\label{eq:masserror}
\end{equation}
and consists of two contributions: The first term represents the divergence error within elements and the second term accounts for the mass balance across element interfaces. The factor $1/2$ is included in the latter term since the error appears on two neighboring elements. In the Poisson equation as described in Equation~\eqref{eq:poissonspat} no control on the second term is included in the DG context. As a consequence, an instability may arise in under-resolved turbulent simulations since marginal resolution generally results in more pronounced velocity discontinuities that increase the impact of the second term in Equation~\eqref{eq:masserror}. This type of instability has recently been examined by Joshi et al.~\cite{Joshi16} and is investigated numerically in Figure~\ref{fig:stab_mass} with a marginally resolved turbulent channel flow simulation. For these computations, the error in mass conservation is, according to Equation~\eqref{eq:masserror}, defined separately as the divergence error
\begin{equation}
\frac{\delta \int_{\Omega_h}|\nabla \cdot \hat{\hat{\bm{u}}}_h| d\Omega}{\int_{\Omega_h}||\hat{\hat{\bm{u}}}_h|| d\Omega}
\label{eq:masserror1}
\end{equation}
with the channel-half width $\delta$ and the continuity error 
\begin{equation}
\frac{\int_{\partial\Omega^{\Gamma}_h}|[\hat{\hat{\bm{u}}}_h] \cdot \bm{n}_{\Gamma}|
d\Gamma}{ \int_{\partial \Omega^{\Gamma}_h}|\{\{\hat{\hat{\bm{u}}}_h \}\}\cdot
\bm{n}_{\Gamma}| d\Gamma}
\label{eq:masserror2}
\end{equation}
measuring loss of mass in between elements. In Figure~\ref{fig:stab_mass}, the ``standard'' variant \emph{VHW} exhibits large divergence and continuity errors which lead to a diverging solution. 

In the present contribution, we review and compare two particular remedies to this issue, which will be discussed in more detail in the subsequent Section~\ref{sec:var}. Their common purpose of reducing the inter-element continuity error is attained by different approaches:
\begin{itemize}
\item \emph{V3b} and \emph{V3c}: We reformulate the right-hand side of the Poisson equation (term $a(q_h,\hat{\bm{u}}_h)$ in Equation~\eqref{eq:poissonspat}) such that a term both for the velocity divergence and the discontinuity at element boundaries is taken into account, in analogy to the error definition listed in Equation~\eqref{eq:masserror}. Steinmoeller et al.~\cite{Steinmoeller13} use the equivalent strong form of this term. However, that work lacks an explanation as well as investigation of this issue. We also show that it may be beneficial to modify the right-hand side of the projection, given as $b(\bm{v}_h,p_h^{n+1})$ in Equation~\eqref{eq:projectionspat}, yielding a complementary improvement of the method.
\item \emph{V4}: Joshi et al.~\cite{Joshi16} propose to penalize velocity jumps at element interfaces within the projection step, which may be seen as a straight-forward measure to the underlying problem, however at the cost of an additional global system to solve.
\end{itemize}

{\bf Remark.} The two modes of instability have been introduced here separately despite an undoubted mutual dependence especially for coarse spatial resolutions and small viscosities as shown by Ferrer et al.~\cite{Ferrer14}. In order to obtain a robust numerical method we see it as a necessity to include a measure for both limits. For example, we cannot determine the primary instability mode under investigation in Steinmoeller et al.~\cite{Steinmoeller13} since they operate in the under-resolved high-Reynolds-number regime and, in addition, do not specify the time step size used for their computations in terms of the CFL number. Indeed, their proposal is similar to \emph{V3b} discussed in the following including measures for both instabilities. 

\subsection{Four Variants}
\label{sec:var}
\subsubsection{Variant 1 (V1)}
In Ferrer et al.~\cite{Ferrer11} instabilities are reported in the limit of small time step sizes when applying the standard formulation to an unsteady Stokes flow problem. This analysis is extended in~\cite{Ferrer14} where it is shown that the time step size has to be larger than a critical time step size~$\nu\Delta t_{\mathrm{lim}} \sim h^{2}/k^{3}$ to ensure stability with the standard scheme. We explicitly note that this lower limit for the time step strongly depends on the spatial resolution characterized by~$h$ and~$k$. Consequently, a high spatial resolution has to be used in case of high-Reynolds-number flows in order to avoid conflicts of this lower bound on the time step with the CFL condition~\eqref{eq:cfl} according to~\cite{Ferrer14}
\begin{equation}
\frac{k}{h}>C_{h/p} \frac{U_{\mathrm{max}}}{\nu} ,
\label{eq:STS_and_CFL_limit}
\end{equation}
where the constant~$C_{h/p}$ is independent of the spatial resolution.

In order to circumvent the problem of instabilities in the small-time-step limit the authors in~\cite{Ferrer14} propose to stabilize the method by multiplying the interior penalty parameter~$\tau_{IP}$ of the discrete pressure Poisson operator with the amplification factor~$1/(\nu \Delta t)$ while the interior penalty parameter of the viscous step remains unchanged. A problematic aspect of this definition may be that this factor is not dimensionless and the physical unit of the penalty term no longer agrees with the unit of the other terms in the bilinear form~\eqref{eq:projectionspat}. For our numerical investigations we therefore consider the following consistent choice of the interior penalty parameter taking the~$1/\Delta t$ proportionality into account
\begin{equation}
\tau_{IP,PPE,V1}=\tau_{IP}\frac{\Delta t_{ref}}{\Delta t},
\end{equation}
where~$\Delta t_{ref}$ is a reference time step size for which the standard formulation is stable.

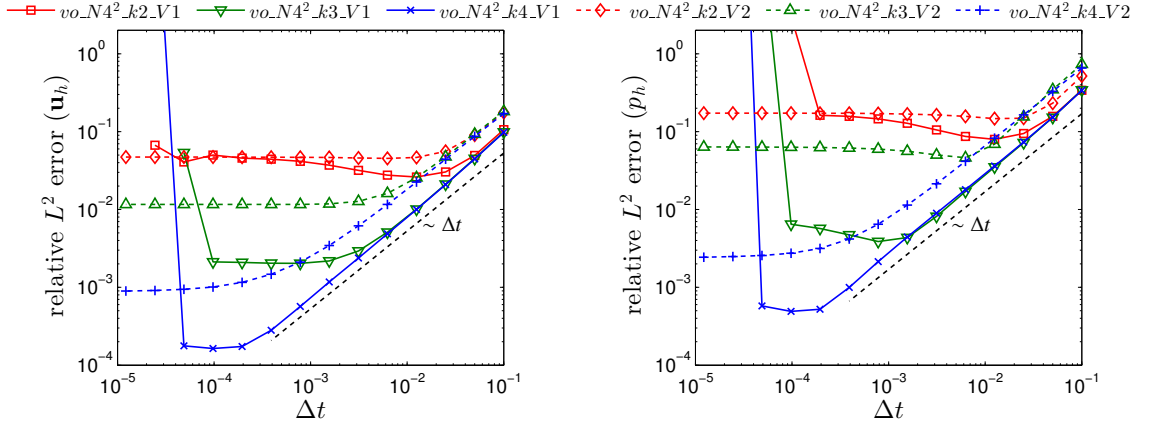
\begin{figure}[htb]
\centering
\includegraphics[trim= 0mm 0mm 0mm 0mm,clip,scale=0.6]{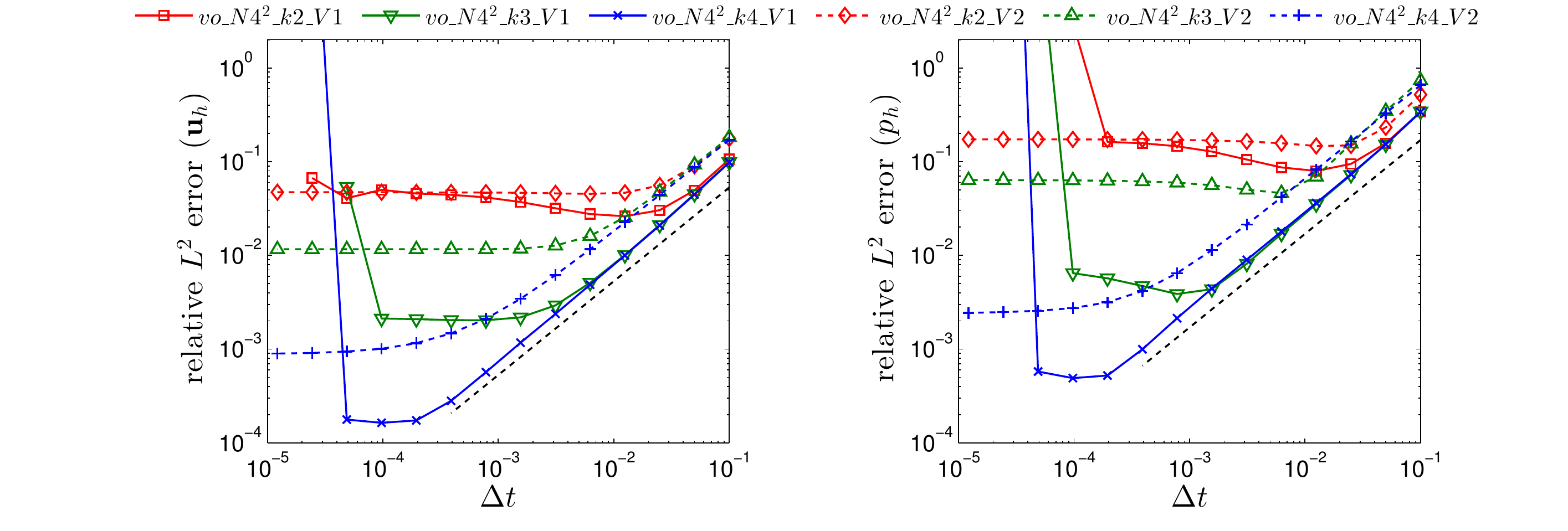}
\begin{picture}(100,0)
\put(15,90){\footnotesize $\sim \Delta t$}
\put(213,90){\footnotesize $\sim \Delta t$}
\end{picture}
\caption{Stability experiments for small time steps using \emph{V1} and \emph{V2} for three spatial discretizations.}
\label{fig:stab_v12}
\end{figure}

Our numerical experiments in Figure~\ref{fig:stab_v12} indicate that the small-time-step limit is relaxed in comparison to Figure~\ref{fig:stab_hw}. All simulations eventually get unstable, though. Moreover, we see the most problematic aspect of this stabilization approach in the drastic increase in computational cost per time step when small time step sizes are employed. Assuming that the condition number of the discrete Laplace operator is proportional to~$\tau_{IP}$~\cite{Hesthaven07}, the cost per time step increases as~$(1/\Delta t)^{1/2}$ when applying an iterative Krylov method to the numerical solution of the pressure Poisson equation (see, e.g., \cite{Shahbazi05}). Finally, this type of stabilization approach does not consider inter-element mass conservation, which would have to be included similar to \emph{V3b}, \emph{V3c} or \emph{V4}.

\subsubsection{Variant 2 (V2)}
The most straight-forward approach to tackle the problematic source term on the right-hand-side of the pressure Poisson equation given in~\eqref{eq:PPE_instabilities} as well as to verify our hypothesis regarding the small-time-step instability would be to exploit the condition $\nabla \cdot \bm{u}^{n-i}=0$ and drop the problematic term entirely. This idea leads to a time integration scheme equivalent to the method proposed in Leriche and Lambrosse~\cite{leriche2000high} and Leriche et al.~\cite{leriche2006numerical} for the unsteady Stokes equations. This decoupling approach is independent of any temporal discretization scheme and is therefore denoted as high-order non-splitting method (HONS) in~\cite{leriche2006numerical} as opposed to the high-order splitting scheme of Karniadakis et al.~\cite{Karniadakis91}. The convective step then reads
\begin{equation}
\hat{\bm{u}}_{\mathrm{d}t}=- \sum_{i=0}^{J-1}\beta_i \nabla \cdot \bm{\mathcal{F}}^c(\bm{u}^{n-i})+\bm{f}^{n+1} .
\label{eq:v2begin}
\end{equation}
Dropping the first term on the right-hand side of Equation~\eqref{eq:PPE_instabilities} yields the following pressure and projection step
\begin{align}
-\nabla^2 p^{n+1}&=-\nabla \cdot \hat{\bm{u}}_{\mathrm{d}t} ,\\
\hat{\bm{u}}&=\frac{\Delta t}{\gamma_0}\hat{\bm{u}}_{\mathrm{d}t}+\sum_{i=0}^{J-1}\frac{\alpha_i}{\gamma_0}\bm{u}^{n-i} \text{ and}\\
\hat{\hat{\bm{u}}}&=\hat{\bm{u}}-  \frac{\Delta t}{ \gamma_0} \nabla p^{n+1} .
\label{eq:v2end}
\end{align}
The viscous step is unaffected by this modification, reading
\begin{equation}
\frac{\gamma_0}{\Delta t}(\bm{u}^{n+1}-\hat{\hat{\bm{u}}}) = \nabla \cdot \bm{\mathcal{F}}^{\nu}(\bm{u}^{n+1}) .
\end{equation}

The numerical results for this scheme are depicted in Figure~\ref{fig:stab_v12} including a comparison to \emph{V1}. We observe that there is no instability with this scheme for small time steps and the error perfectly converges to a constant level where the spatial error is dominant. The results therefore support our hypothesis that the divergence term on the right-hand side of the Poisson equation~\eqref{eq:PPE_instabilities} causes instabilities in the small-time-step limit. However, Figure~\ref{fig:stab_v12} also reveals a major drawback of the modified time integration scheme: the $L^2$ error is generally larger compared to \emph{V1} for large time steps, an observation already reported by Leriche et al.~\cite{leriche2006numerical}. In further numerical experiments we did not obtain optimal convergence rates of order $k+1$ in space for velocity or pressure. If this scheme would be applied to spatially under-resolved high-Reynolds-number flows, an additional stabilization of velocity discontinuities as included in \emph{V3b, V3c} or \emph{V4} would be necessary nonetheless.

\subsubsection{Variant 3 (V3)}
\paragraph{Variant 3a}
Steinmoeller et al.~\cite{Steinmoeller13} propose to post-process the second intermediate velocity $\hat{\hat{\bm{u}}}$ to a point-wise exactly divergence-free velocity field as a means to stabilize the method. Despite the introduction of this idea as an enhancement of the splitting scheme for coarse resolutions and high Reynolds numbers in~\cite{Steinmoeller13}, we demonstrate that this approach also stabilizes in the small-time-step limit.

The post-processing step applied in~\cite{Steinmoeller13} appears costly in three space dimensions, however, since a total of nine different shape functions for each polynomial order of the modal space would be necessary. Instead of projecting the velocity field onto an exactly divergence-free basis we therefore perform this post-processing in an approximate and very efficient way by including a supplementary div-div penalty term in the projection step, which is similar to the frequently used grad-div stabilization in the context of continuous Galerkin (see, e.g.,~\cite{Olshanskii09}) to enhance mass conservation. A similar term is also included in the weak projection by Joshi~\cite{Joshi16} which is discussed in \emph{V4} and a coupled DG solver presented in~\cite{Hedwig15}. The local projection step becomes
\begin{equation}
(\bm{v}_h,\hat{\hat{\bm{u}}}_h)_{\Omega_e} + \underbrace{(\nabla \cdot \bm{v}_h, \tau_D\nabla \cdot \hat{\hat{\bm{u}}}_h)_{\Omega_e}}_{div-div\text{ }penalty} = (\bm{v}_h,\hat{\bm{u}}_h)_{\Omega_e}\underbrace{- (\bm{v}_h,\frac{\Delta t}{\gamma_0} \nabla p_h^{n+1})_{\Omega_e}}_{=b(\bm{v}_h,p_h^{n+1})},
\label{eq:div-div}
\end{equation}
where $\tau_D$ is a penalty parameter. It is noted that this represents a consistent modification of the projection step since it involves the continuity residual~$\nabla \cdot \hat{\hat{\bm{u}}}_h$. Further, the velocity field approaches the point-wise exactly divergence-free one with increasing~$\tau_D$ while the simultaneous degradation of the condition of the matrix system plays a minor role due to the locality of the problem.

We exploit the similarity of the present penalty term to the grad-div stabilization and define the penalty parameter according to \cite{Olshanskii09} for equal-order elements as
\begin{equation}
\tau_D=\zeta_D \Vert \overline{\bm{u}}_h^n\Vert h \Delta t,
\label{eq:tauD}
\end{equation}
where~$\Vert\overline{\bm{u}}_h^n\Vert$ is the norm of the element-wise volume-averaged velocity,~$h=V_e^{1/3}$ is a characteristic element length defined as the cube root of the respective element volume, and the proportionality parameter~$\zeta_D$ may be used to control the final divergence error. Note that the parameter includes a scaling with~$\Delta t$ since the projection equation is multiplied with the factor~$\Delta t$ that originates from the discrete time derivative operator. We further omit the viscous contribution accounted for in~\cite{Olshanskii09} as we focus on convection-dominated flows.

\begin{figure}[htb]
\centering
\includegraphics[trim= 0mm 0mm 0mm 0mm,clip,scale=0.6]{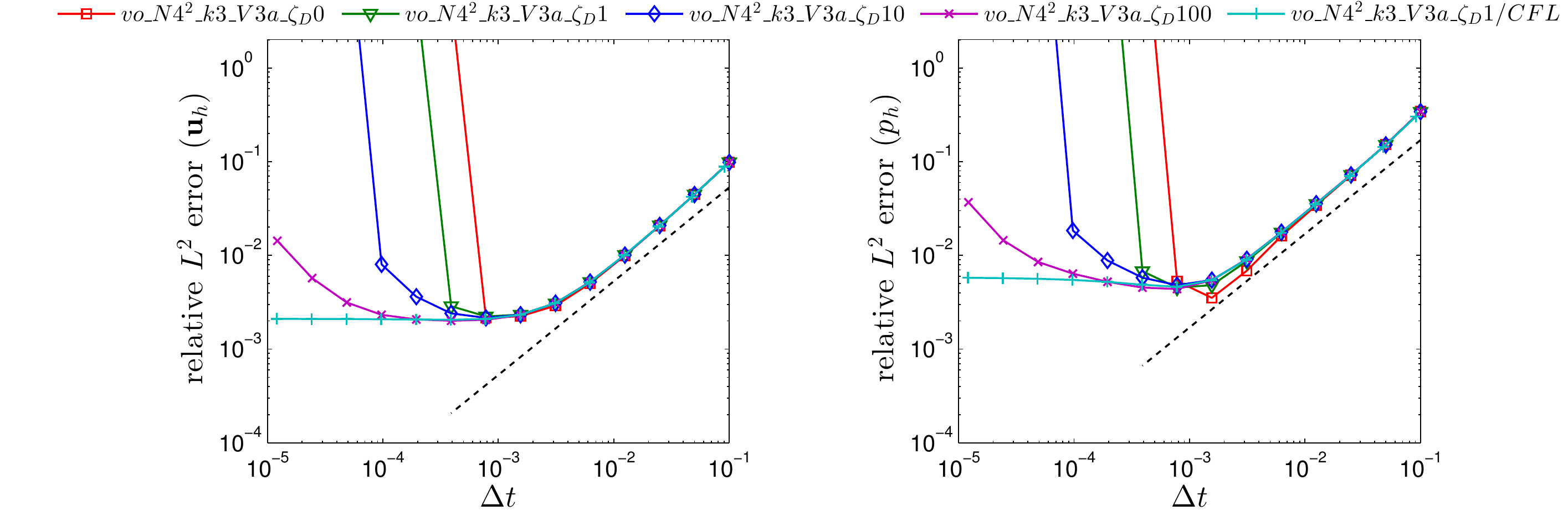}
\begin{picture}(100,0)
\put(15,90){\footnotesize $\sim \Delta t$}
\put(213,90){\footnotesize $\sim \Delta t$}
\end{picture}
\caption{Impact of proportionality parameter $\zeta_D$ included in \emph{V3} and \emph{V4} on small-time-step stability. Legend from left to right: $\zeta_D=0$, $\zeta_D=1$, $\zeta_D=10$, $\zeta_D=100$, $\zeta_D=1/\mathrm{CFL}$.}
\label{fig:stab_v3zeta}
\end{figure}
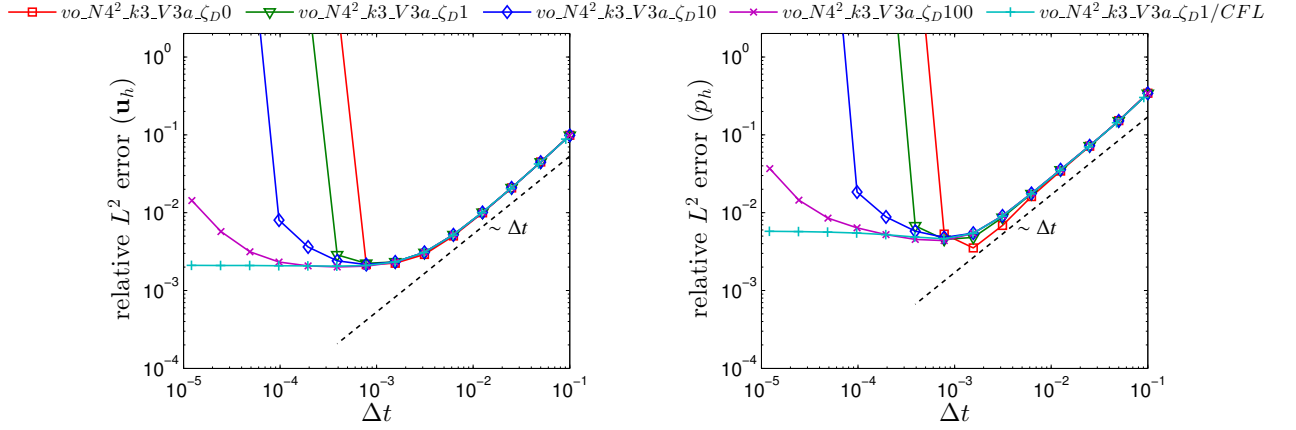

\begin{figure}[htb]
\centering
\includegraphics[trim= 0mm 0mm 0mm 0mm,clip,scale=0.6]{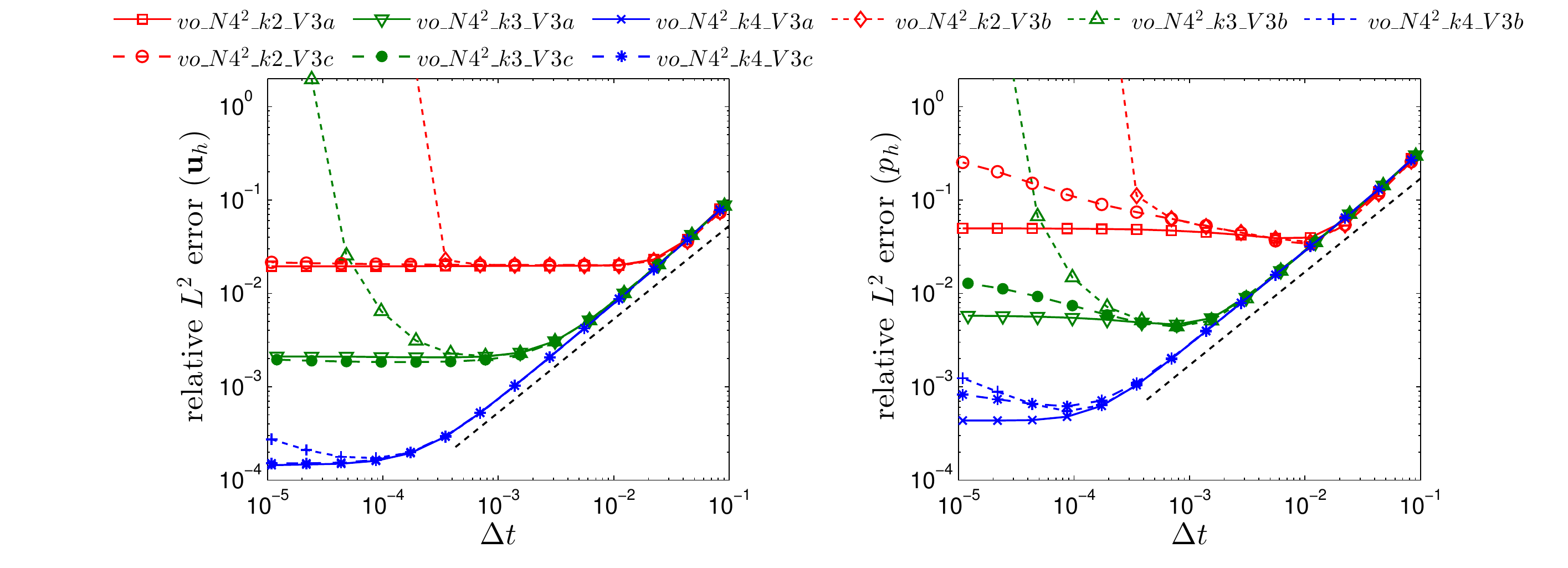}
\begin{picture}(100,0)
\put(15,90){\footnotesize $\sim \Delta t$}
\put(213,90){\footnotesize $\sim \Delta t$}
\end{picture}
\caption{Stability experiments for small time steps using \emph{V3a}, \emph{V3b} and \emph{V3c} for three spatial discretizations.}
\label{fig:sts_v3}
\end{figure}

The impact of the proportionality parameter $\zeta_D$ on the small-time-step limit is investigated in Figure~\ref{fig:stab_v3zeta} by comparing $\zeta_D = \{0,1,10,100,1/\mathrm{CFL}\}$ using the present variant \emph{V3a}. We observe a drastic improvement of the stability behavior of the splitting scheme for increasing penalty parameter and, in particular, stable results are obtained with $\zeta_D = 1/\mathrm{CFL}$ for the whole range of time step sizes considered in this example. From this result we draw the conclusion that the first term on the right-hand side of Equation~\eqref{eq:PPE_instabilities}, which appears as a source term in the pressure Poisson equation, is counter-balanced universally by a scaling of $\zeta_D$ with $1/\mathrm{CFL}$. In the remainder of this article we therefore employ $\zeta_D=\zeta_D^*/\mathrm{CFL}$ with $\zeta_D^*=1$ if not specified otherwise. Further numerical evidence for these arguments is shown in Figure~\ref{fig:sts_v3} where all cases \emph{V3a} exhibit an ideal behavior for small time steps.

\textbf{Remark.} An interesting interpretation of the grad-div stabilization again in the context of continuous Galerkin is also given in \cite{Olshanskii09} where the necessity of this term is related to an insufficient resolution of the pressure field. It is further shown that the grad-div term for continuous Galerkin may be seen as the subgrid component for the pressure. Although the transfer of this idea would certainly be illustrative in the context of the current div-div penalty term, we do not elaborate this concept within this work.

\paragraph{Variant 3b}
Since \emph{V3a} does not contain a measure controlling inter-element mass conservation according to Section~\ref{sec:instabilities_mass}, we show in \emph{V3b} a modification of \emph{V3a}  taking into account this aspect. Steinmoeller et al.~\cite{Steinmoeller13} reformulate the right-hand side of the Poisson equation~\eqref{eq:poissonspat} $a(q_h,\hat{\bm{u}}_h)$ by integration by parts and choosing a central flux formulation, yielding
\begin{equation}
a(q_h,\hat{\bm{u}}_h) = (\nabla q_h,\frac{\gamma_0}{\Delta t} \hat{\bm{u}}_h)_{\Omega_e} - (q_h,\frac{\gamma_0}{\Delta t}\{\{\hat{\bm{u}}_h\}\}^{ND}\cdot \bm{n}_{\Gamma})_{\partial \Omega_e} .
\label{eq:av3}
\end{equation}
This expression may be recast into the strong formulation by integrating by parts once again
\begin{equation}
a(q_h,\hat{\bm{u}}_h) = -(q_h,\frac{\gamma_0}{\Delta t} \nabla \cdot \hat{\bm{u}}_h)_{\Omega_e}+(q_h,\frac{\gamma_0}{\Delta t}\frac{1}{2} [\hat{\bm{u}}_h]^{ND}\cdot \bm{n}_{\Gamma})_{\partial\Omega_e} .
\label{eq:av3_strong}
\end{equation}
While the strong formulation is mathematically equivalent to Equation~\eqref{eq:av3}, it highlights that the right-hand side of the Poisson equation now is of the same structure as the continuity error defined in~\eqref{eq:masserror}, considering terms including the velocity divergence and discontinuity. The additional source term in~\eqref{eq:av3_strong} thus results in a modified pressure field that takes into account mass conservation in between elements, which comes along with a potentially less smooth pressure field in comparison with the ``standard'' variant of $a(q_h,\hat{\bm{u}}_h)$. In numerical investigations presented in Figure~\ref{fig:stab_mass} it is found that this definition of $a(q_h,\hat{\bm{u}}_h)$ indeed improves mass conservation across element boundaries, compared to the standard formulation \emph{V3a}, and results in a constant error level after an initial transient. However, it is observed in Figure~\ref{fig:sts_v3} that this modification of the Poisson equation degrades the stability for small time steps making a further modification of the projection step necessary, which is discussed in the following variant \emph{V3c}. 

\paragraph{Variant 3c}
We also consider partial integration of $b(\bm{v}_h,p_h^{n+1})$ of the projection step
\begin{equation}
b(\bm{v}_h,p_h^{n+1}) = (\nabla \cdot \bm{v}_h,\frac{\Delta t}{\gamma_0} p_h^{n+1})_{\Omega_e} - (\bm{v}_h,\frac{\Delta t}{\gamma_0} \{\{p_h^{n+1}\}\}^{ND} \bm{n}_{\Gamma})_{\partial \Omega_e}.
\label{eq:bv3}
\end{equation}
which increases robustness in our most challenging test cases and yields more accurate results. Again a central flux is used as numerical flux function. We choose the partially integrated version of $a(q_h,\hat{\bm{u}}_h)$ according to Equation~\eqref{eq:av3} and include a div-div penalty in the projection~\eqref{eq:div-div}. The results presented in Figure~\eqref{fig:sts_v3} show that this definition of $b(\bm{v}_h,p_h^{n+1})$ cures the deficiencies observed with \emph{V3b} for small time step sizes. Further, mass conservation in Figure~\ref{fig:stab_mass} exhibits even lower error levels compared to \emph{V3b} both for the divergence and continuity error.

It is noted that this combination of $a(q_h,\hat{\bm{u}}_h)$ and $b(\bm{v}_h,p_h^{n+1})$ represents a similar formulation as presented by Cockburn et al. \cite{Cockburn05} in the framework of a coupled mixed-order DG method. A supplementary pressure stabilization as proposed by Cockburn et al. \cite{Cockburn09} for equal-order coupled DG was not found to be necessary in the context of the present splitting scheme, however.

\subsubsection{Variant 4 (V4)}
\begin{figure}[htb]
\centering
\includegraphics[trim= 0mm 0mm 0mm 0mm,clip,scale=0.6]{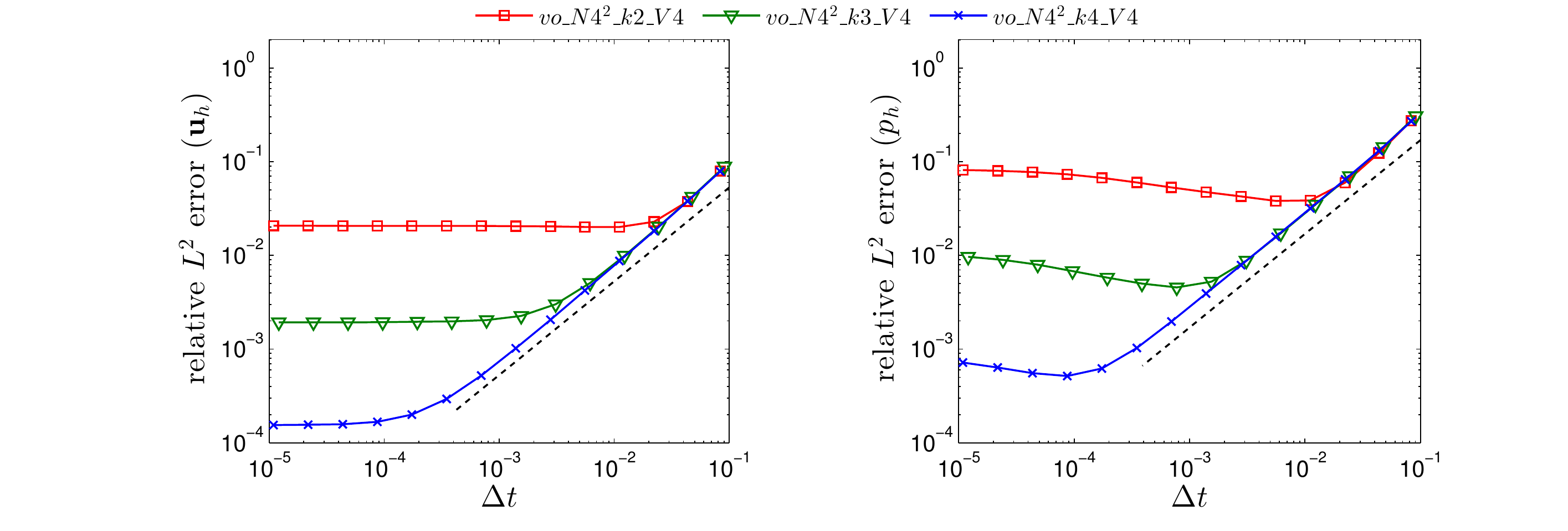}
\begin{picture}(100,0)
\put(15,90){\footnotesize $\sim \Delta t$}
\put(213,90){\footnotesize $\sim \Delta t$}
\end{picture}
\caption{Stability experiments for small time steps using \emph{V4} for three spatial discretizations.}
\label{fig:stab_v4}
\end{figure}
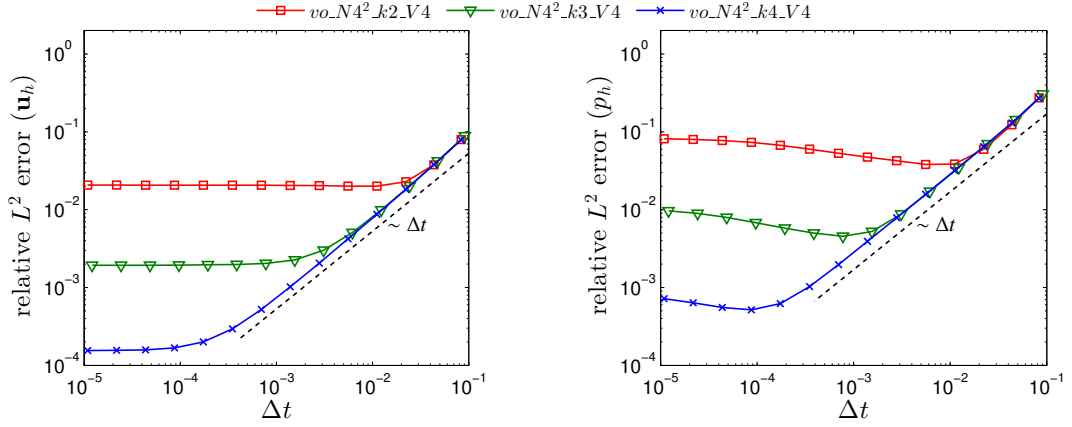

A natural approach to handle both instabilities simultaneously is to include a div-div penalty and a supplementary jump-penalty term controlling both divergence and continuity errors within the projection~\eqref{eq:projectionspat}. Joshi et al.~\cite{Joshi16} have recently proposed a similar idea where these penalty terms are contained in a post-processing step for the intermediate velocity~$\hat{\hat{\bm{u}}}_h$ in order to weakly enforce incompressibility and continuity across element faces, denoted as weak nullspace projection as opposed to the exact null-space projection by Steinmoeller et al.~\cite{Steinmoeller13} (which lead to \emph{V3a}). The projection becomes 
\begin{equation}
(\bm{v}_h,\hat{\hat{\bm{u}}}_h)_{\Omega_e} + \underbrace{(\nabla \cdot \bm{v}_h, \tau_D\nabla \cdot \hat{\hat{\bm{u}}}_h)_{\Omega_e}}_{div-div\text{ }penalty} + \underbrace{(\bm{v}_h, \tau_C[\hat{\hat{\bm{u}}}_h]^{ND})_{\partial\Omega_e}}_{jump-penalty} = (\bm{v}_h,\hat{\bm{u}}_h)_{\Omega_e}\underbrace{- (\bm{v}_h,\frac{\Delta t}{\gamma_0} \nabla p_h^{n+1})_{\Omega_e}}_{=b(\bm{v}_h,p_h^{n+1})}
\label{eq:div-div-veljump}
\end{equation}
with the standard version of $b(\bm{v}_h,p_h^{n+1})$ according to~\eqref{eq:projectionspat} and the continuity-penalty parameter $\tau_C$. We define the latter in analogy to $\tau_D$ as
\begin{equation}
\tau_{C,e}=\zeta_C \Vert\overline{\bm{u}}_h^n\Vert \Delta t,
\label{eq:tauCe}
\end{equation}
with $\zeta_C=\zeta_C^*/\mathrm{CFL}$ similar to $\zeta_D$, $\zeta_C^*=1$ if not specified otherwise. Note that this penalty parameter assures consistent physical units within the projection in contrast to the choice in~\cite{Joshi16}. On internal faces we use the average according to
\begin{equation}
\tau_{C} = \left\{
\begin{array}{ll}
\{\{\tau_{C,e}\}\} &\text{ \hspace{0.5cm}  on } \partial \Omega_e^{\Gamma} \text{ and}\\
\tau_{C,e}^- &\text{ \hspace{0.5cm}  on } \partial \Omega_e^{N} \cup \partial \Omega_e^{D}.
\end{array} \right.
\label{eq:tauC}
\end{equation}
The standard right-hand side of the Poisson equation $a(q_h,\hat{\bm{u}}_h)$ as in Equation~\eqref{eq:poissonspat} is used for this variant. Due to the jump-penalty term, \textit{V4} makes the projection step a more expensive global equation system compared to the purely local projection used in \emph{V3}. Our preliminary investigations in Figure~\ref{fig:stab_v4} indicate ideal behavior for small time steps while maintaining low error levels throughout. According to Figure~\ref{fig:stab_mass}, conservation of mass gives a similar behavior as variants \emph{V3b} and \emph{V3c} which makes this approach a promising alternative to \emph{V3c} despite the additional computational cost.

\subsection{Conclusion on Variants 1--4}
\label{sec:varconclusion}
In this section, we have identified small time steps and coarse spatial resolutions as potential sources of instabilities and discussed a number of remedies. On the one hand, spurious divergence errors are amplified for small time steps leading eventually to an unstable scheme. On the other hand, pronounced velocity discontinuities in under-resolved simulations give rise to excessive violation of the continuity equation also leading to instabilities. While the first issue may be stabilized successfully via a div-div penalty term within the projection, two remedies exhibit promising characteristics for stabilization of marginally-resolved simulations: We have obtained a robust and fast computational method by partial integration of the right-hand side of the Poisson equation and projection, detailed in \emph{V3c}, while a supplementary jump-penalty term within the projection step also yields a promising method at slightly higher computational cost according to \emph{V4}. 

We anticipate at this point that variant \emph{V3c} is our working-horse for turbulent flows, since it combines all highly desirable aspects of a numerical scheme, consisting of stability regarding small time steps and spatially under-resolved simulations, optimal spatial convergence rates as well as low computational cost, thus most of the turbulent flow examples shown below are computed with this variant.

\subsection{Matrix formulation}
\label{sec:mat}
The discussion on the spatial discretization in Sections~\ref{sec:spatial} and~\ref{sec:varrhs} is concluded with the matrix formulation, which is the basis for the presentation of the solution procedures employed for the linear systems in the subsequent Section~\ref{sec:impl}. The matrix formulation for the convective step~\eqref{eq:convstepspat} results in
\begin{equation}\label{eq:advect}
\gamma_0 \hat{\bm{U}} = \sum_{i=0}^{J-1}\alpha_i \bm{U}^{n-i}- \Delta t  \bm{M}^{-1} \sum_{i=0}^{J-1}\beta_i \bm{F}^c(\bm{U}^{n-i})+\bm{F}(t^{n+1}),
\end{equation}
with the block-diagonal mass matrix $\bm{M}$, the evaluation of the convective term for the corresponding time step $\bm{F}^c$, the body-force vector $\bm{F}$ and the respective velocity vectors $\bm{U}$. The matrix form of the pressure Poisson equation~\eqref{eq:poissonspat} is given by
\begin{equation}\label{eq:poissonmat}
\bm{L} \bm{P}^{n+1} = \frac{\gamma_0}{\Delta t}\bm{A}\hat{\bm{U}}- \bm{L}_{BC}\bm{G}_{p}(t^{n+1})
\end{equation}
with the discrete Laplace operator $\bm{L}$, the pressure solution vector $\bm{P}$, the respective variant of the velocity divergence operator $\bm{A}$ according to $a(q_h,\hat{\bm{u}}_h)$ in Section~\ref{sec:var} and boundary terms $\bm{L}_{BC}\bm{G}_{p}$. For the local projection~\eqref{eq:projectionspat}, we get
\begin{equation}\label{eq:project}
\left(\bm{M}+\tau_D \bm{D}+\tau_C \bm{C}\right) \hat{\hat{\bm{U}}} = \left( \bm{M}\hat{\bm{U}} + \frac{\Delta t}{\gamma_0} \bm{B} \bm{P}^{n+1}\right)
\end{equation}
with the block-diagonal div-div penalty operator $\bm{D}$ only considered in {\it V3} and {\it V4}, the jump-penalty terms $\bm{C}$ solely included in {\it V4} as well as the discrete pressure gradient $\bm{B} \bm{P}$ according to $b(\bm{v}_h,p_h^{n+1})$. Finally, the Helmholtz-like equation of the viscous step reads in matrix form
\begin{equation}\label{eq:viscous}
\left(\frac{\gamma_0}{\Delta t}\bm{M}-\bm{F}^{\nu}\right)\bm{U}^{n+1} = \frac{\gamma_0}{\Delta t}\bm{M} \hat{\hat{\bm{U}}} + \bm{F}^{\nu}_{BC}\bm{G}_{\bm{u}}(t^{n+1})
\end{equation}
with the linearized viscous term $\bm{F}^{\nu}$ and the right-hand-side boundary terms $\bm{F}^{\nu}_{BC}\bm{G}_{\bm{u}}$. The present matrix formulation for variants \emph{V1}, \emph{V3} and \emph{V4} extends naturally to the case \emph{V2} by modification of~\eqref{eq:advect} to \eqref{eq:project} according to \eqref{eq:v2begin} to \eqref{eq:v2end}.

\section{Implementation}
\label{sec:impl}
The solver outlined above has been implemented in a high-performance C++ code
based on the deal.II finite element library \cite{Bangerth16}. For the
evaluation of integrals on cells and faces, a highly efficient quadrature
approach based on sum factorization for hexahedra is used
\cite{Kronbichler12,Kormann16}. Sum factorization is a special evaluation
technique aimed at exploiting re-occurring terms in evaluation of unit cell
shape functions and derivatives defined through tensor products. The expansion of the elemental solution is of the form
\begin{equation}
 p_h(\bm{\xi},t)=\sum_{l,m,n=0}^k N_{lmn}^{k}(\bm{\xi})p_{lmn}(t),
\end{equation}
where the shape function $N_{lmn}^k$ is the product of three one-dimensional shape functions of degree $k$ each, $N_{lmn}^k(\bm{\xi}) = \ell_l^k(\xi_1) \ell_m^k(\xi_2) \ell_n^k(\xi_3)$. If this expansion is combined with a quadrature rule constructed by a tensor product of a 1D formula, the interpolation of nodal values on quadrature points as well as the multiplication by test functions and summation for cell- and face-wise quadrature can be realized by a series of operations along one-dimensional stripes. In two space dimensions, the tensorized evaluation exchanges evaluations of the form $(\bm A\otimes \bm B) \bm v$ (matrix-vector) by $\bm B \,\text{mat}(\bm u)\, \bm A^T$ (matrix-matrix), where $A$ and $B$ are matrices containing all 1D shape values or derivatives $\ell_{l=0:k}$ on all 1D quadrature points $x_{q=1:n_q}$ and $\text{mat}(\bm u)$ reshapes the $(k+1)^2$ long vector of nodal values into a $(k+1)\times(k+1)$ matrix. This transformation reduces the evaluation complexity from $\mathcal O((k+1)^{2d})$ operations in the naive evaluation in $d$ dimensions to $\mathcal O(d(k+1)^{d+1})$ operations \cite{Karniadakis13,Kopriva09} and has been applied in a variety of contexts, see e.g.~\cite{Kopriva09,Kronbichler12,Hindenlang12,Bastian14,MayBrown14}.

These fast integration routines are used for all matrix-vector products and weak form evaluations in the steps~\eqref{eq:advect}--\eqref{eq:viscous} listed in Section~\ref{sec:mat}. Cell and face integrals for the nonlinear convective term are evaluated by Gaussian quadrature on $\left \lfloor \frac{3k}{2}\right \rfloor +1$ points per coordinate direction which ensures exact numerical integration on affine cells. This avoids aliasing effects due to inexact quadrature. All other integrals contain terms up to polynomial degree $2k$ and are integrated with a Gauss formula on $n_q = k+1$ points per dimension. Geometries with curved boundaries are interpolated by iso-parametric polynomial mappings through a manifold description of the boundary using facilities of the deal.II library \cite{Bangerth16}.

As usual in DG, the mass matrices $\bm M$ appearing in the algorithm are block-diagonal over the elements. Since we use nodal Lagrange polynomials on Gauss--Lobatto points, the consistent mass matrix is non-diagonal. In order to avoid applying and storing a mass matrix for each cell on non-affine meshes that comes at cost $\mathcal O((k+1)^{2d})$, we use a tensorial sum factorization approach similar to the one used for computing integrals, as described in \cite{Kronbichler15}. In this algorithm, the tensor product matrices for transforming the nodal values to quadrature points and back are inverted individually for the case $n_q = k+1$. The result is a familiar tensorial kernel involving the inverse of the 1D transformation from node values to quadrature points and inverse quadrature weights. An alternative interpretation of this algorithm is the transformation into orthogonal nodal polynomials on Gauss points instead of the Gauss--Lobatto basis, applying the diagonal inverse mass matrix, and transforming back to the Lagrange basis on Gauss--Lobatto points.

\subsection{Solving linear systems}

The linear systems from Equations~\eqref{eq:poissonmat} and \eqref{eq:viscous} are solved by the preconditioned conjugate gradient (CG) method. Since we want to avoid the storage cost and performance penalty of sparse matrices, a matrix-free evaluation of the matrix operators is employed. To this end, the same fast integration kernels as for the right-hand sides and nonlinear terms are used. As demonstrated in \cite{Kronbichler12}, matrix-free kernels outperform sparse matrix-vector by a factor of 2--5 for continuous elements at polynomial degree 2 already, with larger advantages for the matrix-free kernel at higher degrees due to the $\mathcal O((k+1)^{d+1})$ complexity against $\mathcal O((k+1)^{2d})$. For DG-SIP, the cost advantage of matrix-free kernels increases by another factor of two even when assuming optimal sparse matrix storage with compressed index data for the elemental block in DG. We apply problem-tailored preconditioning strategies for the respective equation systems as follows.

\paragraph{Viscous solver} For the viscous matrix, a relatively simple preconditioner based on the inverse mass matrix proves effective, see also e.g.~\cite{Shahbazi07}. This is the case if the viscous term is not too large or time steps are sufficiently small such that the eigenvalue spectrum of the operator $\frac{\gamma_0}{\Delta t}\bm M - \bm{F}^{\nu}$ is close to the one of the mass matrix.

In order to reduce the iteration count of the viscous solver, we use extrapolations of order $J$ of the velocity values as an initial guess by employing the same constants $\beta_i$ as used for extrapolation of the convective term, enabling high absolute accuracy despite low relative tolerances of the linear solver. For the computations shown in the application section, iteration counts are between 2 and 40 iterations, usually around 15.

\paragraph{Poisson solver} For the pressure Poisson matrix, a multigrid preconditioner is selected in order to obtain iteration numbers that are independent of the mesh size $h$. In our code, we choose the preconditioner as a single V-cycle of a geometric multigrid method with a polynomial Chebyshev smoother, cf.~\cite{Adams03}, see also \cite{Sundar15} for a recent evaluation. This particular smoother only needs the matrix diagonal (which is pre-computed once before time stepping) and the action of the matrix operator, thus enabling the use of our highly efficient matrix-vector products. We use a degree of 5 for the Chebyshev polynomial for pre- and post-smoothing (i.e., five matrix-vector products per smoother step) with a smoothing range of 20 on each level. As a coarse solver, we again use the Chebyshev iteration where the number of matrix-vector steps is adapted to the eigenvalue distribution of the matrix such that the usual Chebyshev error estimator reaches a tolerance of $10^{-3}$~\cite{Varga09}. This coarse solver has the advantage that it can leverage parallelism for moderate coarse mesh sizes up to a few hundreds of elements and avoids global communication through inner products inner products. Furthermore, it is preferred over a coarse CG solver because it performs a fixed number of iterations and thus provides a linear operation, crucial for enabling an outer CG solver. Note that we coarsen down to one element if the mesh hierarchy described by the forest-of-tree concept \cite{Burstedde11} allows for that. All eigenvalue estimates for the Chebyshev parameters are done through an initial CG solution on each multigrid level. The level transfer between different grid levels is realized by matrix-free tensorial embedding operations from mother to child cells and vice versa. Since multigrid is only used for preconditioning, the whole V-cycle is done in single precision which is twice as fast as double precision apart from the MPI communication, due to wider vectorization and less memory transfer. The geometric multigrid hierarchy is related to the parallel mesh storage provided by the p4est library \cite{Burstedde11,Bangerth11} where the ownership of coarse level cells is recursively assigned from the owner of the first child on the next finer level \cite{Bangerth16}.

We also employ extrapolation of the pressure values of order $J$ to predict the pressure solution as for the viscous solver to be able to use a low relative tolerance. Iteration counts for the pressure Poisson solver are between 2 and 15 iterations, usually around 8, depending on the mesh, time step size and initial guess.

For pure Dirichlet problems including the channel flow reported below, the pressure is subject to Neumann conditions on the whole boundary and thus, the linear system is singular. To guarantee solvability, a subspace projection to the space of vectors with zero mean is applied to the right-hand side of the Equation~\eqref{eq:poissonmat} and the matrix-vector product in the CG solver. Note that the projection is not applied on the multigrid levels because the input vectors have zero mean and no non-zero mean can be created through the hierarchy. Otherwise, expensive global communication to compute and broadcast the global vector mean would be necessary, affecting scaling beyond 10,000 processor cores.

\paragraph{Projection solver} Finally, we also choose an iterative solver for the div-div projection introduced in \emph{V3}. In the form stated in~\eqref{eq:project}, the individual elements decouple and an iterative conjugate gradient solver is chosen for each element. This again leverages fast matrix-vector products through sum factorization at a complexity of at most $\mathcal O((k+1)^{2d+1})$, obtained when the conjugate gradient solver converges after as many steps as there are matrix rows, rather than $\mathcal O((k+1)^{3d})$ for matrix factorization in a direct solver. We found that CG without preconditioner only converges slowly due to round-off. This is because for large penalty parameters $\tau_D$, the eigenvalues corresponding to divergence-free velocities are by more than a factor of $10^6$ smaller than the penalized ones. However, using the inverse mass matrix for preconditioning resolves the issue: All solenoidal velocities (two thirds of eigenvalue spectrum) have exactly eigenvalue one which is accurately detected by the conjugate gradient method after a few iterations. In all of our experiments, the iteration count has been less than approximately half the local size of the div-div matrix also when $\zeta_D=10^6$. Note that in our implementation where several elements are processed at once \cite{Kronbichler12} via vector instructions (AVX on Intel processors), several systems are solved through vectorized data types with individual CG parameters each.

In case a jump-penalty term is also included in \emph{V4} according to \eqref{eq:div-div-veljump}, the global matrix is no longer decoupled over elements and a global CG solver preconditioned by the inverse mass matrix is used with similar iteration counts as for the local solver.

\begin{figure}[htb]
\centering
\includegraphics[trim= 0mm 0mm 0mm 0mm,clip,scale=0.6]{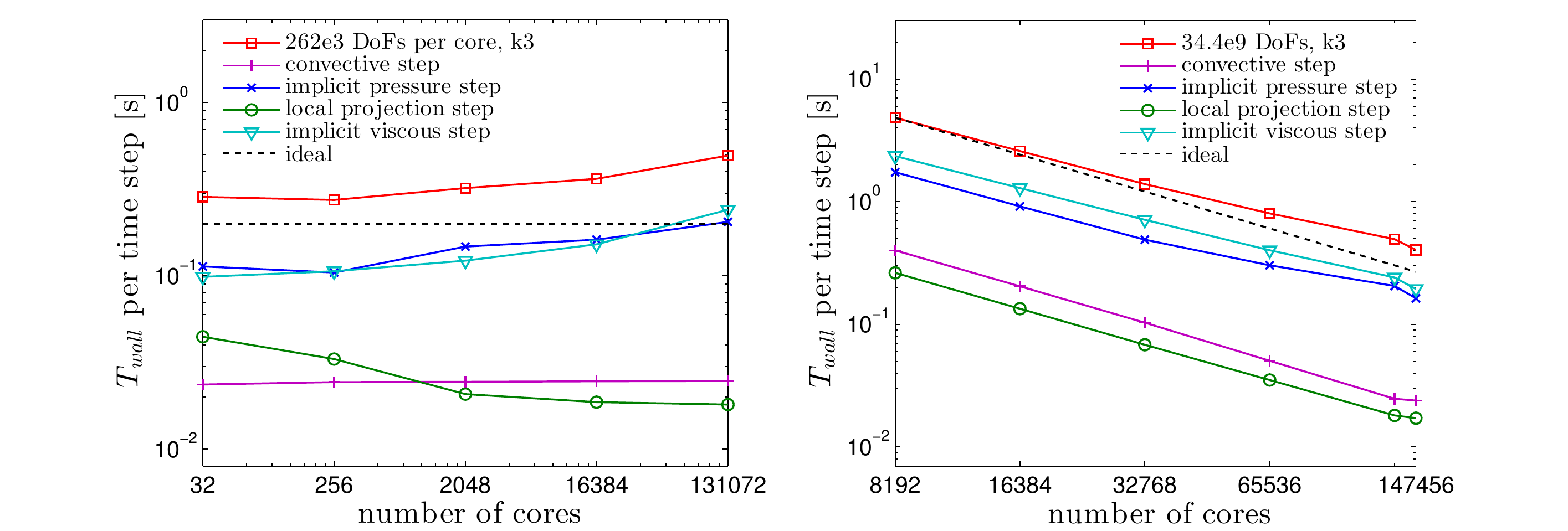}
\includegraphics[trim= 0mm 0mm 0mm 0mm,clip,scale=0.6]{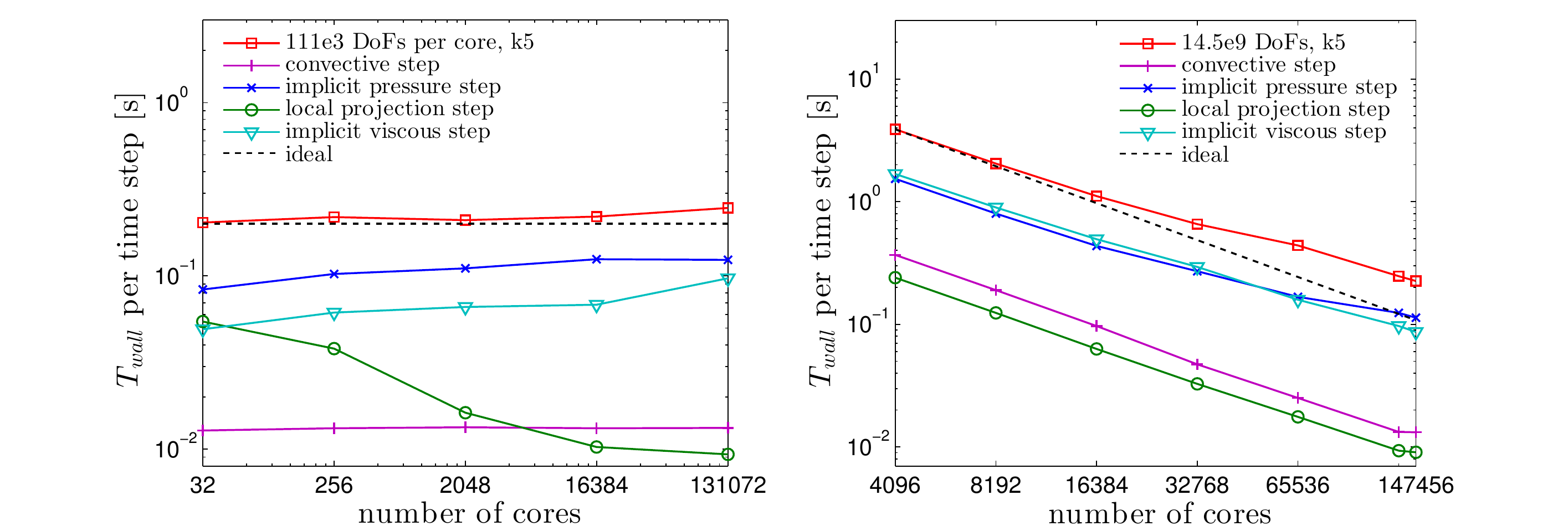}
\caption{Wall time $T_{wall}$ per time step showing weak (left) and strong (right) scaling up to 147,456 CPU cores for third (top) and fifth (bottom) polynomial degree.}
\label{fig:scaling}
\end{figure}

\subsection{Performance evaluation}

Our implementation targets large-scale parallel computations on modern supercomputers. Figure~\ref{fig:scaling} shows results of weak and strong scaling experiments for variant \emph{V3c} on up to 147,456 processor cores on the SuperMUC Phase 1 system (9,216 nodes of dual-socket, eight-core Intel Sandy-Bridge processors at 2.7 GHz each) in Garching, Germany. Given that our code already performs very well in serial due to the particular matrix-free kernels, the observed scaling of all solver components is excellent, including the challenging pressure multigrid algorithm. For example, the pressure Poisson solver processes almost 900,000 degrees of freedom per second and core also on the largest computation with 147k MPI ranks. Further, a comparison of the weak scaling plots in Figure~\ref{fig:scaling} reveals that the computation time per DoF and time step is essentially independent of the polynomial degree, making the use of high orders very attractive for example for laminar flows and direct numerical simulation of turbulent flows.

The weak scaling in Figure~\ref{fig:scaling} shows that the computing time of the local projection step actually reduces as the number of cores is increased. The reason for this somewhat unexpected behavior is that the same material parameters are selected for all problem sizes. As the processor count increases, the resolution increases and thus, the need for penalization of divergence error reduces. This makes the iterative CG solver converge more quickly, reducing the iteration count by almost a factor of 10 when going from 32 to 128k processors.

\section{Verification}
\label{sec:ex_ver}
We verify the code described above and compare the two variants \emph{V3c} and \emph{V4} which emanated as the most promising ones from our discussion in Section~\ref{sec:varrhs} in the following. We commence in \ref{sec:vo} by proving optimal convergence rates in time (Subsection~\ref{sec:vo_t}) and space (Subsection~\ref{sec:vo_h}) for velocity and pressure using the vortex problem already investigated in Section~\ref{sec:varrhs}. We also discuss the efficiency of our code for laminar flows regarding the optimal use of high polynomial orders in Subsection~\ref{sec:performance_k} in the context of this example. A further test case is presented in \ref{sec:fpc} consisting of an unsteady laminar flow past a cylinder that demonstrates the geometrical flexibility of the present approach.

While the examples in this section employ the 2D implementation of the code, we discuss application to large-scale simulations of 3D turbulent channel flow in Section~\ref{sec:dles}.

\subsection{Vortex problem}
\label{sec:vo}
We consider a laminar vortex problem according to~\cite{Hesthaven07} with the analytical solution for velocity and pressure given as
\begin{align}
\bm{u}(\bm{x},t) &= 
\begin{pmatrix}
-\sin\left(2\pi x_2\right) \\
+\sin\left(2\pi x_1\right)
\end{pmatrix}\exp\left(-4\nu \pi^2 t\right)  \text{ and}\\
p(\bm{x},t) &= -\cos\left(2\pi x_1\right)\cos\left(2\pi x_2\right)\exp\left(-8\nu\pi^2 t\right),
\end{align}
defined in the domain $[-0.5,0.5]\times[-0.5,0.5]$ with respective Dirichlet boundary conditions on the inflow and exact Neumann boundary conditions given as $g_p$ and $\bm{h}$ on the outflow boundaries. We choose a viscosity of $\nu=0.025$, the simulation time as $\mathcal{T}=1$ and define the CFL number according to Equation~\eqref{eq:cfl} with the maximum velocity $U=1.4$. The relative $L^2$ error is computed at $t=\mathcal{T}$ and is defined as
\begin{equation}
\frac{{\Vert\ \bm{u}(\bm{x},t=\mathcal{T})-\bm{u}_h(\bm{x},t=\mathcal{T})\Vert}_{\Omega_h}}{{\Vert\bm{u}(\bm{x},t=\mathcal{T})\Vert}_{\Omega_h}}
\end{equation}
for the velocity and
\begin{equation}
\frac{{\Vert\ p(\bm{x},t=\mathcal{T})-p_h(\bm{x},t=\mathcal{T})\Vert}_{\Omega_h}}{{\Vert p(\bm{x},t=\mathcal{T})\Vert}_{\Omega_h}}
\end{equation}
for the pressure.

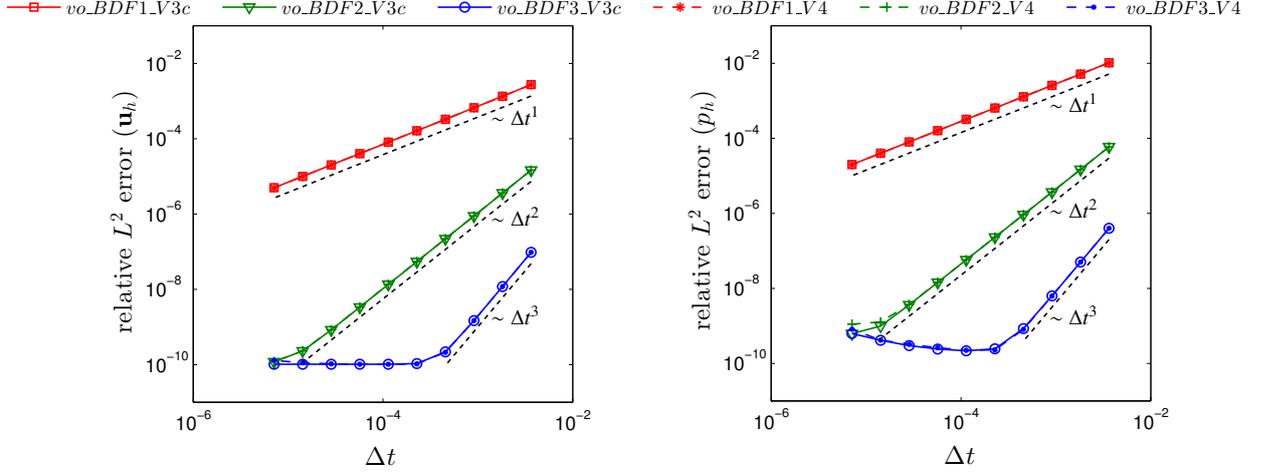
\begin{figure}[t!]
\centering
\includegraphics[trim= 0mm 0mm 0mm 0mm,clip,scale=0.6]{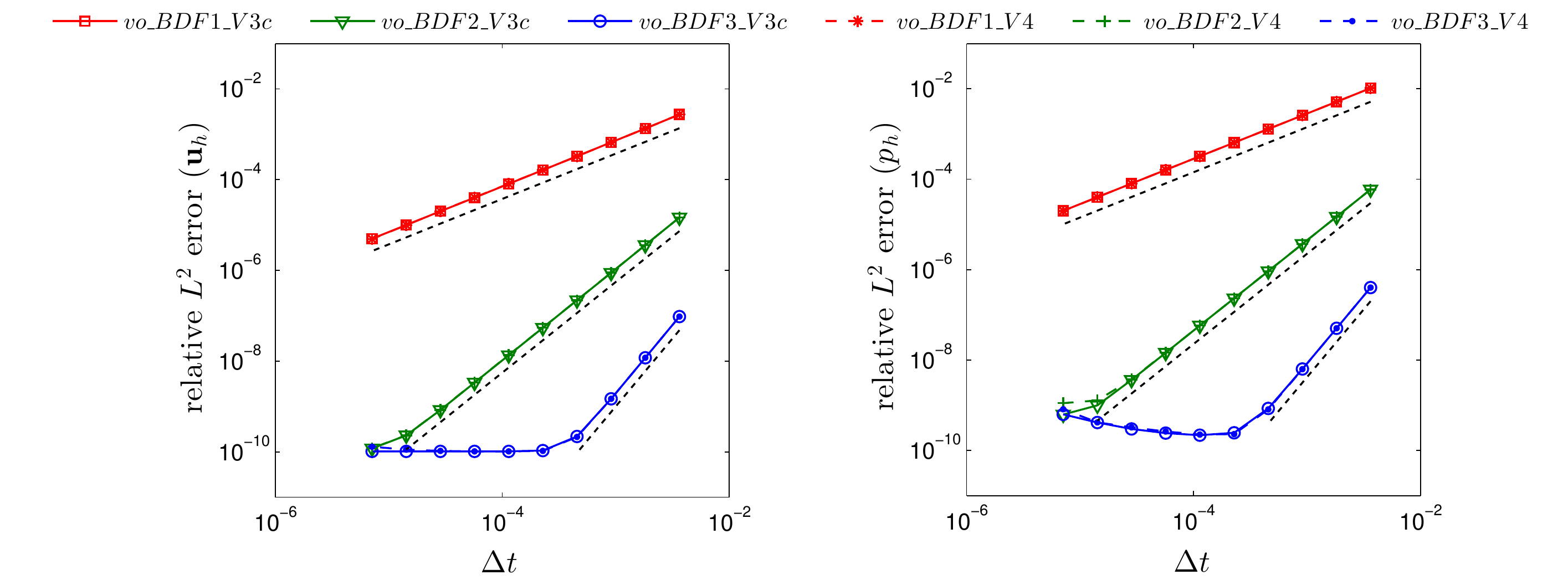}
\begin{picture}(100,0)
\put(15,69){\footnotesize $\sim \Delta t^3$}
\put(224,69){\footnotesize $\sim \Delta t^3$}
\put(15,106){\footnotesize $\sim \Delta t^2$}
\put(224,109){\footnotesize $\sim \Delta t^2$}
\put(15,144){\footnotesize $\sim \Delta t^1$}
\put(224,149){\footnotesize $\sim \Delta t^1$}
\end{picture}
\caption{Temporal convergence of \emph{V3c} and \emph{V4} using $J=J_p=1$ ($BDF1$), $J=J_p=2$ ($BDF2$) and $J=J_p=3$ ($BDF3$).}
\label{fig:v_convergence_time}
\end{figure}

\subsubsection{Temporal convergence}
\label{sec:vo_t}
Simulations are performed for this example to investigate convergence for temporal orders of accuracy $J=J_p=\{1,2,3\}$ in addition to the computations shown in Section~\ref{sec:varrhs} using $J=J_p=1$. The time step is chosen starting from $\mathrm{CFL}=2$ by successive bisection down to $\mathrm{CFL}=0.0039$ for all cases. We have observed that the present scheme allows computations with CFL numbers beyond the common limit of $\mathrm{CFL}=1$ for this flow and that there is no major difference in stability limits comparing the second and third order method, which is in contrast to the conditional stability for $J=3$ reported in~\cite{leriche2006numerical}. The spatial discretization uses $N_e=8^2$ elements and a polynomial degree of $k=7$. The results presented in Figure~\ref{fig:v_convergence_time} are labeled accordingly through the temporal scheme BDF1 for $J=J_p=1$, BDF2 for $J=J_p=2$, BDF3 for $J=J_p=3$ as well as the respective variant \emph{V3c} and \emph{V4}.

The results in Figure~\ref{fig:v_convergence_time} show optimal convergence rates for the respective order until the spatial error becomes predominant. There is essentially no difference between the variants \emph{V3c} and \emph{V4}. Please note that earlier studies frequently present the current splitting scheme with mixed temporal orders, e.g. $J=2$ and $J_p=1$, for which the temporal accuracy according to~\cite{guermond2003velocity} is at most second order for the velocity and order $3/2$ for the pressure. By choosing $J=J_p$ we get optimal convergence rates for both velocity and pressure at virtually no additional computational cost.

\begin{figure}[t!]
\centering
\includegraphics[trim= 0mm 0mm 0mm 0mm,clip,scale=0.6]{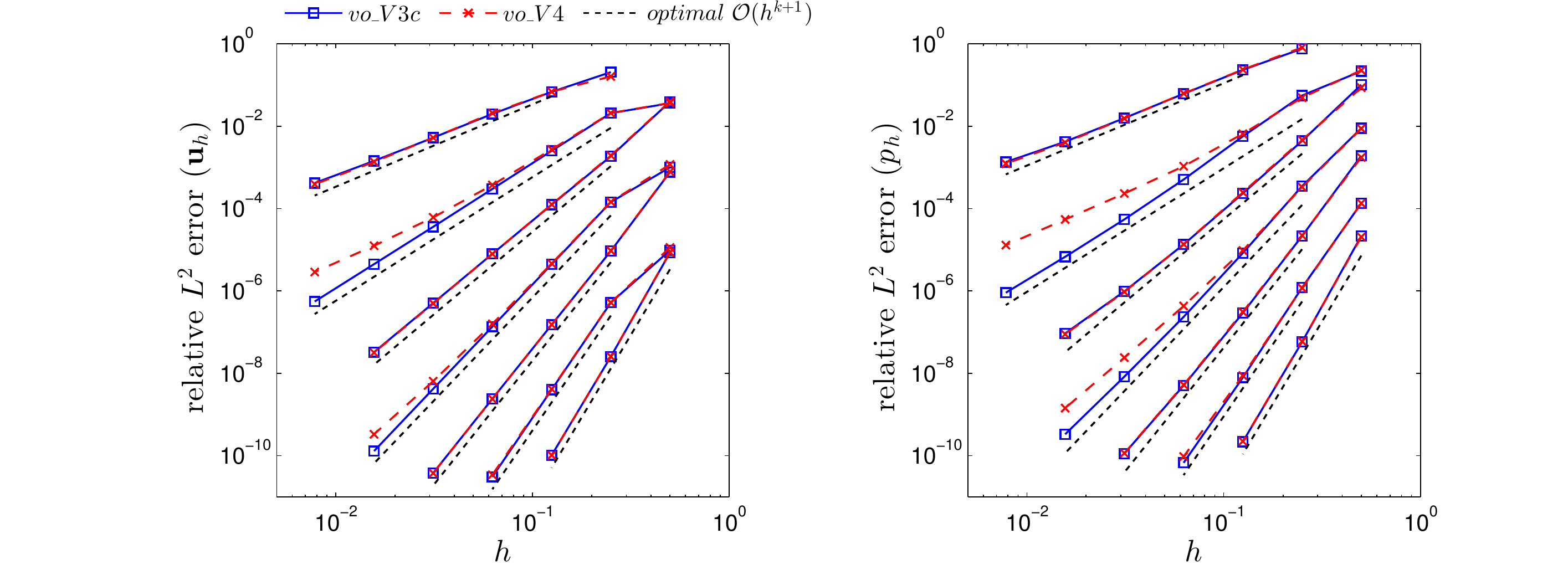}
\begin{picture}(100,0)
\put(-90,126){\footnotesize $k1$}
\put(-90,89){\footnotesize $k2$}
\put(-77,82){\footnotesize $k3$}
\put(-77,52){\footnotesize $k4$}
\put(-59,45){\footnotesize $k5$}
\put(-40,43){\footnotesize $k6$}
\put(-15,43){\footnotesize $k7$}
\put(127,133){\footnotesize $k1$}
\put(127,92){\footnotesize $k2$}
\put(139,88){\footnotesize $k3$}
\put(139,57){\footnotesize $k4$}
\put(158,51){\footnotesize $k5$}
\put(177,48){\footnotesize $k6$}
\put(201,47){\footnotesize $k7$}
\end{picture}
\caption{Spatial convergence for several polynomial degrees $k\{1,2,3,4,5,6,7\}$.}
\label{fig:v_convergence_space}
\vspace{0.7cm}
\includegraphics[trim= 0mm 0mm 0mm 0mm,clip,scale=0.6]{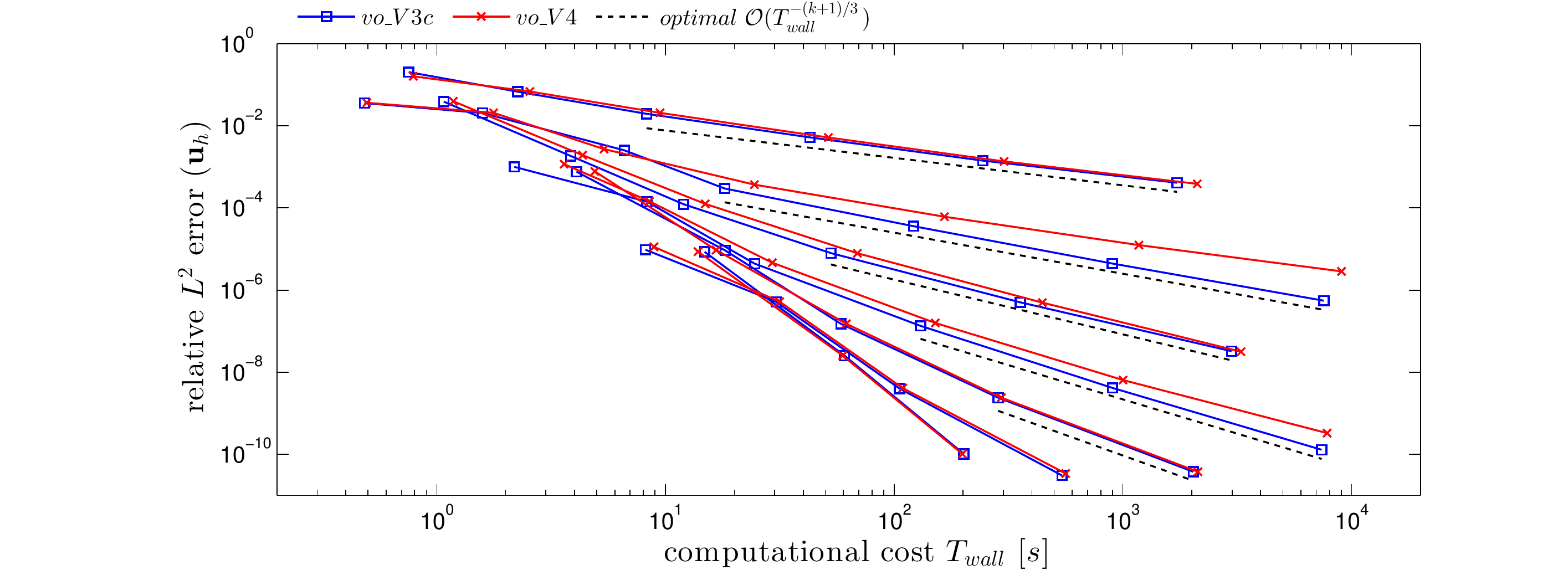}
\begin{picture}(100,0)
\put(195,135){\footnotesize $k1$}
\put(236,99){\footnotesize $k2$}
\put(210,83){\footnotesize $k3$}
\put(235,52){\footnotesize $k4$}
\put(196,45){\footnotesize $k5$}
\put(154,44){\footnotesize $k6$}
\put(123,50){\footnotesize $k7$}
\end{picture}
\caption{Computational cost in terms of the wall time of a serial computation $T_{wall}$ for the vortex problem using several polynomial orders $k\{2,3,4,5,6,7\}$.}
\label{fig:v_performance_space}
\end{figure}

\subsubsection{Spatial convergence}
\label{sec:vo_h}
We proceed with an investigation of the spatial convergence using the same example. $h$-refinement studies are performed for the polynomial degrees $k=\{1,2,3,4,5,6,7\}$ considering both \emph{V3c} and \emph{V4} with the time step chosen according to the CFL condition as $\mathrm{CFL}=0.0625$ and $J=J_p=3$. This way, the spatial error is dominant.

It may be observed in Figure~\ref{fig:v_convergence_space} that optimal convergence rates of order $k+1$ are obtained for \emph{V3c} for all cases both for the velocity and pressure. For variant \emph{V4}, we get sub-optimal convergence rates for the polynomial orders $k=2$ and $4$, while all other cases converge optimally. Further investigation of this issue showed that partial integration of the pressure gradient similar to \emph{V3c} may resolve this issue, we do not consider this further variant in this work, however. For the three to four coarsest refinement levels, the error is not distinguishable between the two variants \emph{V3c} and \emph{V4}.

\subsubsection{Performance evaluation}
\label{sec:performance_k}
Which polynomial degree yields the most efficient algorithm for laminar flows? In order to provide a first answer to this question, Figure~\ref{fig:v_performance_space} shows the error of the vortex problem according to Figure~\ref{fig:v_convergence_space} for $k\ge1$ over the wall time of the simulations $T_{wall}$. All simulations have been conducted in serial on the same computational setup, meaning that the wall time may be interpreted as the computational cost. The graph illustrates that increasing polynomial orders result in a steeper slope, i.e. the error decreases more rapidly if larger computational effort is invested. The expected slope of these curves of order $\mathcal{O}(T_{wall}^{-(k+1)/(d+1)})$ (with $d=2$ for this example) is obtained by considering a decrease of the error with $h^{k+1}$ and an increase in computation time with $T_{wall}\sim h^{-(d+1)}$ due to the number of elements to be evaluated on the one hand as well as the CFL condition on the other hand. The curves depicted in Figure~\ref{fig:v_performance_space} show very good agreement with this slope starting from the fourth refinement level which confirms the optimality of the code for high refine levels. The cause for the discrepancy regarding the coarser meshes is due to the efficiency of the particular multigrid algorithm, which results in decreased performance for a small number of refinement levels than for finer meshes. 

We conclude from this investigation that high-order methods are very efficient if high precision is required while lower polynomial degrees may be advantageous if a fast time to solution at reduced accuracy is desired. Please note that this discussion only presents a rough estimation since we do not tune the CFL condition or the relative error tolerances of the linear solvers for these cases, which would result in reduced computation times for moderate precision and low polynomial degrees. Considering high accuracies it is yet the slope that is more relevant than the absolute wall time. Further, computation times of variant \emph{V4} are slightly elevated in comparison to \emph{V3c} due to the global character of the projection equation system.

\begin{figure}[htb]
\begin{minipage}[b]{0.22\linewidth}
\centering
\begin{picture}(10,0)
\put(-25,30){mesh, coarsest level}
\end{picture}
\end{minipage}
\begin{minipage}[b]{0.75\linewidth}
\centering
\includegraphics[trim= 0mm 68mm 0mm 68mm,clip,scale=0.2]{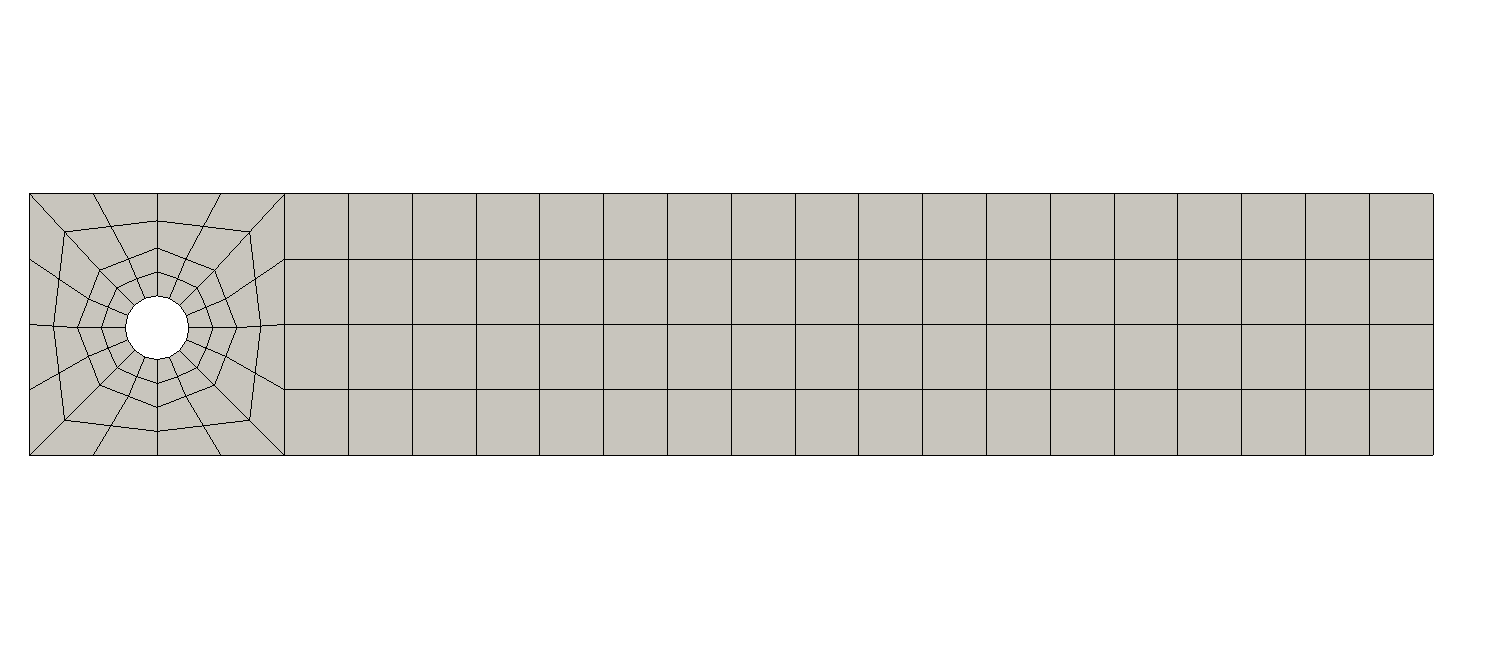}
\end{minipage}
\\
\begin{minipage}[b]{0.22\linewidth}
\centering
\begin{picture}(10,0)
\put(-25,30){velocity magnitude}
\end{picture}
\end{minipage}
\begin{minipage}[b]{0.75\linewidth}
\centering
\includegraphics[trim= 0mm 68mm 0mm 68mm,clip,scale=0.2]{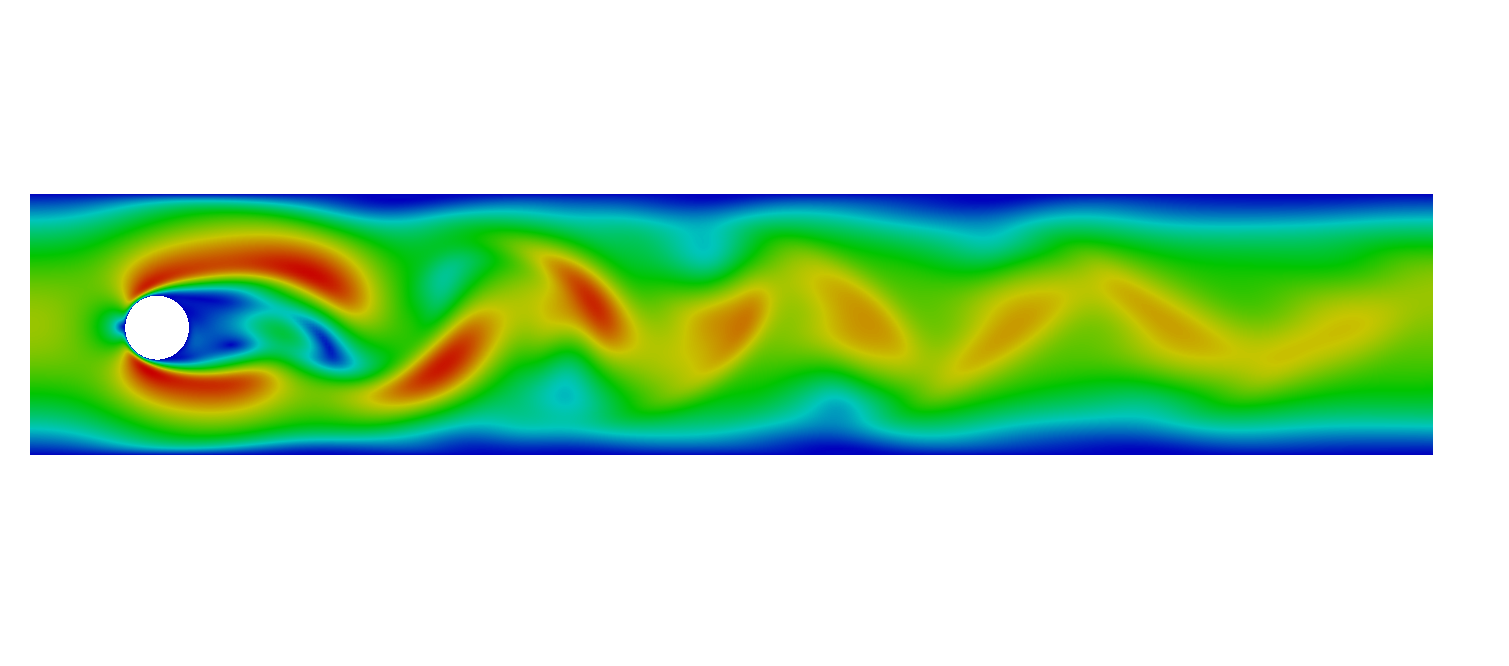}
\end{minipage}
\\
\begin{minipage}[b]{0.22\linewidth}
\centering
\begin{picture}(10,0)
\put(-25,30){vorticity}
\end{picture}
\end{minipage}
\begin{minipage}[b]{0.75\linewidth}
\centering
\includegraphics[trim= 0mm 68mm 0mm 68mm,clip,scale=0.2]{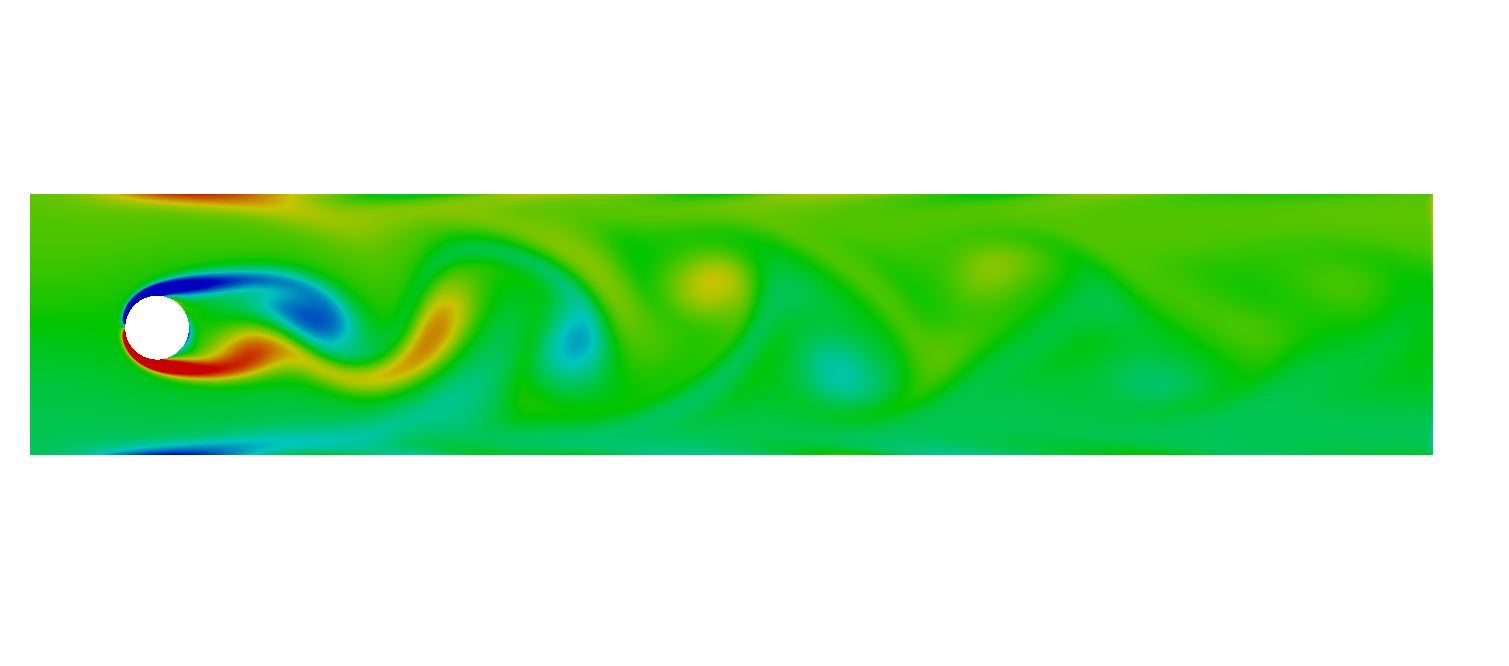}
\end{minipage}
\\
\begin{minipage}[b]{0.22\linewidth}
\centering
\begin{picture}(10,0)
\put(-25,30){pressure}
\end{picture}
\end{minipage}
\begin{minipage}[b]{0.75\linewidth}
\centering
\includegraphics[trim= 0mm 68mm 0mm 68mm,clip,scale=0.2]{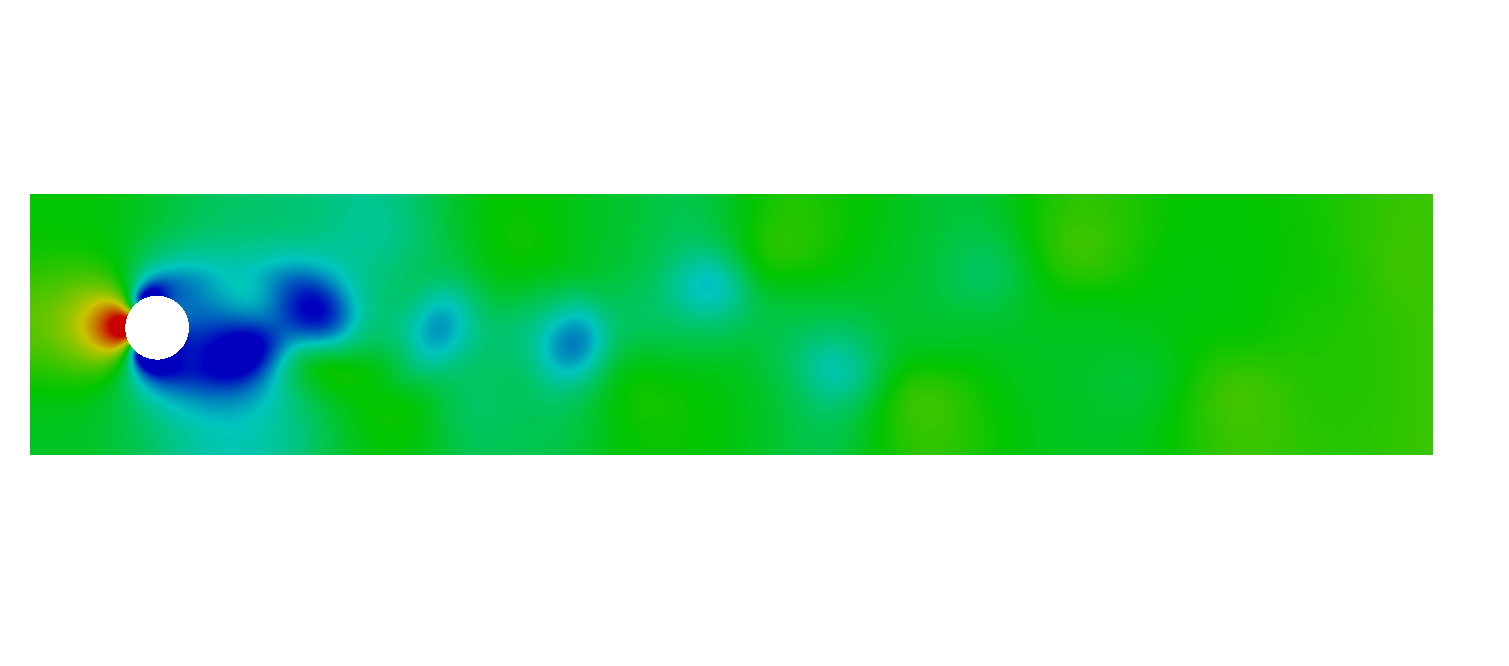}
\end{minipage}
\caption{Flow past cylinder from top to bottom: coarsest mesh level, velocity magnitude, vorticity and pressure snapshot at time $t=5.5$. The latter three have been computed with the finest mesh including $k=6$ and \emph{V3c}. Red indicates high and blue low values.}
\label{fig:fpc_picture}
\end{figure}

\begin{table}[htb]
\caption{Flow past cylinder cases, resolutions and results.}
\label{tab:fpc_cases}
\begin{tabular*}{\textwidth}{l| @{\extracolsep{\fill}} l l l l l l l l}
\hline
\multicolumn{1}{l}{} & & \multicolumn{3}{c}{\emph{V3c}} & \multicolumn{3}{c}{\emph{V4}} \\ \cline{3-5} \cline{6-8}
\multicolumn{1}{l}{$k$} & $n DoFs$ & $c_{D\max}$ & $c_{L\max}$ & $\Delta p_{\text{end}}$ & $c_{D\max}$ & $c_{L\max}$ & $\Delta p_{\text{end}}$ \\ \hline \noalign{\smallskip}
\multirow{3}{*}{4} 
& 10,200 & 2.767868 & 0.437474 & -0.108995 & 2.629942 & 0.331907 & -0.108456  \\
& 40,800 & 2.963629 & 0.487026 & -0.111627 & 2.970182 & 0.491144 & -0.111988  \\
& 163,200 & 2.950792 & 0.478365 & -0.111654  & 2.950915 & 0.478322 & -0.111696 \\ \noalign{\smallskip}
\multirow{3}{*}{5} 
& 14,688 & 2.958448 & 0.512599 & -0.109378 & 2.926823 & 0.512408 & -0.108607 \\
& 58,752 & 2.951065 & 0.480343 & -0.111616 & 2.953314 & 0.476367 & -0.111888 \\
& 235,008 & 2.950454 & 0.477967 & -0.111618 & 2.950350 & 0.477891 & -0.111616 \\ \noalign{\smallskip}
\multirow{3}{*}{6} 
& 19,992 & 2.964686 & 0.504783 & -0.111059 & 2.975162 & 0.514634 & -0.110940 \\
& 79,968 & 2.949137 & 0.478317 & -0.111614 & 2.948503 & 0.478380 & -0.111612 \\
& 319,872 & 2.950829 & 0.477940 & -0.111615 & 2.950798 & 0.477899 & -0.111616 \\ \noalign{\smallskip}
\multirow{3}{*}{7} 
& 26,112 & 2.949849 & 0.487586 & -0.111414 & 2.950902 & 0.484617 & -0.111446 \\
& 104,448 & 2.950198 & 0.477948 & -0.111615 & 2.949882 & 0.47827 & -0.111584 \\ 
& 417,792 & 2.950927 & 0.477941 & -0.111615 & 2.950944 & 0.477899 & -0.111619 \\ \noalign{\smallskip}
2 & $\sim$500,000~\cite{John04}
& 2.95092 & 0.47795 & -0.1116 & 2.95092 & 0.47795 & -0.1116  \\ \noalign{\smallskip}
\multicolumn{2}{l}{lower bound~\cite{Schaefer96}}
& 2.9300 & 0.4700 & -0.1150 & 2.9300 & 0.4700 & -0.1150   \\
\multicolumn{2}{l}{upper bound~\cite{Schaefer96}}
& 2.9700 & 0.4900 & -0.1050 & 2.9700 & 0.4900 & -0.1050  \\ \hline
\end{tabular*}
\end{table}

\subsection{Laminar flow past cylinder}
\label{sec:fpc}
As second laminar benchmark we investigate unsteady vortex shedding in the wake of a cylinder presented as test case 2D-3 in \cite{Schaefer96} with accurate reference data provided in \cite{John04}. The domain is of dimensions $W \times H$ in streamwise and vertical direction, respectively, with $W=2.2$ and $H=0.41$. The cylinder is of diameter $0.1$ and its center point is located at $0.2$ units from the inflow as well as the bottom walls. At the inflow boundary, the velocity is prescribed by
\begin{equation}
g_{\bm{u}1}(x_2,t)=U_m\frac{4x_2(H-x_2)}{H^2}\sin(\pi t /\mathcal{T})
\end{equation}
with $U_m=1.5$ and $\mathcal{T}=8$. At the top and bottom wall as well as the cylinder surface $g_{\bm{u}2}=0$ and no-slip boundary conditions are applied. At the outlet, zero pressure boundary conditions are applied with $g_p=0$ and $\bm{h}=\bm{0}$. The resulting flow exhibits unsteady vortex shedding behind the cylinder as illustrated in Figure~\ref{fig:fpc_picture} at time instant $t=5.5$.

We perform simulations for the spatial polynomial orders $k=\{4,5,6,7\}$ and present three levels of refinement for each polynomial degree. Approximation of the cylindrical geometry is enhanced by mapping the boundary nodes onto the cylinder surface yielding an accurate representation of the geometry for high-order polynomials. The resulting mesh of the coarsest level is displayed in Figure~\ref{fig:fpc_picture}. As temporal discretization we choose the third-order accurate scheme with $J=J_p=3$ and a CFL number of $\mathrm{CFL}=0.25$ based on $U_m$ and the minimum edge length of the respective mesh $h_{\min}$. The solution quality is evaluated according to~\cite{Schaefer96} by the maximum value of the drag and lift coefficients $c_{D\max}$ and $c_{L\max}$ over time as well as the pressure difference between the windward and the lee side of the cylinder at the end of the simulation, denoted $\Delta p_{\text{end}}= \Delta p(t=\mathcal{T})$. A detailed description of how these quantities are computed is given in~\cite{Schaefer96,John04}. For the variant \emph{V3c} the standard penalty parameter $\zeta_D^* = 1$ is selected while the factors $\zeta_C^*=\zeta_D^* = 10$ are necessary for \emph{V4} to obtain a stable scheme.

The results of our simulations are displayed in Table~\ref{tab:fpc_cases} along with reference data by~\cite{John04}, where $c_{D\max}$ is of absolute accuracy $5 e -7$ and $c_{L\max}$ as well as $\Delta p_{\text{end}}$ are of accuracy $1 e -4$. In addition, upper and lower bounds for all three quantities as presented in~\cite{Schaefer96} are also included.

The results shown in Table~\ref{tab:fpc_cases} exhibit excellent agreement with reference data. Especially the pressure difference is already predicted with the same accuracy as the reference data for 40,800 degrees of freedom with \emph{V3c} and reaches two additional digits in precision during refinement. This high level of accuracy may be due to the equal-order approach employed in the present method in contrast to the mixed-order method utilized for the reference computations in~\cite{John04}. The drag coefficient also converges to reference data yielding an accuracy of at least five digits and excellent agreement between \emph{V3c} and \emph{V4}. The lift coefficient converges to 0.4779 using \emph{V3c} and \emph{V4} which is in excellent agreement with reference data given as 0.47795 where the error was specified to be no larger than $1 e -4$.

\section{Application to DNS and LES of turbulent channel flow}
\label{sec:dles}
We demonstrate applicability of the present code to direct numerical simulation (DNS) as well as large-eddy simulation (LES) of turbulent flow. DNS is the most accurate approach to computing turbulent flows, since all turbulent scales are resolved by the numerical method, and is therefore a valuable tool in research and industry, e.g. for computing reference data or investigating selected low-Reynolds-number configurations, see, e.g., \cite{Moin98} for an overview. As discussed in the first Subsection~\ref{sec:dns} for the example of turbulent channel flow, the current approach is very attractive for computing DNS due to its high-order accuracy, speed and high scalability on massively parallel computers. DNS are yet out of reach in the foreseeable future regarding most industrial applications in the moderate to large Reynolds number regime, which is why we show that the present approach is also very well suitable for implicit large-eddy simulation in Subsection~\ref{sec:iles}.

\begin{table}
\caption{Channel flow cases and resolutions. Specification of the mesh quantities $\Delta (\cdot)^+=\Delta (\cdot)_e u_{\tau}/ \nu (k+1)$ is given as usual as the respective element length $\Delta (\cdot)_e$ divided by the number of nodes in each spatial direction per element $k+1$ where $k$ is the polynomial degree. $\Delta x^+$: resolution in $x_1$-direction; $\Delta y_1^+$: first off-wall point in $x_2$-direction; $\Delta y_c^+$: resolution at center in $x_2$-direction; $\Delta z^+$: resolution in $x_3$-direction.}
\label{tab:ch_flows}
\begin{tabular*}{\textwidth}{l @{\extracolsep{\fill}} l l l l l l l l l l}
\hline
Case     & $N_{e,1} \times N_{e,2} \times N_{e,3}$  & $k$ & $n DoFs$ & CFL &$Re_{\tau}$   & $\gamma$ & $\Delta x^+$ & $\Delta  y_1^+$ & $\Delta y_c^+$ & $\Delta z^+$
\\ \hline \noalign{\smallskip}
$ch180\_N32^3\_k5\_V3c\_dns$   & $32 \times 32 \times 32$    & $5$ & $28e6$ & 1 & $180$& $1.4$ & $11.8$& $0.69$&$3.0$& $3.9$\\
$ch590\_N64^3\_k4\_V3c\_dns$   & $64 \times 64 \times 64$    & $4$ & $131e6$ & 0.8 & $590$& $1.65$ & $11.6$& $0.94$&$6.5$& $5.8$\\
\noalign{\smallskip}
$ch180\_N8^3\_k3\_V3c\_les$   & $8 \times 8 \times 8$    & $3$ & $0.13e6$ & 1 & $180$&  $1.8$ & $35.3$& $3.5$&$20.1$& $17.7$\\
$ch180\_N16^3\_k3\_V\{3c,4\}\_les$   & $16 \times 16 \times 16$    & $3$& $1.0e6$ & 1 & $180$& $1.8$ & $17.7$& $1.4$&$10.5$& $8.8$\\
$ch180\_N32^3\_k3\_V3c\_les$   & $32 \times 32 \times 32$    & $3$ & $8.4e6$ & 1 & $180$& $1.8$ & $8.8$& $0.62$&$5.3$& $4.4$\\
$ch590\_N16^3\_k4\_V\{3c,4\}\_les$   & $16 \times 16 \times 16$    & $4$ & $2.0e6$ & 1 & $590$& $2.25$&  $46.3$& $2.0$&$33.1$& $23.2$\\
$ch590\_N32^3\_k4\_V3c\_les$   & $32 \times 32 \times 32$    & $4$ & $16e6$ & 1 & $590$& $2.25$ & $23.2$& $0.85$&$16.9$& $11.6$\\
\hline
\end{tabular*}
\end{table}

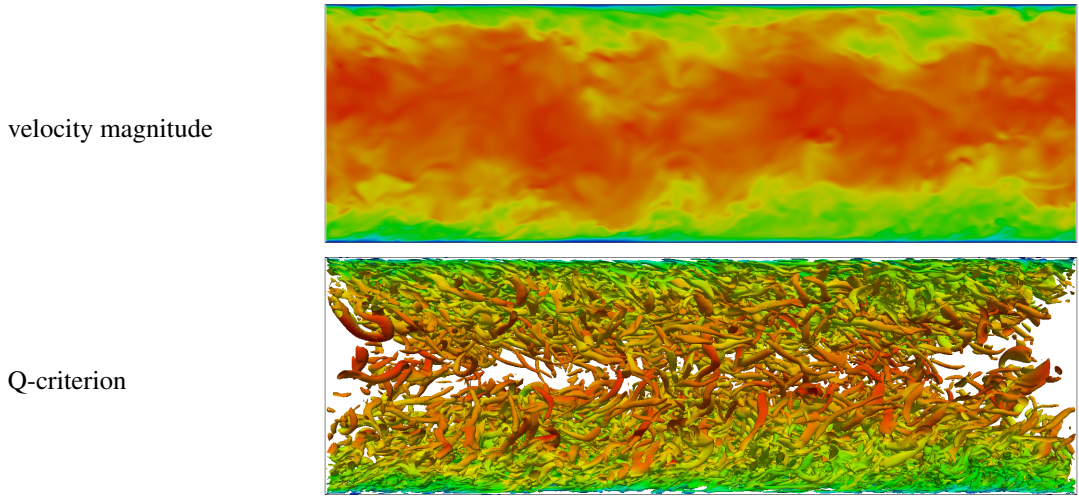
\begin{figure}[htb]
\begin{minipage}[b]{0.2\linewidth}
\centering
\begin{picture}(10,0)
\put(-25,42){velocity magnitude}
\end{picture}
\end{minipage}
\begin{minipage}[b]{0.77\linewidth}
\centering
\includegraphics[trim= 0mm 25mm 0mm 25mm,clip,width=0.79\textwidth]{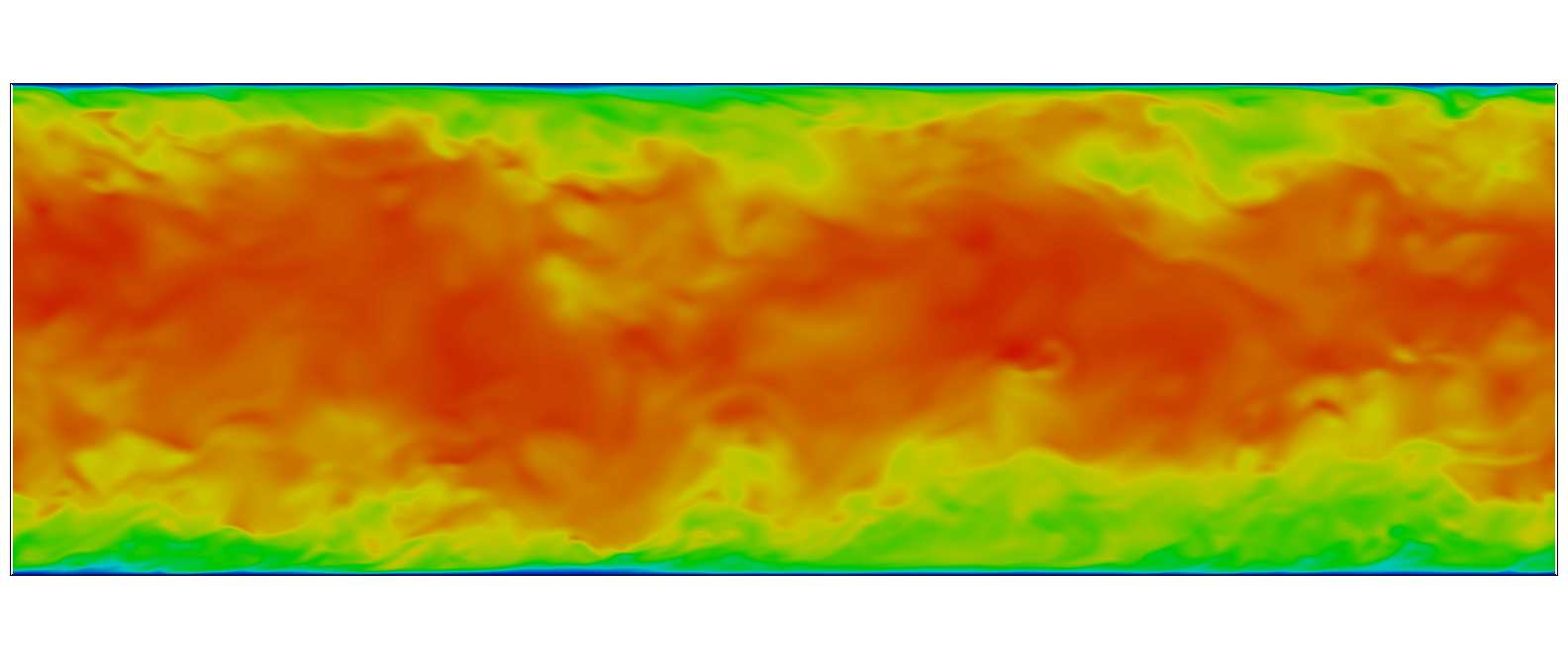}
\end{minipage}
\\
\begin{minipage}[b]{0.2\linewidth}
\centering
\begin{picture}(10,0)
\put(-25,42){Q-criterion}
\end{picture}
\end{minipage}
\begin{minipage}[b]{0.77\linewidth}
\centering
\includegraphics[trim= 0mm 25mm 0mm 25mm,clip,width=0.79\textwidth]{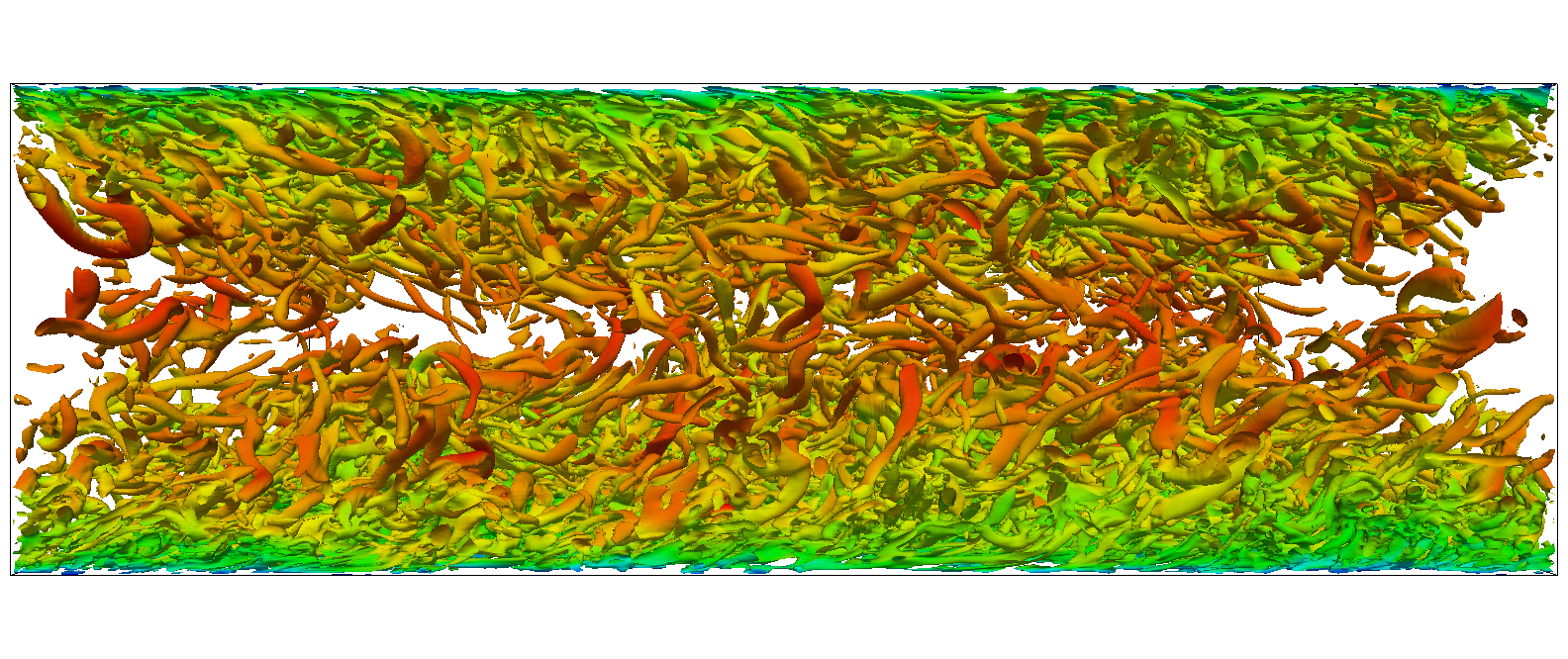}
\end{minipage}
\caption{DNS of turbulent channel flow at $Re_{\tau}=590$: Contour of velocity magnitude (top) and eddies visualized via the Q-criterion, colored by velocity magnitude (bottom). High velocity is colored red and low velocity blue.}
\label{fig:ch_picture}
\end{figure}

\begin{figure}[t!]
\centering
\begin{minipage}[b]{0.497\linewidth}
\centering
\includegraphics[trim= 9mm 0mm 15mm 0mm,clip,width=1.\textwidth]{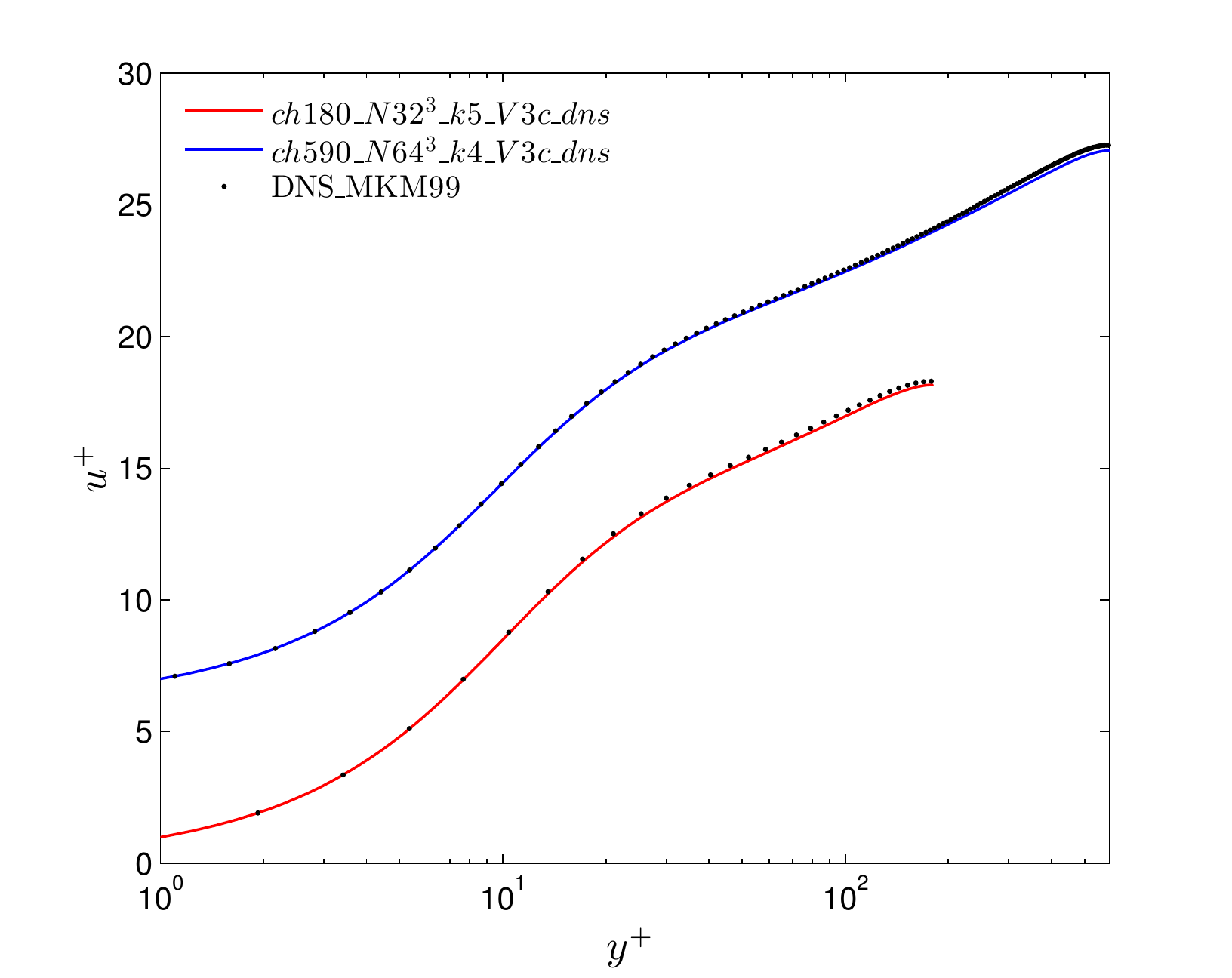}
\end{minipage}
\begin{minipage}[b]{0.497\linewidth}
\centering
\includegraphics[trim= 12mm 0mm 12mm 0mm,clip,width=1.\textwidth]{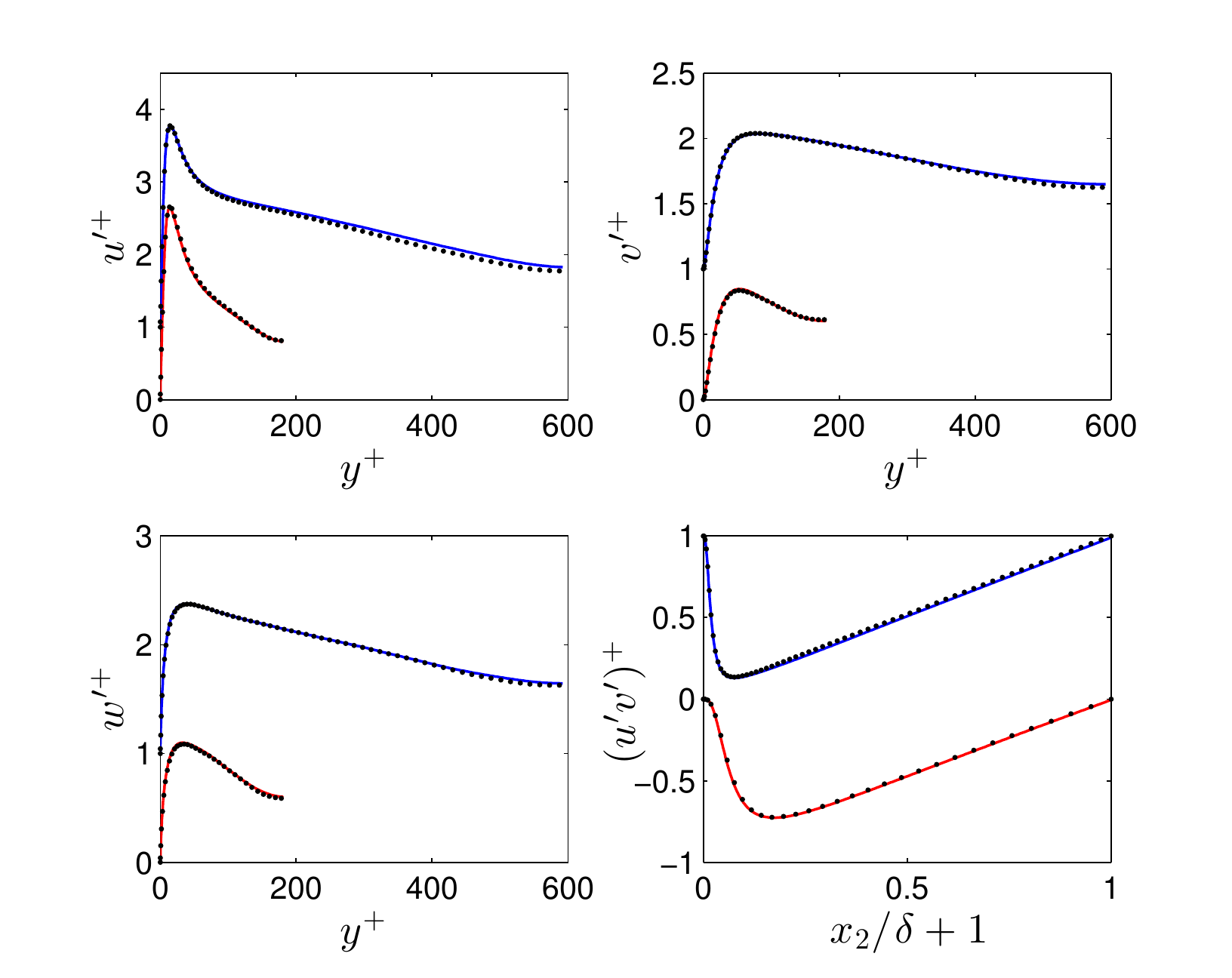}
\end{minipage}
\begin{picture}(100,0)
\put(-155,59){\footnotesize $Re_{\tau}=180$}
\put(-155,94){\footnotesize $Re_{\tau}=590$}
\end{picture}
\caption{DNS of turbulent channel flow at $Re_{\tau}=180$ and 590: Mean velocity $u^+=u_1/u_{\tau}$ (left) and root mean square velocities $u^{\prime+}=\mathrm{RMS}(u_1)/u_{\tau}$, $v^{\prime+}=\mathrm{RMS}(u_2)/u_{\tau}$ and $w^{\prime+}=\mathrm{RMS}(u_3)/u_{\tau}$ as well as Reynolds shear stresses $(u^{\prime}v^{\prime})^{+}=(u_1u_2)/u_{\tau}^2$ (right). For the case $Re_{\tau}=590$, the mean velocity is shifted upwards by six units and all other quantities by one unit for clarity.}
\label{fig:dnsum}
\centering
\begin{minipage}[b]{0.497\linewidth}
\centering
\includegraphics[trim= 9mm 0mm 15mm 0mm,clip,width=1.\textwidth]{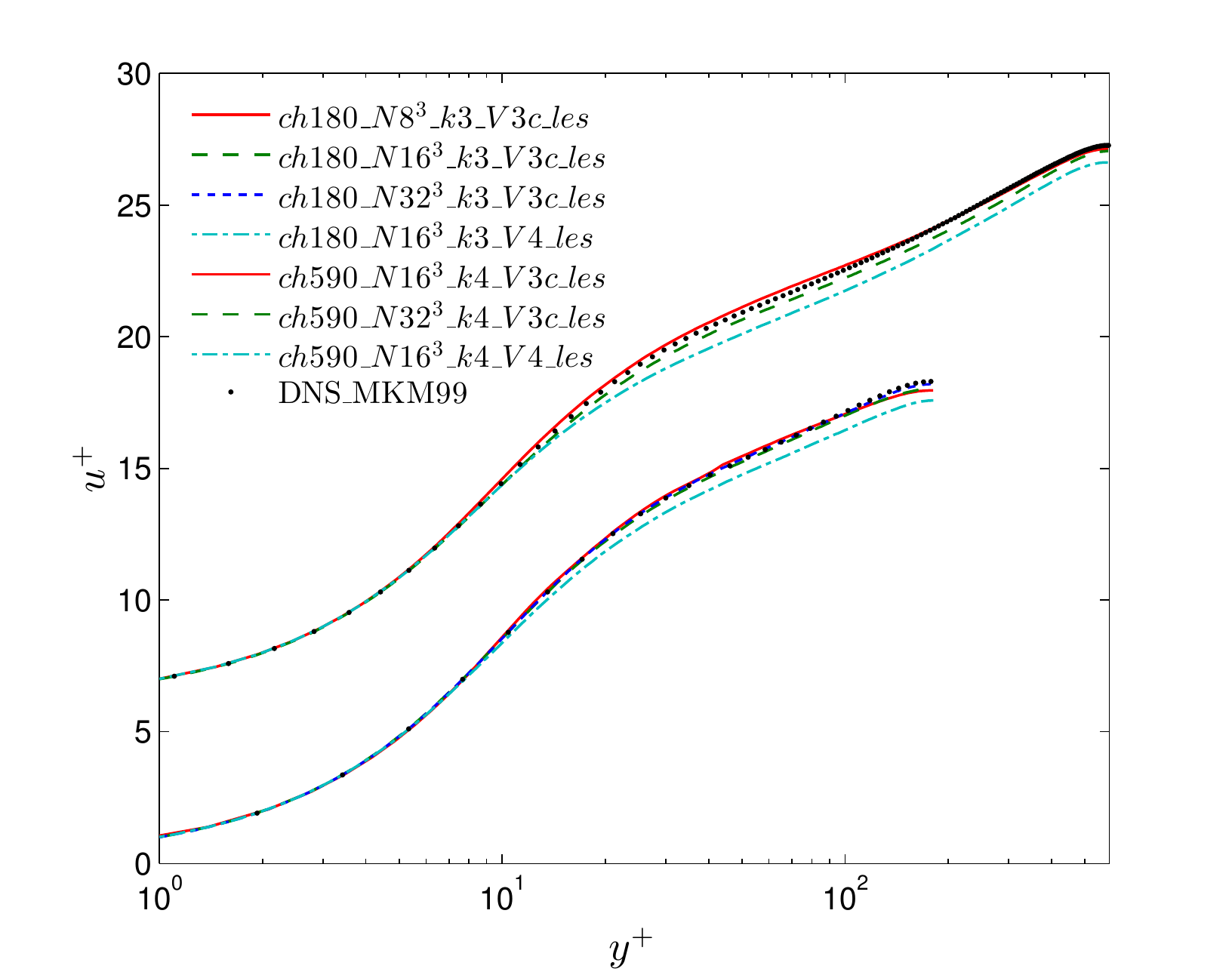}
\end{minipage}
\begin{minipage}[b]{0.497\linewidth}
\centering
\includegraphics[trim= 12mm 0mm 12mm 0mm,clip,width=1.\textwidth]{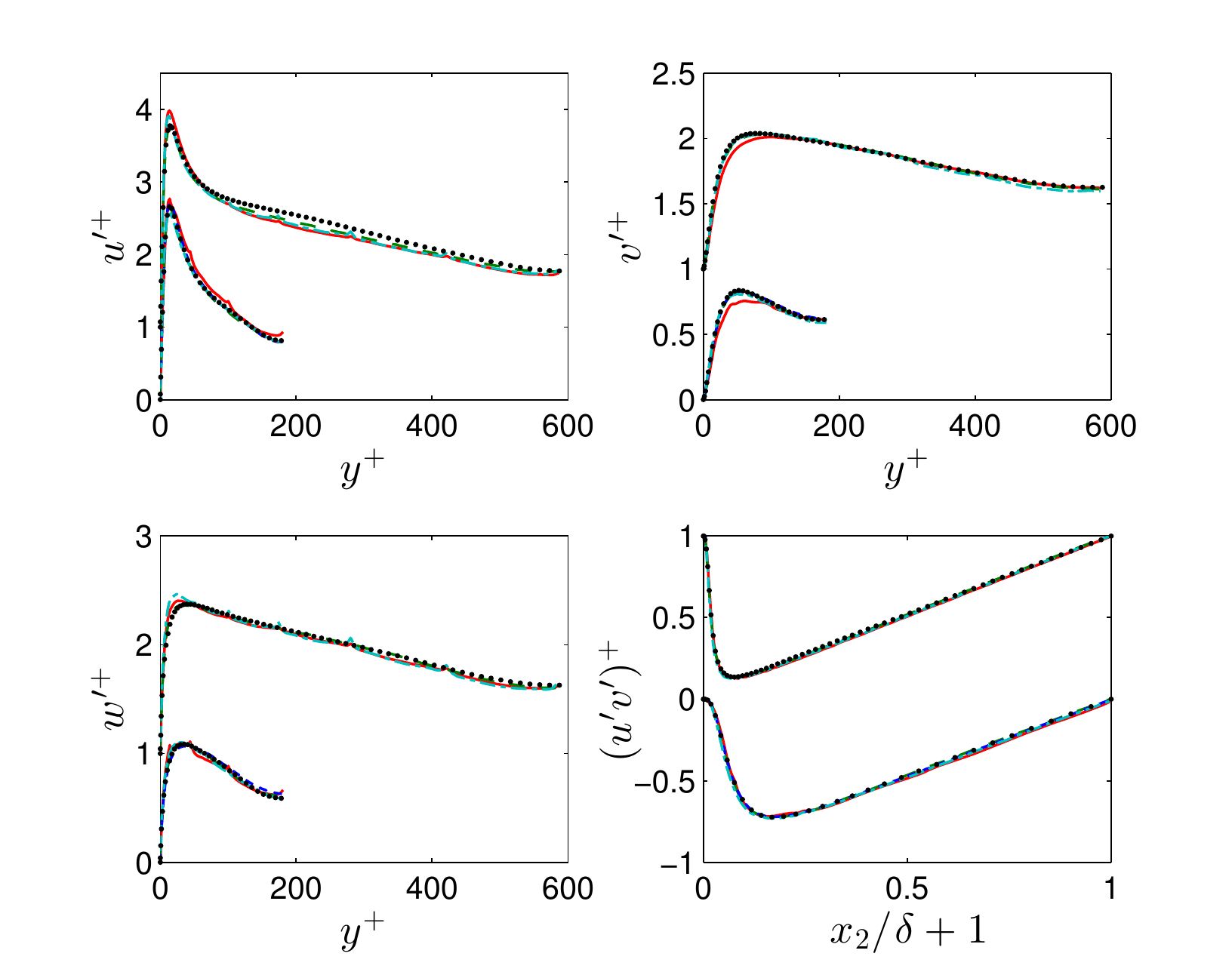}
\end{minipage}
\begin{picture}(100,0)
\put(-155,59){\footnotesize $Re_{\tau}=180$}
\put(-155,94){\footnotesize $Re_{\tau}=590$}
\end{picture}
\caption{ILES of turbulent channel flow at $Re_{\tau}=180$ and 590: Mean velocity $u^+=u_1/u_{\tau}$ (left) and root mean square velocities $u^{\prime+}=\mathrm{RMS}(u_1)/u_{\tau}$, $v^{\prime+}=\mathrm{RMS}(u_2)/u_{\tau}$ and $w^{\prime+}=\mathrm{RMS}(u_3)/u_{\tau}$ as well as Reynolds shear stresses $(u^{\prime}v^{\prime})^{+}=(u_1u_2)/u_{\tau}^2$ (right). For the case $Re_{\tau}=590$, the mean velocity is shifted upwards by six units and all other quantities by one unit for clarity.}
\label{fig:les180590um}
\end{figure}

\subsection{Direct numerical simulation}
\label{sec:dns}
There is a vast number of publications in the field of DNS of turbulent channel flow providing accurate reference data for turbulence modeling research as well as validation of numerical schemes, see for example Moser et al. \cite{Moser99} for friction Reynolds numbers $Re_{\tau}=180$, 395 and 590 and an overview over publications on DNS of the case $Re_{\tau}=180$ by Vreman \cite{Vreman14}. Herein, $Re_{\tau}$ is defined as usual as $Re_{\tau}=u_{\tau} \delta/\nu$ with given channel-half height $\delta$ and friction velocity $u_{\tau}=\sqrt{\tau_{w}/\rho}$, where $\tau_{w}$ is the wall shear stress and $\rho$ the density. In the present work, we perform DNS of turbulent channel flow at $Re_{\tau}=180$ and 590. Computational domain sizes for these flows are specified in \cite{Moser99} as $4\pi \delta \times 2\delta \times \frac{4}{3} \pi \delta$ in streamwise, wall-normal and spanwise direction, respectively, for the case $Re_{\tau}=180$ and $2\pi \delta \times 2\delta \times \pi \delta$ for the case $Re_{\tau}=590$ accordingly. Periodic boundary conditions are considered in the streamwise and spanwise directions and no-slip boundary conditions are imposed at the walls with $\bm{g}_{\bm{u}}=\bm{0}$. The mesh is graded towards the no-slip boundaries to improve resolution of near-wall turbulent structures according to the hyperbolic mesh mapping given as $f$: $[0,1] \to [-\delta, \delta]$:
\begin{equation}
x_2 \mapsto f(x_2)=\delta \frac{\tanh(\gamma (2x_2-1))}{\tanh(\gamma)}
\end{equation}
using the mesh-stretching parameter $\gamma$ according to Table~\ref{tab:ch_flows}. The spatial discretizations employed are similar in resolution to~\cite{Moser99} and listed in Table~\ref{tab:ch_flows}. For the case $Re_{\tau}=180$, $32\times32\times32$ elements of degree 5 are used, resulting in 28.3 million DoFs overall, which is slightly finer than in \cite{Moser99} with $128\times129\times128$ Fourier modes in the periodic directions as well as Chebyshev nodes in the wall-normal direction and 8.5 million DoFs overall. Regarding $Re_{\tau}=590$, a mesh of $64\times64\times64$ elements of degree 4 with 131 million degrees of freedom is used again compared to a mesh of $384\times257\times384$ nodes and 152 million DoFs in \cite{Moser99}. The time step of the BDF3 scheme ($J=J_p=3$) is chosen based on the CFL condition according to $\Delta t = \mathrm{CFL} h_{\mathrm{min}} / U k^2$ where we take $\mathrm{CFL}=1$ for $Re_{\tau}=180$ and $\mathrm{CFL}=0.8$ for $Re_{\tau}=590$, $h_{\mathrm{min}}$ as the minimum edge length and $U=15 u_{\tau}$ representing the estimated maximum velocity occurring in the cells with the shortest edge length, which are located at the no-slip boundaries. Statistics are sampled spatially over homogeneous planes and temporally over approximately 30 flow-through times based on the mean center-line velocity for the case $Re_{\tau}=180$ and 68 flow-through times for $Re_{\tau}=590$. For the DNS computations, we solely consider the variant \emph{V3c}.

The resulting turbulent flow is depicted in Figure~\ref{fig:ch_picture} for the case $Re_{\tau}=590$ via velocity contours as well as the Q-criterion for eddy visualization of one snapshot. The cases are labeled according to Table~\ref{tab:ch_flows} and the results are plotted over the wall coordinate $y^+=y u_{\tau} / \nu$ and $x_2/\delta$, respectively, in Figure~\ref{fig:dnsum}. Herein, the normalized mean velocity is defined as $u^+=u_1/u_{\tau}$ and fluctuations in form of the root mean square velocities as $u^{\prime+}=\mathrm{RMS}(u_1)/u_{\tau}$, $v^{\prime+}=\mathrm{RMS}(u_2)/u_{\tau}$ and $w^{\prime+}=\mathrm{RMS}(u_3)/u_{\tau}$ as well as the Reynolds shear stresses as $(u^{\prime}v^{\prime})^{+}=(u_1 u_2)/u_{\tau}^2$. The curves exhibit excellent agreement with reference data from \cite{Moser99} labeled $DNS\_MKM99$.

\subsection{Implicit large-eddy simulation}
\label{sec:iles}
Spatial resolution requirements of DNS scale with the Reynolds number as $Re^{9/4}$ according to~\cite{Case73}, which makes it too expensive for many industrial applications of high Reynolds number. This requirement is relaxed by solely resolving the larger eddies explicitly by the scheme and modeling the dissipative effect of the smaller ones, which yields spatial resolution requirements of $Re^{9/5}$ for boundary layer flows~\cite{Chapman79}. In this work, we adopt the assumption that the numerical dissipation required to stabilize the scheme is in fact an adequate representation of the physical subgrid scales, which renders a supplementary turbulence model unnecessary. This idea leads to the widely used concept of implicit large-eddy simulation (ILES); see, e.g., \cite{Grinstein07} for an extensive discussion, \cite{Beck14,Wiart15,Bassi16} for examples within compressible DG and \cite{Hughes95} for a rigorous derivation in the variational context.

As a benchmark example, we employ the turbulent channel flow setup similar to the previous subsection with a domain size of $2\pi \delta \times 2\delta \times \pi \delta$ for all computations. All cases considered are listed in Table~\ref{tab:ch_flows} and include three refinement levels for the case $Re_{\tau}=180$ and two refinement levels for $Re_{\tau}=590$. All cases are computed with the variant \emph{V3c} and two representative cases are shown employing variant \emph{V4}, which are labeled accordingly. For the case $ch590\_N16^3\_k4\_V4\_les$, $\zeta^*_C=\zeta^*_D=10$ is chosen in order to obtain a stable numerical method, while the standard values of $\zeta^*_C=\zeta^*_D=1$ were used for all other computations.

The results of all LES cases according to Table~\ref{tab:ch_flows} are depicted in Figure~\ref{fig:les180590um} along with DNS data from~\cite{Moser99}. Regarding variant \emph{V3c}, excellent agreement with reference solutions is observed and coarsening of the discretizations has very little effect on the solution quality. The mean velocity is slightly under-predicted for both cases of \emph{V4}, however.

A comparison of computation times yields an elevated computational cost of the case $ch180\_N16^3\_k3\_V4\_les$ by a multiplicative factor of 2.09 in comparison to $ch180\_N16^3\_k3\_V3c\_les$ using the same computational setup. Analogously, the case $ch1590\_N16^3\_k4\_V4\_les$ completes in 4.25 times the computational cost of $ch590\_N16^3\_k4\_V3c\_les$. These differences are largely due to the global equation system introduced in the projection step in \emph{V4}, which becomes more costly with increasing penalty parameter. It is therefore concluded from these investigations that \emph{V3c} is the most efficient variant discussed in the present paper for simulation of turbulent flows, both for DNS and LES, as it combines high accuracy with best computation times.

\section{Conclusion}
\label{sec:conclusion}
In this paper, we have developed a stable, accurate and efficient numerical scheme for simulation of the incompressible Navier--Stokes equations by reviewing, comparing and extending stabilization techniques proposed in literature. The best stabilization for small time steps is based on a div-div penalty approach that enhances the point-wise divergence-free condition within elements. Under-resolved flows have been stabilized by either including a supplementary jump-penalty term within the projection or by partial integration of the right-hand side of the Poisson equation. The resulting algorithm, especially with the latter variant, exhibits convergence orders equal to polynomial degree plus one in space and three in time both in velocity and pressure and is embedded in a matrix-free implementation.

The high efficiency of this implementation for high polynomial degrees has been demonstrated by two laminar flow examples present in a vortex problem and flow past a cylinder. This characteristic makes the present methodology also very attractive for computation of direct numerical simulation of turbulent flows, as demonstrated in this work with accurate prediction of two turbulent channel flow examples. Finally, we have also shown that this scheme is well suited for implicit large-eddy simulation of turbulent flow.

\section*{Acknowledgments}
The research presented in this paper was partly funded by the German Research Foundation (DFG) under the project ``High-order discontinuous Galerkin for the EXA-scale'' (ExaDG) within the priority program ``Software for Exascale Computing'' (SPPEXA), grant agreement no. KR4661/2-1 and WA1521/18-1. Computational resources on SuperMUC in Garching, Germany, provided by the Leibniz Supercomputing Centre, under the project pr83te are gratefully acknowledged.

\bibliography{dg_paper_v1}

\begin{thebibliography}{10}
\expandafter\ifx\csname url\endcsname\relax
  \def\url#1{\texttt{#1}}\fi
\expandafter\ifx\csname urlprefix\endcsname\relax\def\urlprefix{URL }\fi
\expandafter\ifx\csname href\endcsname\relax
  \def\href#1#2{#2} \def\path#1{#1}\fi

\bibitem{Hindenlang12}
F.~Hindenlang, G.~Gassner, C.~Altmann, A.~Beck, M.~Staudenmaier, C.-D. Munz,
  Explicit discontinuous {G}alerkin methods for unsteady problems, Comput.
  Fluids 61 (2012) 86--93.
\newblock \href {http://dx.doi.org/10.1016/j.compfluid.2012.03.006}
  {\path{doi:10.1016/j.compfluid.2012.03.006}}.

\bibitem{Bassi16}
F.~Bassi, L.~Botti, A.~Colombo, A.~Crivellini, A.~Ghidoni, F.~Massa, On the
  development of an implicit high-order discontinuous {G}alerkin method for
  {DNS} and implicit {LES} of turbulent flows, Eur. J. Mech. B-Fluid 55, Part 2
  (2016) 367 -- 379.
\newblock \href {http://dx.doi.org/10.1016/j.euromechflu.2015.08.010}
  {\path{doi:10.1016/j.euromechflu.2015.08.010}}.

\bibitem{Beck14}
A.~D. Beck, T.~Bolemann, D.~Flad, H.~Frank, G.~J. Gassner, F.~Hindenlang, C.-D.
  Munz, High-order discontinuous {G}alerkin spectral element methods for
  transitional and turbulent flow simulations, Int. J. Numer. Meth. Fluids
  76~(8) (2014) 522--548.
\newblock \href {http://dx.doi.org/10.1002/fld.3943}
  {\path{doi:10.1002/fld.3943}}.

\bibitem{Wiart15}
C.~C. Wiart, K.~Hillewaert, L.~Bricteux, G.~Winckelmans, Implicit {LES} of free
  and wall-bounded turbulent flows based on the discontinuous
  {G}alerkin/symmetric interior penalty method, Int. J. Numer. Meth. Fluids
  78~(6) (2015) 335--354.
\newblock \href {http://dx.doi.org/10.1002/fld.4021}
  {\path{doi:10.1002/fld.4021}}.

\bibitem{Landmann08}
B.~Landmann, M.~Kessler, S.~Wagner, E.~Kr{\"a}mer, A parallel, high-order
  discontinuous {G}alerkin code for laminar and turbulent flows, Comput. Fluids
  37~(4) (2008) 427 -- 438.
\newblock \href {http://dx.doi.org/10.1016/j.compfluid.2007.02.009}
  {\path{doi:10.1016/j.compfluid.2007.02.009}}.

\bibitem{Wang15}
L.~Wang, W.~Kyle~Anderson, T.~Erwin, S.~Kapadia, High-order discontinuous
  {G}alerkin method for computation of turbulent flows, AIAA J. 53~(5) (2015)
  1159--1171.
\newblock \href {http://dx.doi.org/10.2514/1.J053134}
  {\path{doi:10.2514/1.J053134}}.

\bibitem{Bassi14}
F.~Bassi, A.~Ghidoni, A.~Perbellini, S.~Rebay, A.~Crivellini, N.~Franchina,
  M.~Savini, A high-order discontinuous {G}alerkin solver for the
  incompressible {RANS} and $k-\omega$ turbulence model equations, Comput.
  Fluids 98 (2014) 54 -- 68.
\newblock \href {http://dx.doi.org/10.1016/j.compfluid.2014.02.028}
  {\path{doi:10.1016/j.compfluid.2014.02.028}}.

\bibitem{Wiart15b}
C.~C. de~Wiart, K.~Hillewaert, E.~Lorriaux, G.~Verheylewegen, Development of a
  discontinuous galerkin solver for high quality wall-resolved/modelled dns and
  les of practical turbomachinery flows on fully unstructured meshes, ASME
  GT2015-43428, Montreal, Canada.

\bibitem{Hartmann16}
R.~Hartmann, H.~McMorris, T.~Leicht, Curved grid generation and {DG}
  computation for the {DLR-F11} high lift configuration, in: M.~Papadrakakis,
  V.~Papadopoulos, G.~Stefanou, V.~Plevris (Eds.), Proceedings of the ECCOMAS
  Congress 2016, Crete Island, Greece, 5--10 June 2016, 2016.

\bibitem{Giraldo02}
F.~X. Giraldo, J.~S. Hesthaven, T.~Warburton, Nodal high-order discontinuous
  {G}alerkin methods for the spherical shallow water equations, J. Comput.
  Phys. 181~(2) (2002) 499 -- 525.
\newblock \href {http://dx.doi.org/10.1006/jcph.2002.7139}
  {\path{doi:10.1006/jcph.2002.7139}}.

\bibitem{Giraldo08}
F.~X. Giraldo, M.~Restelli, A study of spectral element and discontinuous
  {G}alerkin methods for the {N}avier-–{S}tokes equations in nonhydrostatic
  mesoscale atmospheric modeling: equation sets and test cases, J. Comput.
  Phys. 227~(8) (2008) 3849 -- 3877.
\newblock \href {http://dx.doi.org/10.1016/j.jcp.2007.12.009}
  {\path{doi:10.1016/j.jcp.2007.12.009}}.

\bibitem{Wiart14}
C.~Carton~de Wiart, K.~Hillewaert, M.~Duponcheel, G.~Winckelmans, Assessment of
  a discontinuous {G}alerkin method for the simulation of vortical flows at
  high reynolds number, Int. J. Numer. Meth. Fluids 74~(7) (2014) 469--493.
\newblock \href {http://dx.doi.org/10.1002/fld.3859}
  {\path{doi:10.1002/fld.3859}}.

\bibitem{Collis02}
S.~S. Collis, Discontinuous {G}alerkin methods for turbulence simulation, in:
  Proceedings of the 2002 Center for Turbulence Research Summer Program, 2002,
  p. 115–167.

\bibitem{Marek15}
M.~Marek, A.~Tyliszczak, A.~Bogus\l{}awski, Large eddy simulation of
  incompressible free round jet with discontinuous {G}alerkin method, Int. J.
  Numer. Meth. Fluids 79~(4) (2015) 164--182, fld.4043.
\newblock \href {http://dx.doi.org/10.1002/fld.4043}
  {\path{doi:10.1002/fld.4043}}.

\bibitem{Noventa16}
G.~Noventa, F.~Massa, F.~Bassi, A.~Colombo, N.~Franchina, A.~Ghidoni, A
  high-order discontinuous {G}alerkin solver for unsteady incompressible
  turbulent flows, Comput. Fluids\href
  {http://dx.doi.org/10.1016/j.compfluid.2016.03.007}
  {\path{doi:10.1016/j.compfluid.2016.03.007}}.

\bibitem{Crivellini13}
A.~Crivellini, V.~D'Alessandro, F.~Bassi, High-order discontinuous {G}alerkin
  solutions of three-dimensional incompressible {RANS} equations, Comput.
  Fluids 81 (2013) 122--133.
\newblock \href {http://dx.doi.org/10.1016/j.compfluid.2013.04.016}
  {\path{doi:10.1016/j.compfluid.2013.04.016}}.

\bibitem{Ferrer12b}
E.~Ferrer, R.~H.~J. Willden, A high order discontinuous {G}alerkin –
  {F}ourier incompressible {3D} {N}avier–-{S}tokes solver with rotating
  sliding meshes, J. Comput. Phys. 231~(21) (2012) 7037 -- 7056.
\newblock \href {http://dx.doi.org/10.1016/j.jcp.2012.04.039}
  {\path{doi:10.1016/j.jcp.2012.04.039}}.

\bibitem{Tavelli16}
M.~Tavelli, M.~Dumbser, A staggered space--time discontinuous {G}alerkin method
  for the three-dimensional incompressible {N}avier--{S}tokes equations on
  unstructured tetrahedral meshes, J. Comput. Phys. 319 (2016) 294--323.
\newblock \href {http://dx.doi.org/10.1016/j.jcp.2016.05.009}
  {\path{doi:10.1016/j.jcp.2016.05.009}}.

\bibitem{Cockburn05}
B.~Cockburn, G.~Kanschat, D.~Sch{\"o}tzau, A locally conservative {LDG} method
  for the incompressible {N}avier--{S}tokes equations, Math. Comp. 74~(251)
  (2005) 1067--1095.
\newblock \href {http://dx.doi.org/10.1090/S0025-5718-04-01718-1}
  {\path{doi:10.1090/S0025-5718-04-01718-1}}.

\bibitem{Cockburn09}
B.~Cockburn, G.~Kanschat, D.~Sch{\"o}tzau, An equal-order {DG} method for the
  incompressible {N}avier--{S}tokes equations, J. Sci. Comput. 40~(1-3) (2009)
  188--210.
\newblock \href {http://dx.doi.org/10.1007/s10915-008-9261-1}
  {\path{doi:10.1007/s10915-008-9261-1}}.

\bibitem{Schoetzau02}
D.~Sch{\"o}tzau, C.~Schwab, A.~Toselli, Mixed hp-{DGFEM} for incompressible
  flows, SIAM J. Numer. Anal. 40~(6) (2002) 2171--2194.
\newblock \href {http://dx.doi.org/10.1137/S0036142901399124}
  {\path{doi:10.1137/S0036142901399124}}.

\bibitem{Klein15}
B.~Klein, F.~Kummer, M.~Keil, M.~Oberlack, An extension of the {SIMPLE} based
  discontinuous {G}alerkin solver to unsteady incompressible flows, Int. J.
  Numer. Meth. Fluids 77~(10) (2015) 571--589, fld.3994.
\newblock \href {http://dx.doi.org/10.1002/fld.3994}
  {\path{doi:10.1002/fld.3994}}.

\bibitem{Rhebergen13}
S.~Rhebergen, B.~Cockburn, J.~J. Van Der~Vegt, A space--time discontinuous
  {G}alerkin method for the incompressible {N}avier--{S}tokes equations, J.
  Comput. Phys. 233 (2013) 339--358.
\newblock \href {http://dx.doi.org/10.1016/j.jcp.2012.08.052}
  {\path{doi:10.1016/j.jcp.2012.08.052}}.

\bibitem{Guermond06}
J.~L. Guermond, P.~Minev, J.~Shen, An overview of projection methods for
  incompressible flows, Comput. Methods in Appl. Mech. Eng. 195~(44–47)
  (2006) 6011 -- 6045.
\newblock \href {http://dx.doi.org/10.1016/j.cma.2005.10.010}
  {\path{doi:10.1016/j.cma.2005.10.010}}.

\bibitem{Botti11}
L.~Botti, D.~A.~D. Pietro, A pressure-correction scheme for
  convection-dominated incompressible flows with discontinuous velocity and
  continuous pressure, J. Comput. Phys. 230~(3) (2011) 572 -- 585.
\newblock \href {http://dx.doi.org/10.1016/j.jcp.2010.10.004}
  {\path{doi:10.1016/j.jcp.2010.10.004}}.

\bibitem{Shahbazi07}
K.~Shahbazi, P.~F. Fischer, C.~R. Ethier, A high-order discontinuous {G}alerkin
  method for the unsteady incompressible {N}avier--{S}tokes equations, J.
  Comput. Phys. 222~(1) (2007) 391 -- 407.
\newblock \href {http://dx.doi.org/10.1016/j.jcp.2006.07.029}
  {\path{doi:10.1016/j.jcp.2006.07.029}}.

\bibitem{Lehrenfeld16}
C.~Lehrenfeld, J.~Sch{\"o}berl, High order exactly divergence-free hybrid
  discontinuous {G}alerkin methods for unsteady incompressible flows, Comput.
  Methods in Appl. Mech. Eng. 307 (2016) 339 -- 361.
\newblock \href {http://dx.doi.org/10.1016/j.cma.2016.04.025}
  {\path{doi:10.1016/j.cma.2016.04.025}}.

\bibitem{Karniadakis91}
G.~E. Karniadakis, M.~Israeli, S.~A. Orszag, High-order splitting methods for
  the incompressible {N}avier--{S}tokes equations, J. Comput. Phys. 97~(2)
  (1991) 414 -- 443.
\newblock \href {http://dx.doi.org/10.1016/0021-9991(91)90007-8}
  {\path{doi:10.1016/0021-9991(91)90007-8}}.

\bibitem{Hesthaven07}
J.~S. Hesthaven, T.~Warburton, Nodal discontinuous {G}alerkin methods:
  algorithms, analysis, and applications, Springer, 2007.
\newblock \href {http://dx.doi.org/10.1007/978-0-387-72067-8}
  {\path{doi:10.1007/978-0-387-72067-8}}.

\bibitem{Ferrer12}
E.~Ferrer, A high order discontinuous {G}alerkin-{F}ourier incompressible {3D}
  {N}avier--{S}tokes solver with rotating sliding meshes for simulating
  cross-flow turbines, Ph.D. thesis, University of Oxford (2012).

\bibitem{Steinmoeller13}
D.~T. Steinmoeller, M.~Stastna, K.~G. Lamb, A short note on the discontinuous
  {G}alerkin discretization of the pressure projection operator in
  incompressible flow, J. Comput. Phys. 251 (2013) 480 -- 486.
\newblock \href {http://dx.doi.org/10.1016/j.jcp.2013.05.036}
  {\path{doi:10.1016/j.jcp.2013.05.036}}.

\bibitem{Ferrer14}
E.~Ferrer, D.~Moxey, R.~H.~J. Willden, S.~J. Sherwin, Stability of projection
  methods for incompressible flows using high order pressure-velocity pairs of
  same degree: Continuous and discontinuous {G}alerkin formulations, Commun.
  Comput. Phys. 16 (2014) 817--840.
\newblock \href {http://dx.doi.org/10.4208/cicp.290114.170414a}
  {\path{doi:10.4208/cicp.290114.170414a}}.

\bibitem{EscobarVargas14}
J.~A. Escobar-Vargas, P.~J. Diamessis, T.~Sakai, A spectral quadrilateral
  multidomain penalty method model for high {R}eynolds number incompressible
  stratified flows, Int. J. Numer. Meth. Fluids 75~(6) (2014) 403--425.
\newblock \href {http://dx.doi.org/10.1002/fld.3899}
  {\path{doi:10.1002/fld.3899}}.

\bibitem{Joshi16}
S.~M. Joshi, P.~J. Diamessis, D.~T. Steinmoeller, M.~Stastna, G.~N. Thomsen, A
  post-processing technique for stabilizing the discontinuous pressure
  projection operator in marginally-resolved incompressible inviscid flow,
  Comput. Fluids\href {http://dx.doi.org/10.1016/j.compfluid.2016.04.021}
  {\path{doi:10.1016/j.compfluid.2016.04.021}}.

\bibitem{Karniadakis13}
G.~E. Karniadakis, S.~J. Sherwin, Spectral/hp element methods for computational
  fluid dynamics, Oxford University Press, 2013.
\newblock \href {http://dx.doi.org/10.1093/acprof:oso/9780198528692.001.0001}
  {\path{doi:10.1093/acprof:oso/9780198528692.001.0001}}.

\bibitem{Olshanskii09}
M.~Olshanskii, G.~Lube, T.~Heister, J.~L{\"o}we, Grad--div stabilization and
  subgrid pressure models for the incompressible {N}avier--{S}tokes equations,
  Comput. Methods in Appl. Mech. Eng. 198~(49–52) (2009) 3975 -- 3988.
\newblock \href {http://dx.doi.org/10.1016/j.cma.2009.09.005}
  {\path{doi:10.1016/j.cma.2009.09.005}}.

\bibitem{Arnold82}
D.~N. Arnold, An interior penalty finite element method with discontinuous
  elements, SIAM J. Numer. Anal. 19~(4) (1982) 742--760.
\newblock \href {http://dx.doi.org/10.1137/0719052}
  {\path{doi:10.1137/0719052}}.

\bibitem{Kronbichler12}
M.~Kronbichler, K.~Kormann, A generic interface for parallel finite element
  operator application, Comput. Fluids 63 (2012) 135--147.
\newblock \href {http://dx.doi.org/10.1016/j.compfluid.2012.04.012}
  {\path{doi:10.1016/j.compfluid.2012.04.012}}.

\bibitem{Kormann16}
K.~Kormann, M.~Kronbichler, Efficient matrix-free implementations for
  discontinuous {G}alerkin methods, In preparation.

\bibitem{Ferrer11}
E.~Ferrer, R.~H.~J. Willden, A high order discontinuous {G}alerkin finite
  element solver for the incompressible {N}avier--{S}tokes equations, Comput.
  Fluids 46~(1) (2011) 224 -- 230.
\newblock \href {http://dx.doi.org/10.1016/j.compfluid.2010.10.018}
  {\path{doi:10.1016/j.compfluid.2010.10.018}}.

\bibitem{Orszag86}
S.~A. Orszag, M.~Israeli, M.~O. Deville, Boundary conditions for incompressible
  flows, J. Sci. Comput. 1~(1) (1986) 75--111.
\newblock \href {http://dx.doi.org/10.1007/BF01061454}
  {\path{doi:10.1007/BF01061454}}.

\bibitem{leriche2000high}
E.~Leriche, G.~Labrosse, High-order direct {S}tokes solvers with or without
  temporal splitting: numerical investigations of their comparative properties,
  SIAM J. Sci. Comput. 22~(4) (2000) 1386--1410.
\newblock \href {http://dx.doi.org/10.1137/S1064827598349641}
  {\path{doi:10.1137/S1064827598349641}}.

\bibitem{Klein13}
B.~Klein, F.~Kummer, M.~Oberlack, A {SIMPLE} based discontinuous {G}alerkin
  solver for steady incompressible flows, J. Comput. Phys. 237 (2013) 235 --
  250.
\newblock \href {http://dx.doi.org/10.1016/j.jcp.2012.11.051}
  {\path{doi:10.1016/j.jcp.2012.11.051}}.

\bibitem{Hillewaert13}
K.~Hillewaert, Development of the discontinuous {G}alerkin method for
  high-resolution, large scale {CFD} and acoustics in industrial geometries,
  Ph.D. thesis, Univ. de Louvain (2013).

\bibitem{Riviere08}
B.~Riviere, Discontinuous {G}alerkin methods for solving elliptic and parabolic
  equations: theory and implementation, Society for Industrial and Applied
  Mathematics, Philadelphia, PA, USA, 2008.

\bibitem{Hartmann07}
R.~Hartmann, Adjoint consistency analysis of discontinuous {G}alerkin
  discretizations, SIAM J. Numer. Anal. 45~(6) (2007) 2671--2696.
\newblock \href {http://dx.doi.org/10.1137/060665117}
  {\path{doi:10.1137/060665117}}.

\bibitem{leriche2006numerical}
E.~Leriche, E.~Perchat, G.~Labrosse, M.~O. Deville, Numerical evaluation of the
  accuracy and stability properties of high-order direct {S}tokes solvers with
  or without temporal splitting, J. Sci. Comput. 26~(1) (2006) 25--43.
\newblock \href {http://dx.doi.org/10.1007/s10915-004-4798-0}
  {\path{doi:10.1007/s10915-004-4798-0}}.

\bibitem{Shahbazi05}
K.~Shahbazi, An explicit expression for the penalty parameter of the interior
  penalty method, J. Comput. Phys. 205~(2) (2005) 401 -- 407.
\newblock \href {http://dx.doi.org/10.1016/j.jcp.2004.11.017}
  {\path{doi:10.1016/j.jcp.2004.11.017}}.

\bibitem{Hedwig15}
M.~Hedwig, P.~W. Schr{\"o}der, A grad-div stabilized discontinuous {G}alerkin
  based thermal optimization of sorption processes via phase change materials,
  Tech. Rep.~5, Institut f{\"u}r Numerische und Angewandte Mathematik,
  Georg-August-Universit{\"a}t G{\"o}ttigen (2015).

\bibitem{Bangerth16}
W.~Bangerth, T.~Heister, L.~Heltai, G.~Kanschat, M.~Kronbichler, M.~Maier,
  B.~Turcksin, The deal.{II} library, version 8.3, Archive of Numerical
  Software 4~(100).
\newblock \href {http://dx.doi.org/10.11588/ans.2016.100.23122}
  {\path{doi:10.11588/ans.2016.100.23122}}.

\bibitem{Kopriva09}
D.~A. Kopriva, Implementing spectral methods for partial differential
  equations: algorithms for scientists and engineers, Springer, 2009.

\bibitem{Bastian14}
P.~Bastian, C.~Engwer, D.~G\"oddeke, O.~Iliev, O.~Ippisch, M.~Ohlberger,
  S.~Turek, J.~Fahlke, S.~Kaulmann, S.~M\"uthing, D.~Ribbrock, {EXA-DUNE}:
  Flexible {PDE} solvers, numerical methods and applications, in: L.~Lopes,
  J.~\v{Z}ilinskas, A.~Costan, R.~G. Cascella, G.~Kecskemeti, E.~Jeannot,
  M.~Cannataro, L.~Ricci, S.~Benkner, S.~Petit, V.~Scarano, J.~Gracia,
  S.~Hunold, S.~L. Scott, S.~Lankes, C.~Lengauer, J.~Carretero, J.~Breitbart,
  M.~Alexander (Eds.), Euro-Par 2014: Parallel Processing Workshops, Lecture
  Notes in Computer Science, Springer, 2014.
\newblock \href {http://dx.doi.org/10.1007/978-3-319-14313-2_45}
  {\path{doi:10.1007/978-3-319-14313-2_45}}.

\bibitem{MayBrown14}
D.~A. May, J.~Brown, L.~Le~Pourhiet, p{T}atin3{D}: High-performance methods for
  long-term lithospheric dynamics, in: J.~M. Kunkel, T.~Ludwig, H.~W. Meuer
  (Eds.), Supercomputing (SC14), New Orleans, 2014, pp. 1--11.
\newblock \href {http://dx.doi.org/10.1109/SC.2014.28}
  {\path{doi:10.1109/SC.2014.28}}.

\bibitem{Kronbichler15}
M.~Kronbichler, S.~Schoeder, C.~M{\"u}ller, W.~A. Wall, Comparison of implicit
  and explicit hybridizable discontinuous {G}alerkin methods for the acoustic
  wave equation, Int. J Numer. Meth. Eng. 106~(9) (2016) 712--739, nme.5137.
\newblock \href {http://dx.doi.org/10.1002/nme.5137}
  {\path{doi:10.1002/nme.5137}}.

\bibitem{Adams03}
M.~Adams, M.~Brezina, J.~Hu, R.~Tuminaro, Parallel multigrid smoothing:
  polynomial versus {G}auss--{S}eidel, J. Comput. Phys. 188 (2003) 593--610.
\newblock \href {http://dx.doi.org/10.1016/S0021-9991(03)00194-3}
  {\path{doi:10.1016/S0021-9991(03)00194-3}}.

\bibitem{Sundar15}
H.~Sundar, G.~Stadler, G.~Biros, Comparison of multigrid algorithms for
  high-order continuous finite element discretizations, Numer. Linear Algebra
  Appl. 22 (2015) 664--680.
\newblock \href {http://dx.doi.org/10.1002/nla.1979}
  {\path{doi:10.1002/nla.1979}}.

\bibitem{Varga09}
R.~S. Varga, Matrix iterative analysis, 2nd Edition, Springer, Berlin, 2009.

\bibitem{Burstedde11}
C.~Burstedde, L.~C. Wilcox, O.~Ghattas, \texttt{p4est}: {S}calable algorithms
  for parallel adaptive mesh refinement on forests of octrees, SIAM J. Sci.
  Comput. 33~(3) (2011) 1103--1133.
\newblock \href {http://dx.doi.org/10.1137/100791634}
  {\path{doi:10.1137/100791634}}.

\bibitem{Bangerth11}
W.~Bangerth, C.~Burstedde, T.~Heister, M.~Kronbichler, Algorithms and data
  structures for massively parallel generic finite element codes, ACM Trans.
  Math. Softw. 38~(2).
\newblock \href {http://dx.doi.org/10.1145/2049673.2049678}
  {\path{doi:10.1145/2049673.2049678}}.

\bibitem{guermond2003velocity}
J.-L. Guermond, J.~Shen, Velocity-correction projection methods for
  incompressible flows, SIAM J. Numer. Anal. 41~(1) (2003) 112--134.
\newblock \href {http://dx.doi.org/10.1137/S0036142901395400}
  {\path{doi:10.1137/S0036142901395400}}.

\bibitem{John04}
V.~John, Reference values for drag and lift of a two-dimensional time-dependent
  flow around a cylinder, Int. J. Numer. Meth. Fluids 44~(7) (2004) 777--788.
\newblock \href {http://dx.doi.org/10.1002/fld.679}
  {\path{doi:10.1002/fld.679}}.

\bibitem{Schaefer96}
M.~Sch{\"a}fer, S.~Turek, F.~Durst, E.~Krause, R.~Rannacher, Benchmark
  computations of laminar flow around a cylinder, Springer, 1996.

\bibitem{Moin98}
P.~Moin, K.~Mahesh, Direct numerical simulation: a tool in turbulence research,
  Annu. Rev. Fluid Mech. 30~(1) (1998) 539--578.
\newblock \href {http://dx.doi.org/10.1146/annurev.fluid.30.1.539}
  {\path{doi:10.1146/annurev.fluid.30.1.539}}.

\bibitem{Moser99}
R.~D. Moser, J.~Kim, N.~N. Mansour, Direct numerical simulation of turbulent
  channel flow up to ${R}e_{\tau}=590$, Phys. Fluids 11~(4) (1999) 943--945.
\newblock \href {http://dx.doi.org/10.1063/1.869966}
  {\path{doi:10.1063/1.869966}}.

\bibitem{Vreman14}
A.~W. Vreman, J.~G.~M. Kuerten, Comparison of direct numerical simulation
  databases of turbulent channel flow at ${R}e_{\tau} = 180$, Phys. Fluids
  26~(1).
\newblock \href {http://dx.doi.org/10.1063/1.4861064}
  {\path{doi:10.1063/1.4861064}}.

\bibitem{Case73}
K.~M. Case, F.~J. Dyson, E.~A. Frieman, C.~E. Grosch, F.~W. Perkins, Numerical
  simulation of turbulence, Tech. Rep. AD-774 161, Stanford Research Institute
  (1973).

\bibitem{Chapman79}
D.~R. Chapman, Computational aerodynamics development and outlook, AIAA J.
  17~(12) (1979) 1293--1313.

\bibitem{Grinstein07}
F.~F. Grinstein, L.~G. Margolin, W.~J. Rider, Implicit large eddy simulation:
  computing turbulent fluid dynamics, Cambridge University Press, 2007.

\bibitem{Hughes95}
T.~J.~R. Hughes, Multiscale phenomena: {G}reen's functions, the
  {D}irichlet-to-{N}eumann formulation, subgrid scale models, bubbles and the
  origins of stabilized methods, Comput. Methods in Appl. Mech. Eng. 127~(1)
  (1995) 387--401.
\newblock \href {http://dx.doi.org/10.1016/0045-7825(95)00844-9}
  {\path{doi:10.1016/0045-7825(95)00844-9}}.

\end{thebibliography}

\end{document}